\documentclass[a4paper,reqno,11pt,UTF8]{amsart}
\usepackage{mathrsfs}
\usepackage{txfonts}
\usepackage{amssymb}
\usepackage{amsfonts,amsmath,amsthm,amssymb,stmaryrd,esint,extarrows,chemarrow}
\usepackage[numbers,sort&compress]{natbib}
\usepackage{booktabs,tabularx}
\usepackage{soul, color, xcolor}
\usepackage{subfigure}
\usepackage{float}
\usepackage{mathdots}
\usepackage{graphicx}
\usepackage{hyperref}

\newtheorem{Theorem}{Theorem}[section]
\newtheorem{Lemma}{Lemma}[section]
\newtheorem{Proposition}{Proposition}[section]
\newtheorem{Remark}{Remark}[section]

\newtheorem{Definition}{Definition}[section]

\numberwithin{equation}{section}
\usepackage[left=1 in, right=1 in,top=1 in, bottom=1 in]{geometry}

\def\XXint#1#2#3{{\setbox0=\hbox{$#1{#2#3}{\int}$ }
\vcenter{\hbox{$#2#3$ }}\kern-.6\wd0}}

\DeclareMathOperator{\diver}{div}
\DeclareMathOperator{\cur}{curl}
\def\r3{\mathbb{R}^3}

\DeclareMathOperator{\gl}{\mathcal{G}}
\newcommand{\na}{\nabla}
\newcommand{\la}{\lambda}

\newcommand{\De}{\Delta}

\newcommand{\pa}{\partial}
\newcommand{\ga}{\gamma}

\def\leq{\leqslant}
\def\geq{\geqslant}

\let\sss= \scriptscriptstyle

\makeatletter
\@namedef{subjclassname@2020}{\textup{2020} Mathematics Subject Classification}
\makeatother

\allowdisplaybreaks
\begin{document}
\bibliographystyle{plain}

\title[I\MakeLowercase{ncompressible} N\MakeLowercase{on-isothermal} V\MakeLowercase{iscoelastic} F\MakeLowercase{lows}]{G\MakeLowercase{lobal weak solutions of incompressible viscoelastic heat-conducting flows}}


\author{Yong Wang}
\address{Yong Wang \newline School of Mathematical Sciences, Key Laboratory of Mathematical Foundation and Computation in Nonlinear Sciences of Guangdong Higher Education Institutes, and Guangdong Key Laboratory of Applied Mathematics and Artificial Intelligence, South China Normal University, Guangzhou, Guangdong 510631, China}
\email{wangyongxmu@163.com}

\author{Zhouping Xin}
\address{Zhouping Xin \newline The Institute of Mathematical Sciences, The Chinese University of Hong Kong, Shatin, N.T., Hong Kong, China}
\email{zpxin@ims.cuhk.edu.hk}

\author{Yanting Yuan$^\ast$}
\address{Yanting Yuan \newline School of Mathematical Sciences, Key Laboratory of Mathematical Foundation and Computation in Nonlinear Sciences of Guangdong Higher Education Institutes, and Guangdong Key Laboratory of Applied Mathematics and Artificial Intelligence, South China Normal University, Guangzhou, Guangdong 510631, China}
\email{ytyuan2404@163.com}

\thanks{$^\ast$Corresponding author: Yanting Yuan (ytyuan2404@163.com).}

\begin{abstract}
In this paper, we establish the global existence of weak solutions to the Cauchy problem for incompressible viscoelastic heat-conducting flows with discontinuous initial data in two and three dimensions. Compared with the isothermal case, the temperature effects bring essential difficulties, manifested in the energy equation by highly singular nonlinear terms that pose a significant obstacle to closing conventional energy estimates. To overcome these difficulties, we develop a refined a priori estimate framework of incompressible viscoelastic heat-conducting system. By introducing time-weighted estimates, employing refined analysis in Hardy spaces and Besov spaces, and carefully exploiting the structure of effective viscous fluxes together with the  energy equation, we successfully establish the global existence of weak solutions under small energy conditions on the initial data. As a byproduct, our analysis results generalize those for the ``intermediate weak" solutions of incompressible viscoelastic equations obtained by Hu-Lin [Comm. Pure Appl. Math. 69 (2) (2016)].
\end{abstract}


\maketitle

\setcounter{tocdepth}{3}

\pagenumbering{Roman}
\tableofcontents

\pagenumbering{arabic}

\section{Introduction}

\

In fluid dynamics research, heat-conducting effects play a crucial role in many practical flow problems. Temperature variations not only affect the thermodynamic properties of fluids but also significantly alter their viscoelastic response, particularly in complex fluids such as polymers (see Chapter 3 of \cite{Larson1999}) and glassy liquids (see Chapter 4 of \cite{Larson1999}). Therefore, establishing and studying viscoelastic heat-conducting flow models has significant physical meaning and practical value. The incompressible viscoelastic heat-conducting system (1.1) considered in this paper incorporates an energy equation with the heat-conduction term based on the classical viscoelastic model \cite{Lin-Liu-Zhang2005}, which allows for a more realistic description of energy transfer and dissipation during actual flow processes.

Specifically, we consider the following incompressible viscoelastic heat-conducting system:
\begin{align}\label{1.1'}
\begin{cases}
u_t+\diver(u\otimes u)=-\nabla P+\mu\Delta u+c^2\diver(\mathbb{F}{\mathbb{F}}^\top),\\
\mathbb{F}_t+{u}\cdot \nabla \mathbb{F}=\nabla { u}\mathbb{F},\\
\mathcal{E}_t+\diver(\mathcal{E} { u})=\diver(\mathbb{T}{ u})+\kappa\Delta\theta,\\
\diver {u}=0, \quad (x,t)\in \mathbb{R}^n\times\mathbb{R}^+,n=2,3,
\end{cases}
\end{align}
which is derived in the appendix.
Here ${u(x,t)}$, $\mathbb{F}(x,t)$, $P(x,t)$, $e(x,t)$ and $\theta(x, t)$ represent, respectively, the velocity, the deformation gradient, the pressure, the specific internal energy and the absolute temperature, and $\mathcal{E}=\tfrac{1}{2}|{u}|^2+ e+\frac12 c^2|\mathbb{F}|^2$ is the specific total energy. The positive coefficients $\mu$, $c$, and $\kappa$ denote the viscosity, the speed of propagation of shear waves and the heat conductivity, respectively. The total stress tensor $\mathbb{T}=-P\mathbb{I}+\mathbb{T}_v+\mathbb{T}_e$,
where $\mathbb{T}_v$ and $\mathbb{T}_e$ denote the viscous and elastic stress tensors, respectively, and $\mathbb{I}=(\delta_{i,j})_{i,j=1}^n$ is the identity matrix. For simplicity, we assume that $\mathbb{T}_v$ satisfies Newton's law of viscosity, that is,
\begin{align*}
\mathbb{T}_v=2\mu\mathbb{D}({u}),\quad \mathbb{D}({u})=\tfrac{1}{2}(\nabla {u}+(\nabla {u})^\top),
\end{align*}
and we choose
\begin{align*}
\mathbb{T}_e=c^2 \mathbb{F}\mathbb{F}^\top.
\end{align*}
In order to simplify the presentation, we shall consider the case
\begin{equation}\label{cv}
e=C_v\theta,
\end{equation}
where $C_v>0$ is the
specific heat at constant volume.

The system \eqref{1.1'} satisfies the first and second laws of thermodynamics, and the elastic stress tensor is determined via the entropy inequality, thereby ensuring the physical plausibility. Moreover, the model unifies the classical heat-conducting Navier-Stokes system with the incompressible isothermal viscoelastic flows in form. When the deformation gradient is disregarded, the system naturally reduces to the classical heat-conducting fluid equations; when the energy equation is omitted, the system simplifies to the incompressible viscoelastic equations (see \eqref{yi} below). These two aspects together demonstrate the system's natural extensibility in its structural design.

In this paper, we investigate the Cauchy problem of system \eqref{1.1'} with the initial data
\begin{equation}\label{1.2'}
(u(x, t), \mathbb{F}(x, t),\mathcal{E}(x, t))\mid_{t=0}=(u_{0}(x), \mathbb{F}_{0}(x),\mathcal{E}_0(x)).
\end{equation}
Here $\mathcal{E}_0(x)=\tfrac{1}{2}|{u_0}|^2+ C_v\theta_0+\frac12 c^2|\mathbb{F}_0|^2$.

We also assume that  $\mathbb F_0$ satisfies the following constraints:
\begin{equation}\label{1.3'}
\quad\det\mathbb F_0=1,\quad\mathrm{div}(\mathbb F_0^\top)=0,\quad  {\mathbb F}_0^{lk}\partial_{x_l}{\mathbb F}_0^{ij}= {\mathbb F}_0^{lj}\partial_{x_l}{\mathbb F}_0^{ik}.
\end{equation}
For the classical solutions $(u,\mathbb F,\theta)$ of the Cauchy problem \eqref{1.1'}--\eqref{1.3'}, it is well-known that
\begin{equation*}
\det\mathbb F=1,\quad\mathrm{div}({\mathbb F}^\top)=0,\quad  {\mathbb F}^{lk}\partial_{x_l}{\mathbb F}^{ij}= {\mathbb F}^{lj}\partial_{x_l}{\mathbb F}^{ik}
\end{equation*}
are valid for all $t>0$ (see   \cite{Chen-Zhang2006,Lin-Liu-Zhang2005,Lei-Liu-Zhou2008}).

In case that the deformation gradient $\mathbb F$ is absent, equations \eqref{1.1'} will be reduced to the classical Navier-Stokes equations for heat-conducting fluids.  For the Navier-Stokes equations in the incompressible, heat-conducting case, Lions  \cite{Lions1996} proved the global existence of weak solutions in the spirit of Leray. On the other hand, in the compressible case, Matsumura and Nishida \cite{Matsumura-Nishida1980} first proved the global existence of smooth solutions when the initial data are close to a nonvacuum equilibrium in $H^3(\mathbb R^n)$.  It is worth mentioning that Huang-Li-Xin \cite{Huang-Li-Xin2012} and Li-Xin \cite{Huang-Li2018} have refined their results by attaining the global existence of classical solutions with smooth initial data of small energy but possibly large oscillations and vacuum.
In the case of large initial data, the global existence of finite energy weak solutions was established by Lions \cite{Lions1998} when the adiabatic exponent $\ga\geq \frac{3n}{n+2}(n=2,3)$. Later, Feireisl et al. \cite{Feireisl-Novotny-Petzeltov2001} extended the global existence to the case of $\ga>\frac n2$. Furthermore, Jiang-Zhang \cite{Jiang-Zhang2001,Jiang-Zhang2003} extended it to the case that $\ga>1$ but only for the spherically symmetric or axisymmetric solutions.

Then, Hoff \cite{Hoff1995_1} first proposed a class of ``intermediate weak" solutions. These solutions are defined as follows: they are small in the weak $L^2$ norm, have nonnegative and essentially bounded initial densities, but allow for discontinuities. The rationality of ``intermediate weak" solutions lies in their ability to bridge the gap between the small-smooth solutions obtained by Matsumura and Nishida \cite{Matsumura-Nishida1980} and large-energy weak solutions obtained by Lions \cite{Lions1998}, especially in the context of discontinuous initial data. Compared to the initial data of small-smooth solutions, the initial data of ``intermediate weak" solutions allows discontinuous across a hypersurface of $\mathbb R^n$. Meanwhile, compared to large-energy weak solutions, which possess very little regularity, ``intermediate weak" solutions retain important qualitative features such as the smoothness of vorticity and the continuity of the effective viscous flux $F = (\mu + \lambda) \text{div} \, u - P(\rho) + P(\bar{\rho})$. The key insight lies in decomposing the velocity gradient into smoother components. This is achieved particularly through the use of the effective viscous flux, which compensates for the lack of regularity in the density. This allows one to close the energy estimates and obtain pointwise control on the density, which would otherwise be inaccessible for large-energy weak solutions. Moreover, the solution exhibits instantaneous regularity for positive times: the velocity becomes H\"{o}lder continuous, the vorticity and the effective viscous flux gain $H^1$ regularity, and discontinuities are convected along particle paths. Thus, ``intermediate weak" solutions are not merely mathematical artifacts, and they reflect a physically consistent and mathematically rigorous framework for handling compressible flows with discontinuous initial data.

Hoff \cite{Hoff1995_1} proved the global existence of weak solutions of the Navier-Stokes equation for compressible, isothermal flow ($P(\rho)= a\rho$) in two and three dimensions:
\begin{align*}
\begin{cases}
\rho_t+\diver(\rho u)=0,\\
(\rho u)_t+\diver(\rho u\otimes u)=-\nabla P(\rho)+\mu\Delta u+(\mu+\lambda)\na\diver u,
\end{cases}
\end{align*}
under the constraints on the initial data and viscosity coefficients that
\begin{itemize}
\item the initial density $\rho_0$ is close to a positive constant in $L^2(\mathbb R^n)\cap L^\infty(\mathbb R^n)$;
\item the initial velocity $u_0$ is small in $L^2(\mathbb R^n)$ and bounded in $L^{2^n}(\mathbb R^n)$;
\item $\mu>0$, $\mu+\lambda\geq0$, $\frac{\mu+\lambda}{\mu}<b_n$, where $b_n=1+\sqrt{2}$ when $n=2$ and $b_n=\frac{1+\sqrt{3}}{\sqrt{3}}$ when $n=3$.
\end{itemize}
Then, Hoff \cite{Hoff1995_2} extended the above existence theory to general polytropic pressures $P(\rho)= a\rho^\gamma, \gamma> 1$.
Later, Hoff \cite{Hoff1997} established the global existence of weak solutions of the  Navier-Stokes equations for compressible, heat-conducting fluids :
\begin{align*}
\begin{cases}
\rho_t+\diver(\rho u)=0,\\
(\rho u)_t+\diver(\rho u\otimes u)=-\nabla P(\rho,e)+\mu\Delta u+(\mu+\lambda)\na\diver u,\\
(\rho \mathcal E_1)_t+\diver(\rho \mathcal E_1u+Pu)=\Delta(\kappa \theta+\frac12\mu |u|^2)+\mu\diver((\na u) u)+\lambda \diver((\diver u)u),
\end{cases}
\end{align*}
where $\mathcal E_1=\frac12|u|^2+e$, under the constraints on the initial data and viscosity  coefficients that
\begin{itemize}
\item the initial density $\rho_0$ is close to a positive constant in $L^2(\mathbb R^n)\cap L^\infty(\mathbb R^n)$;
\item the initial velocity $u_0$ is small in $H^s(\mathbb R^n)\cap L^4(\mathbb R^n)$, where $s=0$ when $n=2$ and $s\in(\frac13,\frac12)$ when $n=3$;
\item the initial temperature $\theta_0$ is close to a positive constant in $L^2(\mathbb R^n)$ and has a positive lower bound;
\item $\mu>0$, $0<\mu+\lambda<(1+\sqrt{2}){\mu}$ when $n=2$ and $\mu>0$, $\frac13\mu<\mu+\lambda<(1+\frac{\sqrt{15}}{3})\mu$ when $n=3$.
\end{itemize}

Recently, many articles have improved the conditions required for proving the existence of  ``intermediate weak" solutions  in \cite{Hoff1995_1,Hoff1997}. Those results are mainly manifested in the following aspects:
\begin{itemize}
\item the condition that the initial density $\rho_0$ must have a positive lower bound is improved in \cite{Hoff2005} to allow for initial vacuum.
\item the smallness of the initial density perturbations $\rho_0-\tilde{\rho}$ $(\tilde{\rho}>0)$ in $L^\infty(\mathbb R^n)$ is relaxed to the boundedness in $L^\infty(\mathbb R^n)$ (see \cite{Suen-Hoff2012,Suen2020,Wu-Tan2018});
\item the smallness of the initial velocity $u_0$ in $H^s(\mathbb R^n)$ is relaxed to the boundedness in $ H^s(\mathbb R^n)$ (cf. methods and techniques in \cite{Huang-Li-Xin2012});
\item many articles have partially relaxed the restrictions on viscosity coefficients, as can be seen in \cite{Hoff2005,Suen-Hoff2012,Suen2020,Wu-Tan2018, Hoff2002}. For instance in \cite{Suen2020}:
    \begin{align*}
    \mu>0,\ 0<\mu+\lambda<\frac{\mu}{4}\quad \mbox{in}\ \mathbb{R}^3.
    \end{align*}
\end{itemize}

What's more, Hoff's techniques were applied and  generalized to other systems, such as the incompressible viscoelastic system \cite{Hu-Lin2016}, the two-phase flow model \cite{Li-Liu-Ye2022}, the compressible magnetohydrodynamics \cite{Suen-Hoff2012} and the nematic liquid crystal flow \cite{Wu-Tan2018}. It is worth noticing that those ``intermediate weak" solutions would be more regular than the one developed by Lions \cite{Lions1998}, and the uniqueness and continuous dependence of solutions have been constructed by Hoff \cite{Hoff2006} and Suen \cite{Suen2020}.

The other equations closely related to (1.1) are the incompressible viscoelastic equations, which lack the energy equation:
\begin{align}\label{yi}
\begin{cases}
u_t+\diver(u\otimes u)=-\nabla P+\mu\Delta u+c^2\diver(\mathbb{F}{\mathbb{F}}^\top),\\
\mathbb{F}_t+{u}\cdot \nabla \mathbb{F}=\nabla { u}\mathbb{F},\\
\diver {u}=0, \quad (x,t)\in \mathbb{R}^n\times\mathbb{R}^+,n=2,3.
\end{cases}
\end{align}
A great number of works have been devoted to the well-posedness of solutions to incompressible viscoelastic equations. For classical solutions of \eqref{yi}, Lin et al. \cite{Lin-Liu-Zhang2005} proved the global existence of small solutions by introducing auxiliary vector fields instead of transport variables in two dimensions, while Lei and Zhou \cite{Lei-Zhou2005} obtained similar results by directly applying incompressible limits on deformation gradients. Corresponding results have been established in \cite{Lei-Liu-Zhou2008} and \cite{Chen-Zhang2006} for the three-dimensional case. For ``intermediate weak" solutions of \eqref{yi}, Hu and Lin \cite{Hu-Lin2016} proved the global existence only in two dimensions under the conditions that
\begin{itemize}
\item the initial deformation gradient $\mathbb F_0$ is close to the identity matrix in $L^2(\mathbb R^2)\cap L^\infty(\mathbb R^2)$;
\item  the initial velocity is small in $L^2(\mathbb R^2)$ and bounded in $L^p(\mathbb R^2)$ for some $p>2$.
\end{itemize}
For compressible viscoelastic folws, the global existence and asymptotic behaviors of strong solutions have been established, one can refer to \cite{Qian-Zhang2010,Wei-Lin-Yao2016,Hu-Wu2013,Hu-Wang2010,Hu-Ou-Wang-Lu2023} and the references cited therein. Besides, Hu proved the global existence of  weak solutions with discontinuous initial data in   \cite{Hu2018}. The global existence of weak solutions with large initial data for viscoelastic equations \eqref{yi} is still an open problem.

In this paper, we will prove the global existence of ``intermediate weak" solutions for the Cauchy problem \eqref{1.1'}--\eqref{1.3'}.  Before stating our main result, we would like to give some notations that will be used throughout this paper.\\

\noindent
\textbf{Notations:}

We use subscripts to denote differentiation with respect to $x$; thus $u_j^k={\partial u^k}/{\partial x_j}$, etc.
For $s\in\mathbb{R}$, the $H^s$ norm of a function $f$ is given by
\begin{equation*}
\|f\|_{H^s({\mathbb{R}^n})}^2\triangleq\int_{\mathbb{R}^n}\left(1+|\xi|^2\right)^s\lvert\mathcal F({f})(\xi)\rvert^2\;d\xi,
\end{equation*}
where $\mathcal F$ and $\mathcal F^{-1}$ denote the Fourier transform and the inverse Fourier transform, respectively. The symbol $\mathcal{H}^p(\mathbb{R}^n)$ with $0<p<+\infty$ denotes the Hardy space.
We adopt the following usual notations for H\"{o}lder seminorm: for $v: \mathbb{R}^n\rightarrow\mathbb{R}^m$, $\alpha \in (0,1]$,
\begin{equation*}
\langle v\rangle^\alpha\triangleq\sup_{\overset{x_1,\;x_2\in\mathbb{R}^n}{x_1\neq x_2}}\frac{|v(x_2)-v(x_1)|}{|x_2-x_1|^\alpha};
\end{equation*}
and for $v:Q\subseteq\mathbb{R}^n\times[0,\infty)\to\mathbb{R}^m$, $\alpha_1,  \alpha_2\in (0,1]$,
\begin{equation*}
\langle v\rangle_{Q}^{\alpha_{1},\alpha_{2}}\triangleq\sup_{\overset{(x_{1},t_{1}),\;(x_{2},t_{2})\in Q}{(x_{1},t_{1})\neq(x_{2},t_{2})}}
\frac{|v(x_{2},t_{2})-v(x_{1},t_{1})|}{|x_{2}-x_{1}|^{\alpha_{1}}+|t_{2}-t_{1}|^{\alpha_{2}}}.
\end{equation*}
For a given solution $(u,\mathbb F,\theta)$, we define $A(T)\triangleq A_1(T)+A_2(T)+A_3(T)+A_4(T)$, with
\begin{align}\label{A}
A_1(T)\triangleq&\sup\limits_{0\leq t\leq T}\int_{\mathbb{R}^{n}}(|u(x, t)|^2+|\mathbb{F}(x, t)-\mathbb I|^2+\sigma(t)|\nabla u(x, t)|^2+\sigma(t)^n|\mathcal{P}\dot{u}(x, t)|^2\big)\;dx\nonumber\\
&+\int_0^T\int_{\mathbb{R}^n}(|\nabla u(x, t)|^2+\sigma(t)|\mathcal{P}\dot {u}(x, t)|^2+\sigma(t)^n|\nabla\mathcal{P}\dot{u}(x, t)|^2)\;dxdt,\nonumber\\
A_2(T)\triangleq&\sup\limits_{0\leq t\leq \sigma(T)}\int_{\mathbb{R}^{n}}|u(x, t)|^4\;dx+\int_0^{\sigma(T)}\int_{\mathbb{R}^n}|u(x,t)|^{2}|\nabla u(x, t)|^2\;dxdt,\nonumber\\
A_3(T)\triangleq&\sup\limits_{0\leq t\leq T}\int_{\mathbb{R}^n}(\sigma(t)^{1-s}|\nabla u(x, t)|^2+\sigma(t)^{2-s}|\mathcal P\dot u(x, t)|^2)\;dx\nonumber\\
&+\int_0^T\int_{\mathbb{R}^n}(\sigma(t)^{1-s}|\mathcal{P}\dot{u}(x, t)|^2+\sigma(t)^{2-s}|\na\mathcal{P}\dot{u}(x, t)|^2)\;dxdt\nonumber,\\
A_4(T)\triangleq&\sup\limits_{0\leq t\leq T}\int_{\mathbb{R}^n}(|\theta(x,t)-\tilde{\theta}|^2+\sigma(t)^2|\na \theta(x,t)|^2+\sigma(t)^4|\dot \theta(x,t)|^2)\;dx\nonumber\\
&+\int_0^T\int_{\mathbb{R}^n}(|\na \theta(x,t)|^2+\sigma(t)^2 |{\dot\theta(x,t)}|^2+\sigma(t)^4 |{\na\dot\theta(x,t)}|^2)\;dxdt,
\end{align}
where $\dot{f}=\frac {Df}{Dt}\triangleq f_t+u\cdot\nabla f$ is the material derivative of $f$, $\sigma(t)=\operatorname{min}\{1, t\}$, and the operator $\mathcal{P}\triangleq \mathbb{I}+\nabla(-\Delta)^{-1}\operatorname{div}$ denotes the Helmholtz projection; and
\begin{equation}\label{B}
B(T)\triangleq\sup\limits_{0\leq t\leq T}\|\mathbb{F}(x, t)-\mathbb I\|_{\dot B_{p,1}^{\frac np}(\mathbb{R}^n)}^2.
\end{equation}

We introduce the following auxiliary variables associated with the system \eqref{1.1'}, which are known as the effective viscous fluxes $\mathcal{G}_1$ and $\mathcal{G}_2$:
\begin{gather*}
\mathcal{G}_1\triangleq\nabla u-
(-\Delta)^{-1}\nabla{\mathcal{P}}\operatorname{div}(\mathbb{F}\mathbb{F}^{\top}-\mathbb I),\\
\mathcal{G}_2\triangleq\cur u-
(-\Delta)^{-1}\cur\diver(\mathbb{F}\mathbb{F}^{\top}-\mathbb I),
\end{gather*}
where $\cur u\triangleq\na u-(\na u)^\top$. Then it is easy to check from $\eqref{1.1'}_1$ that
\begin{equation*}
\Delta \mathcal{G}_1=\na{\mathcal{P}}\dot{u}\quad \text{and} \quad \Delta \gl_2=\cur(\mathcal P\dot u).
\end{equation*}
Concerning the initial data $(u_0, \mathbb F_0, \theta_0)$, we assume that there is a positive number $N$, which may be arbitrarily large such that for $q>4$,
\begin{equation}\label{sn}
\|u_0\|_{H^s({\mathbb{R}^n})}+\|u_0\|_{L^q({\mathbb{R}^n})}\leq N,
\end{equation}
where $s\in(0,1)$ when $n=2$ and $s\in(\frac{8}{15},1)$ when $n=3$,
and we denote
\begin{equation}\label{ep}
\varepsilon_{0}\triangleq\|\mathbb F_{0}-\mathbb I\|_{{\dot B_{p,1}^{\frac np}(\mathbb{R}^n)}}^2+\|\mathbb F_{0}-\mathbb I\|_{L^2({\mathbb R}^n)}^2+
\|u_{0}\|_{L^2({\mathbb R}^n)}^2+\|\theta_0-\tilde{\theta}\|_{L^2({\mathbb R}^n)}^2,
\end{equation}
where $\tilde{\theta}$ is a fixed positive constant and $\dot B_{p,1}^{\frac np}({\mathbb R}^n)$ is the homogeneous Besov space with $p>2$ when $n=2$ and $p\in(3,6]$ when $n=3$.

A weak solution of \eqref{1.1'}--\eqref{1.3'} is defined in the same way as that in  \cite{Hoff1997} and \cite{Hu-Lin2016}.
\begin{Definition}
A triplet of functions $(u, \mathbb{F}, \mathcal{E})$ is a weak solution of the Cauchy problem \eqref{1.1'}--\eqref{1.3'}, if $u, \na u, \mathbb{F}, \theta, \na\theta \in L_{\mathrm{loc}}^1 (\mathbb{R}^n\times\mathbb{R}^+)$ and if for all test functions  $\varphi,{\psi},  {\phi} \in C_c^{\infty}(\mathbb{R}^n\times\mathbb{R}^{+})$ with $ \operatorname{div} \varphi=0$,
\begin{equation}
\int_{\mathbb{R}^n} {u}_{0} \cdot {\varphi}(x, 0) \;d x  +\int_{0}^{\infty} \int_{\mathbb{R}^n}[{u} \cdot {\varphi}_{t}+({u} \otimes {u}-c^2{\mathbb{F} \mathbb{F}}^{\top}): \nabla {\varphi}] \;d x d t
=\int_{0}^{\infty} \int_{\mathbb{R}^n} \mu\nabla {u}: \nabla {\varphi}\; d x d t,\label{we1}
\end{equation}
\begin{equation}
\int_{\mathbb{R}^n}\big(\mathbb{F}_{j}\big)_{0} \cdot {\psi}(x, 0) \;d x+\int_{0}^{\infty} \int_{\mathbb{R}^n}\big[\mathbb{F}_{j} \cdot {\psi}_{t}+\big({\mathbb{\mathbb{F}}}_{j} \otimes {u}-{u} \otimes {\mathbb{F}}_{j}\big): \nabla {\psi}\;\big]\, dx dt=0,\label{we2}
\end{equation}
where $\mathbb{F}_{j}, j=1,...,n$ are columns of \;$\mathbb{F}$, and
\begin{equation}\label{we3}
\int_{\mathbb{R}^n}\mathcal{E}_0 {\phi}(x, 0) \;d x+\int_{0}^{\infty} \int_{\mathbb{R}^n} (\mathcal{E}\phi_t+\mathcal{E}u\cdot\na\phi)\; dxdt=\int_{0}^{\infty} \int_{\mathbb{R}^n} (\mathbb T u \cdot\nabla\phi+\kappa  \na\theta\cdot\na\phi) \;dxdt.
\end{equation}
\end{Definition}

Now we state our main result on the existence of the global-in-time weak solution.
\begin{Theorem}\label{The1.1}
Let $n=2$ or $3$, fix the system parameters $\mu$, $c$, $\kappa$ and $C_v$ as in \eqref{1.1'}--\eqref{cv}, and choose constants $s$, $q$, $p$, $\tilde{\theta}$ and $N$ as in \eqref{sn}--\eqref{ep}. 
Then there exist positive constants $d$, $C$, $M$ and $\alpha$, depending only on the above parameters and constants, such that if the initial data $(u_0,\mathbb{F}_0,\mathcal{E}_0)$ satisfies $\varepsilon_{0}\leq d$,
then the Cauchy problem \eqref{1.1'}--\eqref{1.3'} has a global weak solution $(u,\mathbb{F},\mathcal{E})$  satisfying the following conditions:
\begin{gather}
A_1(T)+A_2(T)+A_4(T)+B(T)\leq C\varepsilon_{0}^\alpha,\quad A_3(T)\leq M,\quad \forall\ T>0,\label{ab}\\
u(\cdot,T),\ \theta(\cdot,T)-\tilde{\theta}\in H^1(\mathbb{R}^n),\quad \forall\ T>0,\label{H}\\
\|u\|_{L^\infty(\mathbb{R}^n\times[\tau,\infty))}+\langle u\rangle_{\mathbb{R}^{n} \times[\tau, \infty)}^{\beta, \frac{\beta}{n+2\beta}} \leq C(\tau,\beta) \varepsilon_{0}^\alpha, \quad \tau >0,\label{Ho}\\
\|\theta-\tilde{\theta}\|_{L^\infty(\mathbb{R}^n\times[\tau,\infty))}+\langle \theta\rangle_{\mathbb{R}^{n} \times[\tau, \infty)}^{\beta, \frac{\beta}{n+2\beta}}\leq C(\tau,\beta) \varepsilon_{0}^\alpha, \quad \tau >0,\label{Ho2}
\end{gather}
where $\beta\in (0,1)$ when $n=2$ and $\beta\in (0,\frac12]$ when $n=3$, and $C(\tau,\beta)$ may depend additionally on a positive lower bound for $\tau$ and $\beta$. 
Moreover, there is the following large-time behavior:
\begin{equation}
\lim_{t\to\infty}\big(\|(u(\cdot,t), \theta(\cdot,t)-\tilde{\theta})\|_{L^{p_1}(\mathbb{R}^n)}+\|\mathbb F(\cdot,t)-\mathbb I\|_{L^{p_2}(\mathbb{R}^n)}\big)=0\label{lb}
\end{equation}
holds for $p_1\in (2,\infty]$ and $p_2\in(2,\infty)$.
\end{Theorem}
\begin{Remark}
We shall prove the above theorem under the assumption that $\mu=c=C_v=\kappa=1$ 
for convenience. For a more general case not \eqref{cv} where $\theta$ is a smooth function of the internal energy $e\in(0,\infty)$, satisfying $\theta(e)/e\leq C$ and $\theta'(e)\geq C$  for a positive constant $C$, Theorem \ref{The1.1} is still valid.
\end{Remark}
\begin{Remark}
In the case that $u_0$ is small in $L^4$ (the additional  constraint for $u_0$ in \cite{Hoff1997}) instead of $u_0$ is bounded in $L^q(\mathbb R^n)$ for $q>4$, Theorem \ref{The1.1} also holds.
\end{Remark}

{\it Main ideas,  difficulties and strategies in analysis.} We now outline the main ideas in proving Theorem \ref{The1.1}, focusing on the difficulties that arise and the strategies we develop to overcome them. These strategies, while motivated by \cite{Hoff1997,Hu2018,Hu-Lin2016}, necessarily depart from those works in several respects, leading to weaker assumptions on initial data and wider ranges of parameters---precisely the improvements that constitute the main results of this paper.

First, we need to establish the a priori estimates for local smooth solutions of the Cauchy problem \eqref{1.1'}--\eqref{1.3'}, which consist of the energy estimates for $u$ (see Lemma \ref{A1}), the pointwise bounds (to be precise, the $\dot B_{p,1}^{n/p}$ norm) for $\mathbb F$ (see Lemma \ref{bt}), the energy estimates for $\theta$ (see Lemma \ref{theta}) and the pointwise bounds for $u$ and $\theta$ (see Lemma \ref{pt}). For $n=2$, the estimates for $u$ and $\mathbb F$ were already obtained by Hu-Lin \cite{Hu-Lin2016} for the isothermal case. However, extending these to $n=3$ and to the heat-conducting case  introduces new difficulties.

Consider first the estimate for $u$. In \cite{Hu-Lin2016}, the quantity $\iint\sigma|\dot u|^2dxdt$ appears naturally. Estimating it requires control of $\iint|u|^4dxdt$. For $n=2$, this follows from the standard energy estimate \eqref{nl} and the Gagliardo-Nirenberg inequality \eqref{GN}:
\begin{equation*}
\iint|u|^4dxdt\leq\int\left(\int |u|^2\;dx\right)\left(\int |\na u|^2\;dx\right)dt\leq \varepsilon^2.
\end{equation*} For $n=3$, the same inequality gives
\begin{equation}\label{121}
\iint|u|^4dxdt\leq\int\left(\int |u|^2\;dx\right)^\frac12\left(\int |\na u|^2\;dx\right)^\frac32dt.
\end{equation}
There seems to be no mechanism which enforces a bound on \eqref{121}, leading to the failure of estimating
$\iint\sigma|\dot u|^2dxdt$.  Fortunately, if we estimate $\iint\sigma|\mathcal P\dot u|^2dxdt$, the term $\iint|u|^4dxdt$ will not occur, thereby promoting the success of closing the estimates.

A second difficulty in closing the estimates for $A_1(T), A_2(T), B(T)$ is the lack of adequate smoothing rate in $L^2$ for $\mathcal P\dot u$ and $\na \mathcal P\dot u$ (see the definition of $A_3(T)$). To obtain the estimates of $A_3(T)$, we employ a Riesz-Thorin interpolation argument, which relies only on the boundedness (not smallness) of $u_0$ in $H^s(\mathbb R^n)$. This is another departure from \cite{Hoff1997,Hu2018}, where smallness in $H^s$ was required. As we will show, terms such as \eqref{388} and \eqref{269} reduce to products of $A_3(T)$ with $A_1(T)$ or $A_2(T)$, hence remain bounded by $C\varepsilon^\alpha$ even when $A_3(T)$ is merely bounded.

It is  worth pointing out that the effective viscous fluxes ${\mathcal G}_1$ and ${\mathcal G}_2$ play a crucial role in closing the estimates. For one thing, we need to control  the higher-order terms, for example, $\iint \sigma^n|\na u|^4dxdt$, owing to the convection term $u\cdot\na u$. Noticing that the Sobolev embedding $W^{2,2}\hookrightarrow W^{1,4}$ is invalid for $u$ in this paper since $\Delta u$ may not be locally integrable when $\na u$ is discontinuous across hypersurfaces of $\mathbb R^n$. Observing that $\na u$ is a linear combination of $\mathcal G_1$ and $(-\Delta)^{-1}\na\diver(\mathbb F\mathbb F^\top-\mathbb I)$, by virtue of their regularities, we ultimately obtain that $\iint \sigma^n|\na u|^4dxdt$ can be bounded by some terms appearing in the definitions of $A(T)$ and $B(T)$. For another, the effective viscous flux $\mathcal G_i$ provides a new way to derive the pointwise bound for $\mathbb F$ in the absence of the $L^\infty$ bound on $\na u$. To be specific, it holds that
\begin{align*}
\frac{D}{Dt}(\mathbb F-\mathbb I)+\mathbb F-\mathbb I=\mathcal{G}\mathbb F+\mathcal{R}\mathbb F-(\mathbb F-\mathbb I)(\mathbb F-\mathbb I),
\end{align*}
the derivation of which can be seen in \eqref{qd}. Since the Riesz transform is not convenient to work directly in $L^\infty(\mathbb R^n)$ for $\mathbb F$, instead we replace $L^\infty(\mathbb R^n)$ with the Besov space $\dot B_{p,1}^{n/p}(\mathbb R^n) \subset L^\infty(\mathbb R^n).$ Then the estimate for $B(T)$ can be computed carefully by substituting for $\mathcal G$ in terms of $\mathcal P\dot u$ in Lemma \ref{bt}.
Here we improve upon \cite{Hu2018}: when $n=2$, the smallness of initial velocity in $ H^s(\mathbb R^n)$ is relaxed to the boundedness in $H^s(\mathbb R^n)$, and the range of $s$ and $p$ is improved from $(0,\frac12)$, $(\frac{-1+\sqrt{1+16s}}{2s},4)$ to $(0,1)$, $(2,+\infty)$ (see Lemma \ref{bt}), respectively. Specifically, in this paper, we use the Gagliardo-Nirenberg inequality to estimate $\|\mathcal{G}_i\|_{\dot{B}_{p,1}^{n/p}(\mathbb R^n)}$ in \eqref{383} instead of the Hausdorff-Young's inequality in Lemma 4.1 of \cite{Hu2018}, such that we can relax the range of $s$ and $p$. What's more, in \eqref{388}, we equip the term $\|\na u\|^2_{L^2(\mathbb R^n)}$ with $t$ instead of $t^{1-s}$ (see Lemma 4.1 case two in \cite{Hu2018}). Since the term $\sup\limits_{t\in[0,1]}t\|\na u\|^2_{L^2(\mathbb R^n)}$ is controlled by $A_1(T)$, we can get desired estimate of  $\int_0^T\|\mathcal G_i\|_{{\dot{B}_{p,1}^{{n/p}}(\mathbb R^n)}}(t)\;dt$ by the estimate of $A_1(t)$ in \eqref{a12} and the boundedness of $A_3(T)$ in Lemma \ref{A3}. Hence, the smallness of initial velocity in $ H^s(\mathbb R^n)$ is  also not necessary in this process.

Turning to $\theta$, the energy estimates are more delicate. The basic energy estimate gives the desired bound for $u$ directly, but the analogous low-order estimate for $\theta$ requires a case analysis. For $T\leq 1$, since the term $|\na u|^2$ in $\eqref{1.1'}_3 $ is highly singular at the initial layer near $t=0$, we need to use the equations \eqref{1.1'} to replace it with other terms as follows:
\begin{equation*}
|\na u|^2=\mathcal{P}\diver((\mathbb{F}-\mathbb I)(\mathbb{F}-\mathbb I)^\top)u+\diver(\mathbb{F}-\mathbb I)u
+\Delta(\frac12|u|^2)-u\mathcal{P}\dot{u}.
\end{equation*}
After integration by parts, the integral about $\na u$ turns into $\int |u|^4dx$ and other routine terms (see \eqref{370} for example). Hence we have to derive the bounds for $\int |u|^4dx$, which together with Lemma \ref{2.3} give the reason that we assume that $u_0$ is additionally  bounded in $L^q(\mathbb R^n)(q>4)$. By H\"{o}lder's inequality,\begin{equation*}
\int_{\mathbb{R}^{n}}|u_0|^4\;dx\leq \left(\int_{\mathbb{R}^{n}}|u_0|^2\;dx\right)^\frac{q-4}{q-2}\left(\int_{\mathbb{R}^{n}}|u_0|^q\;dx\right)^\frac
{2}{q-2}\leq C\varepsilon^\alpha,
\end{equation*}
so the $L^4$ smallness assumed in \cite{Hoff1997} follows from the weaker $L^q$ boundedness. For $T>1$, we use a standard entropy argument to achieve our goal. The higher-order estimates of $\theta$ are similar to $u$. However, it should be noted here that in Lemma \ref{2.5}, careful computation gives that the higher-order term, $\iint \sigma^4 \Delta \dot \theta u\cdot \na\theta\;dxdt$, would be canceled, so that we can close the estimates. Then, the pointwise bounds for $u$ and $\theta$ are  similar to the ones in Hoff \cite{Hoff1997} (see Theorem 4.1 therein).

After establishing the a priori estimates, we need to construct  approximate solutions with mollified approximating initial data. It can be shown that these approximate solutions also satisfy the established a priori estimates.
Finally, we make full use of the a priori estimates to show that the approximate solutions converge to the weak solution of \eqref{1.1'}--\eqref{1.3'}.

The rest of this paper is devoted to proving Theorem \ref{The1.1}. In Section 2, we introduce the relevant knowledge of homogeneous Besov spaces and some important inequalities  that are used in the later analysis. In Section 3, we derive the necessary a priori bounds for smooth, local-in-time solutions to  \eqref{1.1'}--\eqref{1.3'} that are needed to extend the local solution to the global one.
With all the a priori estimates, the proof of Theorem \ref{The1.1} is completed in Section 4. In the appendix, we derive the incompressible viscoelastic system heat-conducting including the homogeneous or inhomogeneous case.

\section{Preliminaries}

\

In this section, we recall the properties of homogeneous Besov spaces and some auxiliary inequalities  that will be frequently used later.

\subsection{Littlewood-Paley theory and Besov spaces}

\

For the convenience of the readers, we recall some basic facts on the Littlewood-Paley
theory, one may check \cite{Bahouri-Chemin-Danchin2011}  for more details. Let us first introduce the Littlewood-Paley decomposition.
The symbol $\mathcal S(\mathbb R^n)$ stands for the Schwartz space of smooth functions over $\mathbb R^n$ whose derivatives of all order
decay at infinity. The set $\mathcal S'(\mathbb R^n)$ of temperate distributions is the dual set of $\mathcal{S}(\mathbb{R}^n)$ for the usual pairing.
Let $(\varphi,\chi)$ be a couple of smooth functions valued in $[0,1]$, such that $\varphi$ is supported in the shell $\{\xi\in{\mathbb R}^n|\; \frac34\leq|\xi|\leq\frac83\}$, $\chi$ is supported in the ball $\{\xi\in{\mathbb R}^n|\; |\xi|\leq\frac43\}$ and
\begin{gather*}
\forall\xi\in\mathbb{R}^n,\quad\chi(\xi)+\sum_{j\in\mathbb{N}}\varphi(2^{-j}\xi)=1,\\
\forall \xi \in \mathbb{R}^{d} \backslash\{0\}, \sum_{j \in \mathbb{Z}} \varphi\left(2^{-j} \xi\right)=1, \\
\left|j-j^{\prime}\right| \geq 2 \Rightarrow \operatorname{supp} \varphi\left(2^{-j} \cdot\right) \cap \operatorname{supp} \varphi\left(2^{-j^{\prime}} \cdot\right)=\emptyset, \\
j \geq 1 \Rightarrow \operatorname{supp} \chi \cap \operatorname{supp} \varphi\left(2^{-j} \cdot\right)=\emptyset.
\end{gather*}
Denoting ${h}= \mathcal{F}^{-1}\varphi$ and $\widetilde h = \mathcal{F}^{-1}\chi $, the homogeneous dyadic blocks $\dot\Delta_j$ and the homogeneous low-frequency cut-off operators $\dot S_j$  are defined for all $j\in\mathbb Z$ and $f\in\mathcal S'(\mathbb R^n)$ by
\begin{align*}
\dot{\Delta}_{j}f=\varphi(2^{-j}D)f=2^{jn}\int_{\mathbb{R}^n}h(2^jy)f(x-y) dy,\\
\dot{S}_{j}f=\chi(2^{-j}D)f=2^{jn}\int_{\mathbb{R}^n}\widetilde{h}(2^jy)f(x-y) dy .
\end{align*}
We denote by $\mathcal S'_h(\mathbb R^n)$ the space of temperate distributions $f$ such that
\begin{equation*}
\lim\limits_{j\to-\infty}\|\dot{S}_jf\|_{L^\infty}=0.
\end{equation*}

For $f\in \mathcal S'_h(\mathbb R^n)$, it holds that $
f=\sum_j\dot{\Delta}_jf.$ In addition, due to the definition of $\varphi$ and $\chi$, one can easily verify that
\begin{gather}
\dot{\Delta}_{j'}\dot{\Delta}_{j}f=0\quad \text{if} \quad|j'-j|\geq2,\nonumber\\
\dot{\Delta}_{j'}(\dot S_{j-1}f\dot{\Delta}_{j}f)=0 \quad \text{if} \quad|j'-j|\geq5.\label{0}
\end{gather}

We now give the definition of homogeneous Besov spaces:
\begin{Definition}
Let $s\in \mathbb R, 1\leq p,r\leq\infty$. The homogeneous Besov space $\dot B_{p,r}^s(\mathbb R^n)$ is defined by
\begin{equation*}
\dot{B}_{p,r}^{s}(\mathbb R^n)=\{f\in \mathcal S'_h(\mathbb R^n):\|f\|_{\dot{B}^{s}_{{p,r}}(\mathbb R^n)}<+\infty\},
\end{equation*}
where
\begin{equation*}
\|f\|_{\dot B_{p,r}^s(\mathbb R^n)}=\left(\sum_{q\in\mathbb Z}(2^{qs}\|\dot\Delta_qf\|_{L^p(\mathbb R^n)})^r\right)^{\frac1r}.
\end{equation*}
\end{Definition}
Next we recall the homogeneous Bony's decomposition:
\begin{equation}\label{bony}
fg=\dot{T}_fg+\dot{T}_gf+\dot{R}(f,g),
\end{equation}
with
\begin{equation*}
\dot{T}_fg=\sum_q\dot{S}_{q-1}f \dot{\Delta}_qg,\quad \dot{R}(f,g)=\sum_{|p-q|\leq1}\dot{\Delta}_pf \dot{\Delta}_qg.
\end{equation*}
Then there is the following lemma.
\begin{Lemma}\label{TB}
Let $s,s_1,s_2>0$ and  $1\leq p,p_1,p_2,r,r_1,r_2\leq \infty$. Assume that
\begin{equation*}
\frac1{p'}\triangleq\frac1{p_1}+\frac1{p_2}\leq 1 \quad\text{and} \quad \frac1{r'}\triangleq\frac1{r_1}+\frac1{r_2}\leq 1.
\end{equation*}
Then there exists a constant $C$, depending only on $n$ and $s$, such that for any $u_1,v_1\in\dot B_{p,r}^s(\mathbb R^n)\cap L^\infty(\mathbb R^n)$ and $(u_2,v_2)\in\dot B_{p_1,r_1}^{s_1}(\mathbb R^n)\times \dot B_{p_2,r_2}^{s_2}(\mathbb R^n)$,
\begin{gather*}
\|\dot{T}_{u_1}{v_1}\|_{\dot B_{p,r}^s(\mathbb R^n)}\leq C\|u_1\|_{L^\infty(\mathbb R^n)}\|v_1\|_{\dot B_{p,r}^s(\mathbb R^n)},\nonumber\\
\|u_1v_1\|_{\dot B_{p,r}^s(\mathbb R^n)}\leq C (\|u_1\|_{L^\infty(\mathbb R^n)}\|v_1\|_{\dot B_{p,r}^s(\mathbb R^n)}+\|v_1\|_{L^\infty(\mathbb R^n)}\|u_1\|_{\dot B_{p,r}^s(\mathbb R^n)}),\nonumber\\
\|\dot{R}(u_2,v_2)\|_{\dot B_{p',r'}^{s_1+s_2}(\mathbb R^n)}\leq C\|u_2\|_{\dot B_{p_1,r_1}^{s_1}(\mathbb R^n)}\|v_2\|_{\dot B_{p_2,r_2}^{s_2}(\mathbb R^n)}.
\end{gather*}
\end{Lemma}
The following Bernstein lemma will be frequently used.
\begin{Lemma}\label{ber}
Let $k\in\mathbb N, 1 \leq p \leq q \leq +\infty$. Assume that $f \in L^{p}(\mathbb{R}^{n})$. Then there exist constants $C_1, C_2 $ independent of $f, j$ such that
\begin{gather*}
\operatorname{supp}\mathcal F(f)\subseteq\{|\xi|\leq A_02^j\}\Rightarrow
\sup\limits_{|\alpha|=k}\|\partial^\alpha f\|_{L^q(\mathbb R^n)}\leq C_12^{jk+nj(\frac1p-\frac1q)}\|f\|_{L^p(\mathbb R^n)},
\\\operatorname{supp}\mathcal F(f)\subseteq\{A_{1}2^{j}\leq|\xi|\leq A_{2}2^{j}\}\Rightarrow
C_{2}2^{jk}\|f\|_{L^p(\mathbb R^n)}\leq \sup_{|\alpha|=k}\|\partial^{\alpha}f\|_{L^p(\mathbb R^n)}\leq C_22^{jk}\|f\|_{L^p(\mathbb R^n)}.
\end{gather*}
\end{Lemma}

\subsection{Some important inequalities}

\

First, the
well-known Gagliardo-Nirenberg inequality (cf. \cite[THEOREM, p.125]{Nirenberg1959}) is stated as follows:
\begin{Lemma}
Let $m$ be a positive integer and $1\leq  q, r \leq+\infty$. Then there is a constant $C(n,\alpha, m,q,r)$ such that for any integer $\alpha$ with $0\leq \alpha< m$,
\begin{equation}\label{GN}
\|\na^\alpha f\|_{L^{p}(\mathbb{R}^n)} \leq C(n,\alpha, m,q,r)\|\nabla^mf\|_{L^{q}(\mathbb{R}^n)}^{\vartheta}\|f\|_{L^{r}(\mathbb{R}^n)}^{1-\vartheta},
\end{equation}
where $\frac{\alpha}{m} \leq \vartheta \leq 1 $ and
\begin{equation*}
\frac{1}{p}-\frac \alpha n=\left(\frac{1}{q}-\frac{m}{n}\right) \vartheta+\frac{1}{r}(1-\vartheta).
\end{equation*}
Notice that there are the following exceptional cases:

(1) If $\alpha=0$, $qm<n$ and $r=+\infty$, then one additionally needs that either $f(x)\to 0$ as $|x|\to+\infty$ or $|f|^{\tilde{q}}\in L^1(\mathbb R^n)$ for some $\tilde{q}\in (0,+\infty).$

(2) If $1<q<+\infty$, and $m-\alpha-\frac{n}{q}$ is a nonnegative integer, then \eqref{GN} holds only for $\vartheta$ satisfying $\frac\alpha m \leq \vartheta < 1$.
\end{Lemma}

Then we give the Sobolev's embedding inequalities; see \cite[THEOREM 4, p.280]{Evans2010}.
\begin{Lemma}
If $p>n$, then there is a constant $C(n,p)$ such that
\begin{gather}
\|f\|_{L^\infty(\mathbb R^n)}\leq C(n,p)\|f\|_{W^{1,p}(\mathbb R^n)},\label{em1}\\
\langle f\rangle^{\beta}\leq C(n,p)\|\nabla f\|_{L^p(\mathbb{R}^{n})}\label{em2}
\end{gather}
for given $f\in W^{1,p}(\mathbb R^n)$ and $\beta=1-\frac n p$.
\end{Lemma}


Next, we introduce the following well-known estimate (cf. \cite[THEOREM II.1, p.250]{Coifman-Lions-Meyer-Semmes1993}):
\begin{Lemma}\label{RP}
If $v\in L^2({\mathbb{R}^n})$, $w\in H^1({\mathbb{R}^n})$ and $\operatorname{div}v=0$, then $ v\cdot\nabla w\in\mathcal{H}^1({\mathbb{R}^n}) $. In addition, there is  a constant $C(n)$ such that
\begin{equation*}\label{hd}
\|v\cdot\nabla w\|_{\mathcal{H}^1({\mathbb{R}^n})}\leq C(n)\|v\|_{L^2({\mathbb{R}^n})}\|\nabla w\|_{L^2({\mathbb{R}^n})}.
\end{equation*}
\end{Lemma}
Finally, we give the classical commutator estimate (cf. \cite[Lemma 2.97, p.110]{Bahouri-Chemin-Danchin2011}):
\begin{Lemma}\label{jh}Let  $f$  be a  $C^{1}$  function on  $\mathbb{R}^{n}$  such that $ (1+|\cdot|) \mathcal F(f) \in L^{1}(\mathbb{R}^{n})$. There exists a constant  $C$  such that for any Lipschitz function $g$ with gradient in  $L^{p}(\mathbb{R}^{n})$  and any function $h$ in  $L^{q}(\mathbb{R}^{n})$, we have, for any positive  $\lambda$ ,
\begin{equation*}
\|[f(\lambda^{-1} D), g] h\|_{L^{r}(\mathbb{R}^{n})} \leq C \lambda^{-1}\|\nabla g\|_{L^{p}(\mathbb{R}^{n})}\|h\|_{L^{q}(\mathbb{R}^{n})} \quad \text { with } \quad \frac{1}{r}\triangleq\frac{1}{p}+\frac{1}{q},
\end{equation*}
where $f(D)g\triangleq \mathcal F^{-1}(f\mathcal F g)$.
\end{Lemma}


\section{A priori estimates}

\

In this section, we would establish a number of a priori estimates for local-in-time smooth solutions, corresponding roughly to \eqref{ab}. Since these processes are rather long and technical, we will omit details which are identical to or nearly identical to arguments given in  \cite{Hu-Lin2016}. To begin with, let $T>0$ be fixed and assume that $(u,\mathbb{F}, \theta)$ is a smooth solution of \eqref{1.1'}--\eqref{1.3'} defined on $\mathbb{R}^n\times [0,T]$, then we have the following key a priori estimates:
\begin{Proposition}\label{gj}
Fix constants $q>4$, $N>0$ (not necessarily small), $\tilde{\theta}>0$, $s\in(0,1)$ when $n=2$ and $s\in(\frac{8}{15},1)$ when $n=3$, $p>2$ when $n=2$ and $p\in(3,6]$ when $n=3$, then there are positive constants $\tilde{\delta}$, $b$, $d$, $C$, $M$ and $\alpha$ depending only on $p, q, N, s$ and $\tilde{\theta}$ such that if $\varepsilon_{0}\leq d$ and
\begin{gather}
A_1(T)+A_2(T)+A_4(T)+B(T)\leq \tilde{\delta},\label{12}\\
|(u(x,t),\theta(x,t))-(0,\tilde{\theta})|\leq b,\quad \min\{1,T\}\leq t\leq T,\label{uc}
\end{gather}
the following estimates hold:
\begin{gather*}
A_1(T)+A_2(T)+A_4(T)+B(T)\leq C\varepsilon_{0}^\alpha, \quad A_3(T)\leq M,\label{2.3'}\\
\|u\|_{L^\infty(\mathbb{R}^n\times[\tau,T])}+\langle u\rangle_{\mathbb{R}^{n} \times[\tau, T]}^{\beta, \frac{\beta}{n+2\beta}} \leq C(\tau,\beta) \varepsilon_{0}^\alpha, \quad \tau >0,\label{ug}\\
\|\theta-\tilde{\theta}\|_{L^\infty(\mathbb{R}^n\times[\tau,T])}+\langle \theta\rangle_{\mathbb{R}^{n} \times[\tau, T]}^{\beta, \frac{\beta}{n+2\beta}}\leq C(\tau,\beta) \varepsilon_{0}^\alpha, \quad \tau >0,\label{eg}
\end{gather*}
where $\beta\in (0,1)$ when $n=2$ and $\beta\in (0,\frac12]$ when $n=3$, $C(\tau,\beta)$ is a constant that depends additionally on  $\tau$ and $\beta$.
\end{Proposition}
\begin{proof}
This proposition is a direct result of combining Lemma \ref{A1}, Lemma \ref{theta} and Lemma \ref{pt}.
\end{proof}

\begin{Lemma}\label{A1}
Assume that the hypotheses and notations of Proposition \ref{gj} are in force. Then
\begin{equation}
A_1(T)+A_2(T)+B(T)\leq C\varepsilon_{0}^\alpha, \quad A_3(T)\leq M.\label{3.2'}
\end{equation}
\end{Lemma}
\begin{proof}
The proof of this lemma is carried out in Lemma \ref{A1d}, Lemma \ref{2.2}, Lemma \ref{2.3}, Lemma \ref{A3} and  Lemma \ref{bt}.
\end{proof}

To begin with, we derive the following standard energy estimate for $(u,\mathbb F)$ and preliminary $L^2-$ bounds for $\na u$, $\mathcal P\dot u$ and $\na\mathcal P\dot u$.
\begin{Lemma}\label{A1d}
Assume that the hypotheses and notations of Proposition \ref{gj} are in force. Then
\begin{equation}\label{2.4}
A_1(T)\leq C\left(\varepsilon_{0}+\int_0^T\int_{\mathbb{R}^{n}}\sigma^n|\nabla u|^4\;dxdt+\sum\left|\int_0^T\int_{\mathbb{R}^{n}}\sigma u_i^lu_j^mu_k^r\;dxdt\right|\right).
\end{equation}
\end{Lemma}
\begin{proof}
The proof of this lemma is exactly similar to section 2 in Hu-Lin \cite{Hu-Lin2016}, which consists of three separate energy-type estimates.

\textit{Step} 1. Adding $(\ref{1.1'}) _1$ multiplied by $u$ and $(\ref{1.1'})_2$ multiplied by $(\mathbb F-\mathbb I)$ leads to
\begin{align*}
&\frac{1}{2} (|u|^2+|\mathbb F-\mathbb I|^2)_t+\frac{1}{2}u\cdot\nabla(|u|^2+|\mathbb F-\mathbb I|^2)-\Delta u\cdot u+\nabla P\cdot u=(\mathbb F^{ik}\mathbb F^{jk}u^i)_j-u_j^k\mathbb F^{jk}.
\end{align*}
Since $\operatorname{div}u=0 $ and $\operatorname{div}\mathbb F^{\top}=0$, we obtain following after integrating the above equation over $\mathbb{R}^n$:
\begin{equation*}
\frac{1}{2} \frac{d}{dt} \int_{\mathbb{R}^{n}}(|u|^2+|\mathbb F-\mathbb I|^2)\;dx+\int_{\mathbb{R}^{n}}|\nabla u|^2\;dx=0.
\end{equation*}
Integrating the above equation over $[0,t]$ for $t\in [0,T]$, one gets
\begin{equation}\label{nl}
\sup\limits_{0\leq t\leq T}\int_{\mathbb{R}^{n}}(|u|^2+|\mathbb F-\mathbb I|^2)\;dx+\int_0^T\int_{\mathbb{R}^{n}}|\nabla u|^2\;dxdt\leq
C\int_{\mathbb{R}^{n}}(|u_0|^2+|\mathbb F_0-\mathbb I|^2)\;dx\leq C\varepsilon_0.
\end{equation}

\textit{Step} 2: Next we show the estimate for $\sigma \int_{\mathbb{R}^{n}}|\nabla u|^2\;dx$ and $\int_0^T\int_{\mathbb{R}^{n}}\sigma |\mathcal{P}\dot {u}|^2\;dxdt$.
Applying $\mathcal{P}$ to (\ref{1.1'})$_1$, and noticing $\mathcal{P}\nabla P=0$ together with $\mathcal{P}\Delta u=\Delta u$ imply
\begin{equation}\label{pzy}
\mathcal{P}\dot{u}=\Delta u+\mathcal{P}\operatorname{div}(\mathbb F\mathbb F^{\top}).
\end{equation}
Multiplying the above equation by $\sigma\mathcal{P}\dot{u}$ and integrating the resulting equality over $\mathbb{R}^n\times [0, t]$ for $t\in[0,T]$, we have
\begin{align*}
\int_0^t\int_{\mathbb{R}^{n}}\sigma |\mathcal{P}\dot {u}|^2\;dxdt'=  \int_0^t\int_{\mathbb{R}^{n}}\Delta u\cdot \sigma\mathcal{P} \dot{u}\;dxdt'+\int_0^t\int_{\mathbb{R}^{n}}\mathcal{P}\operatorname{div}(\mathbb F\mathbb F^{\top})\cdot \sigma\mathcal{P} \dot{u}\;dxdt'\triangleq I_1+I_2.
\end{align*}
Assume $\mathcal{O}=\mathbb Id-\mathcal{P}$, then it is easy to check that operators $ \mathcal{O} $ and $ \mathcal{P} $ are orthogonal.

For $I_1$, the orthogonality between $ \mathcal{P} $ and $ \mathcal{O} $ yields
\begin{align*}
I_1&=\int_0^t\int_{\mathbb{R}^{n}}\sigma \Delta u\cdot(u_{t'}+u\cdot \nabla u-\mathcal{O}(u\cdot \nabla u))\;dxdt'\\
&=\int_0^t\int_{\mathbb{R}^{n}}\sigma \Delta u\cdot u_{t'}\;dxdt'+\int_0^t\int_{\mathbb{R}^{n}}\sigma \Delta u\cdot(u\cdot\nabla u)\;dxdt'\\
&\triangleq I_{11}+I_{12}.
\end{align*}
Then it follows from integration by parts and \eqref{nl} that
\begin{align*}
I_{11}&=-\int_0^t(\frac{d}{dt'}\int_{\mathbb{R}^{n}}\frac{\sigma}{2}|\nabla u|^2\;dx-\frac{1}{2}\int_{\mathbb{R}^{n}}\sigma_{t'}|\nabla u|^2\;dx)dt'\\
&=-\frac{\sigma}{2}\int_{\mathbb{R}^{n}}|\nabla u|^2\;dx+\frac{1}{2}\int_0^\sigma\int_{\mathbb{R}^{n}}|\nabla u|^2\;dxdt'\\
&\leq -\frac{\sigma}{2}\int_{\mathbb{R}^{n}}|\nabla u|^2\;dx+C\varepsilon_{0}.
\end{align*}
In view of $\operatorname{div}u=0 $, we can bound $I_{12}$ as
\begin{align*}
I_{12}&=\int_0^t\int_{\mathbb{R}^{n}}\sigma \Delta u\cdot(u\cdot\nabla u)\;dxdt'
=\int_0^t\int_{\mathbb{R}^{n}}\sigma u_{ii}^lu^mu_m^l\;dxdt'
\leq C\sum\left|\int_0^t\int_{\mathbb{R}^{n}}\sigma u_i^lu_j^mu_k^r\;dxdt'\right|.  	
\end{align*}
Thus
\begin{align*}
I_1\leq C\varepsilon_0-\frac{\sigma}{2}\int_{\mathbb{R}^{n}}|\nabla u|^2\;dx+ C\sum\left|\int_0^t\int_{\mathbb{R}^{n}}\sigma u_i^lu_j^mu_k^r\;dxdt'\right|.
\end{align*}

For $I_2$, the orthogonality between $ \mathcal{P} $ and $ \mathcal{O} $ implies
\begin{align*}
I_2&=\int_0^t\int_{\mathbb{R}^{n}}\sigma\mathcal{P}\operatorname{div}(\mathbb F\mathbb F^{\top})\cdot \mathcal{P} \dot{u}\;dxdt'\\
&=\int_0^t\int_{\mathbb{R}^{n}}\sigma\mathcal{P}\operatorname{div}(\mathbb F\mathbb F^{\top})\cdot (\mathbb Id-\mathcal{O}) \dot{u}\;dxdt'\\
&=\int_0^t\int_{\mathbb{R}^{n}}\sigma \operatorname{div}(\mathbb F\mathbb F^{\top})\cdot \dot{u}\;dxdt'-\int_0^t\int_{\mathbb{R}^{n}}\sigma\mathcal{O} \operatorname{div}(\mathbb F\mathbb F^{\top})\cdot\dot u\;dxdt'\\
&\triangleq I_{21}+I_{22}.
\end{align*}
Since $\operatorname{div}\mathbb{F}^{\top}=0, $ $\operatorname{div}(\mathbb{F}\mathbb{F}^{\top})=\operatorname{div}
((\mathbb{F}-\mathbb I)(\mathbb{F}^{\top}-\mathbb I+\mathbb I))=\operatorname{div}((\mathbb{F}-\mathbb I)(\mathbb{F}-\mathbb I)^\top)+\operatorname{div}(\mathbb{F}-\mathbb I). $ Thus,
$I_{21}$ can be computed after using integration by parts that
\begin{align} \label{Io}
I_{21}=&\int_0^t\int_{\mathbb{R}^{n}}\sigma \operatorname{div}(\mathbb F\mathbb F^{\top})\cdot \dot{u}\;dxdt'\nonumber\\
=&\int_0^t\int_{\mathbb{R}^{n}}\sigma (\operatorname{div}(\mathbb F-\mathbb I)+\operatorname{div}((\mathbb F-\mathbb I)(\mathbb F-\mathbb I)^\top))\cdot \dot{u}\;dxdt'\nonumber\\
=&-\int_0^t\int_{\mathbb{R}^{n}}\sigma (\mathbb F-\mathbb I):\nabla u_{t'}\;dxdt'-\int_0^t\int_{\mathbb{R}^{n}}\sigma (\mathbb F-\mathbb I):\nabla(u\cdot \nabla u)\;dxdt'\nonumber\\
&-\int_0^t\int_{\mathbb{R}^{n}}\sigma (\mathbb F-\mathbb I)(\mathbb F-\mathbb I)^\top:\nabla u_{t'}\;dxdt'-\int_0^t\int_{\mathbb{R}^{n}}\sigma (\mathbb F-\mathbb I)(\mathbb F-\mathbb I)^\top:\nabla(u\cdot \nabla u)\;dxdt'\nonumber\\
\triangleq &J_1+J_2+J_3+J_4.
\end{align}
For $J_1$, we have, using integration by parts,
\begin{align}
J_1&=-\int_0^t\int_{\mathbb{R}^{n}}\sigma (\mathbb F-\mathbb I):\nabla u_{t'}\;dxdt'\nonumber\\
&=\int_0^t\int_{\mathbb{R}^{n}}\sigma (\mathbb F-\mathbb I)_{t'}:\nabla u\;dxdt'+\int_0^\sigma\int_{\mathbb{R}^{n}} (\mathbb F-\mathbb I):\nabla u\;dxdt'-\int_{\mathbb{R}^{n}}\sigma (\mathbb F-\mathbb I):\nabla u\;dx\nonumber\\
&\triangleq J_{11}+J_{12}+J_{13}.
\end{align}
and analogously,
\begin{align}
J_3=&\int_0^t\int_{\mathbb{R}^{n}}\sigma ((\mathbb F-\mathbb I)(\mathbb F-\mathbb I)^\top)_{t'}:\nabla u\;dxdt'+\int_0^\sigma \int_{\mathbb{R}^{n}}(\mathbb F-\mathbb I)(\mathbb F-\mathbb I)^\top:\nabla u\;dxdt'\nonumber\\
&-\int_{\mathbb{R}^{n}}\sigma (\mathbb F-\mathbb I)(\mathbb F-\mathbb I)^\top:\nabla u\;dx\nonumber\\
\triangleq &J_{31}+J_{32}+J_{33}.
\end{align}
Based on integration by parts, (\ref{1.1'})$_2$, \eqref{12}, \eqref{nl} and   Young's inequality, we deduce that
\begin{align}\label{315}
J_{11}+J_{2}=&\int_0^t\int_{\mathbb{R}^{n}}\sigma (\mathbb F-\mathbb I)_{t'}:\nabla u\;dxdt' \nonumber-\int_0^t\int_{\mathbb{R}^{n}}\sigma (\mathbb F-\mathbb I):\nabla(u\cdot \nabla u)\;dxdt'\\
=&\int_0^t\int_{\mathbb{R}^{n}}\sigma (-u\cdot \nabla(\mathbb F-\mathbb I)+\nabla u\mathbb F):\nabla u\;dxdt'-\int_0^t\int_{\mathbb{R}^{n}}\sigma (\mathbb F-\mathbb I):\nabla(u\cdot \nabla u)\;dxdt' \nonumber\\
=&\int_0^t\int_{\mathbb{R}^{n}}\sigma \operatorname{div}u(\mathbb F-\mathbb I):\nabla u\;dxdt'+\int_0^t\int_{\mathbb{R}^{n}}\sigma (\mathbb F-\mathbb I):(u\cdot \nabla )\nabla u\;dxdt'\nonumber\\
&+\int_0^t\int_{\mathbb{R}^{n}}\sigma \nabla u\mathbb F:\nabla u\;dxdt'-\int_0^t\int_{\mathbb{R}^{n}}\sigma (\mathbb F-\mathbb I):((u\cdot \nabla)\nabla u+\nabla u\nabla u)\;dxdt' \nonumber\\
=&\int_0^t\int_{\mathbb{R}^{n}}\sigma \nabla u\mathbb F:\nabla u\;dxdt'-\int_0^t\int_{\mathbb{R}^{n}}\sigma (\mathbb F-\mathbb I):(\nabla u\nabla u)\;dxdt' \nonumber\\
\leq&C\int_0^t\int_{\mathbb{R}^{n}}|\nabla u|^2\;dxdt'\leq C\varepsilon_{0} .
\end{align}
It follows in a similar way to \eqref{315} that
\begin{equation} \label{Ic}
J_{31}+J_{4}\leq C\varepsilon_{0}.
\end{equation}
In view of \eqref{nl}, (\ref{Io})-(\ref{Ic}), H\"{o}lder's inequality and Young's inequality, one gets
\begin{align*}
I_{21}&\leq C\left(\varepsilon_0+\sigma\int_{\mathbb{R}^{n}}|\nabla u||\mathbb F-\mathbb I|\;dx+\int_0^\sigma\int_{\mathbb{R}^{n}}|\nabla u||\mathbb F-\mathbb I|\;dxdt'\right)\\
&\leq \frac{1}{4}\sigma\int_{\mathbb{R}^{n}}|\nabla u|^2\;dx+C\varepsilon_0.
\end{align*}
For $I_{22}$, since the Riesz transform is a bounded operator on the Hardy space $ \mathcal{H}^{1} (\mathbb R^n)$, \eqref{12} and Lemma \ref{RP}, we have
\begin{align*}
I_{22} &=- \int_0^t\int_{\mathbb{R}^{n}}\sigma\mathcal{O} \operatorname{div}(\mathbb F\mathbb F^{\top})\cdot\dot u\;dxdt'\\
&=\int_0^t\int_{\mathbb{R}^{n}}\sigma\nabla(-\Delta)^{-1}\operatorname{div} \operatorname{div}(\mathbb F\mathbb F^{\top})\cdot\dot u\;dxdt'\\
& \leq C \int_{0}^{t}\left\|\mathbb F\mathbb F^\top\right\|_{B M O(\mathbb R^n)}\|\nabla\nabla(-\Delta)^{-1}((\nabla u)^\top:\nabla u )\|_{\mathcal{H}^{1}(\mathbb R^n)}\;d t' \\
& \leq C\int_{0}^{t}\left\|\mathbb F\mathbb F^\top\right\|_{L^\infty(\mathbb R^n)}\|(\nabla u)^\top:\nabla u \|_{\mathcal{H}^{1}(\mathbb R^n)}\;dt' \\
& \leq C \int_{0}^{t}\|\nabla u\|_{L^{2}(\mathbb R^n)}^{2}\;dt'\\
&\leq C\varepsilon_0.
\end{align*}
Thus, collecting the estimates about $I_{21}$ and $I_{22}$ together gives
\begin{equation*}
I_2\leq  C\varepsilon_0+\frac{1}{4}\sigma\int_{\mathbb{R}^{n}}|\nabla u|^2\;dx.
\end{equation*}
Combining the bounds for $I_1$ and $I_2$ leads to
\begin{align}\label{318}
\sup_{0\leq t\leq T}\sigma\int_{\mathbb{R}^n}|\nabla u|^2\;dx+\int_{0}^{T}\int_{\mathbb{R}^n}\sigma|\mathcal{P}\dot{u}|^2\;dxdt\leq C\left(\varepsilon_{0}+\sum\left|\int_0^T\int_{\mathbb{R}^{n}}\sigma u_i^lu_j^mu_k^r\;dxdt\right|\right).
\end{align}

\textit{Step} 3. Finally, we estimate the terms $\sigma^n\int_{\mathbb{R}^{n}}|\mathcal{P}\dot {u}|^2\;dx$ and  $\int_0^T\int_{\mathbb{R}^{n}}\sigma^n|\nabla\mathcal{P}\dot{u}|^2\;dxdt$. Operating $\partial_t+u\cdot \nabla$ to $\eqref{pzy}$ yields that
\begin{align}\label{jifen}
 (\mathcal{P}\dot{u})_t+(u\cdot\nabla)\mathcal{P}\dot{u}=&\Delta u_t+(u\cdot\nabla)\Delta u+(\mathcal{P}\operatorname{div}(\mathbb F\mathbb F^{\top}))_t+(u\cdot\nabla)\mathcal{P}\operatorname{div}(\mathbb F\mathbb F^{\top})\nonumber\\
=&\Delta u_t+(u\cdot\nabla)\Delta u+\Delta u\operatorname{div}u+(\mathcal{P}\operatorname{div}(\mathbb F\mathbb F^{\top}))_t\nonumber+(u\cdot\nabla)\mathcal{P}\operatorname{div}(\mathbb F\mathbb F^{\top})+\mathcal{P}\operatorname{div}(\mathbb F\mathbb F^{\top})\operatorname{div}u\nonumber\\
=&\Delta u_t+\operatorname{div}(\Delta u\otimes u)+\mathcal{P}\operatorname{div}(\mathbb F\mathbb F^{\top})_t+\operatorname{div}(\mathcal{P}\operatorname{div}(\mathbb F\mathbb F^{\top})\otimes u),
\end{align}
where we have used $\diver u=0$.
Multiplying (\ref{jifen}) by $\sigma^n\mathcal P \dot u$ and integrating the resulting equation over $\mathbb{R}^n\times [0, t]$ for $ t\in[0,T]$,
noticing
\begin{align*}
\int_{0}^{t}\int_{\mathbb{R}^{n}}\sigma^n (\mathcal{P}\dot{u})_{t'}\cdot\mathcal{P}\dot{u}\;dxdt'
&=\frac{1}{2}\int_{0}^{t}\frac{d}{dt'}\int_{\mathbb{R}^{n}}\sigma^n |\mathcal{P}\dot{u}|^2\;dxdt'-\frac{n}{2}\int_{0}^{t}\int_{\mathbb{R}^{n}}\sigma^{n-1}\sigma'|\mathcal{P}\dot{u}|^2\;dxdt'\\
&=\frac{1}{2}\sigma^n\int_{\mathbb{R}^{n}} |\mathcal{P}\dot{u}|^2\;dxdt'-\frac{n}{2}\int_{0}^{t}\int_{\mathbb{R}^{n}}\sigma^{n-1}\sigma'|\mathcal{P}\dot{u}|^2\;dxdt',
\end{align*}
which together with
\begin{align*}
\int_{0}^{t}\int_{\mathbb{R}^{n}}\sigma^n(u\cdot\nabla\mathcal{P}\dot{u})\cdot\mathcal{P}\dot{u}\;dxdt'
=-\frac{1}{2}\int_{0}^{t}\int_{\mathbb{R}^{n}}\sigma^n\operatorname{div}u|\mathcal{P}\dot{u}|^2\;dxdt'
=0
\end{align*}
implies that
\begin{align*}
\frac{1}{2}\int_{\mathbb{R}^{n}}\sigma^n |\mathcal{P}\dot{u}|^2\;dxdt'
=&\int_{0}^{t}\int_{\mathbb{R}^{n}}\sigma^n(\Delta u_{t'}+\operatorname{div}(\Delta u\otimes u))\cdot\mathcal P \dot u\;dxdt'+\frac{n}{2}\int_{0}^{t}\int_{\mathbb{R}^{n}}\sigma^{n-1}\sigma'|\mathcal{P}\dot{u}|^2\;dxdt'\\
&+\int_{0}^{t}\int_{\mathbb{R}^{n}}\sigma^n(\mathcal{P}\operatorname{div}(\mathbb F\mathbb F^{\top})_{t'}+\operatorname{div}(\mathcal{P}\operatorname{div}(\mathbb F\mathbb F^{\top})\otimes u))\cdot\mathcal P \dot u\;dxdt'\\
\triangleq&K_1+K_2+K_3.
\end{align*}
Using integration by parts to write $K_1$ as
\begin{align*}
K_{1} & =\int_{0}^{t} \int_{\mathbb{R}^{n}} \sigma^{n} \mathcal{P} \dot{u} \cdot\left(\Delta u_{t'}+\operatorname{div}(\Delta u \otimes u)\right) dxdt' \\
& =-\int_{0}^{t} \int_{\mathbb{R}^{n}} \sigma^{n}\left(\nabla \mathcal{P} \dot{u}: \nabla u_{t'}+\nabla \mathcal{P} \dot{u}: (\Delta u \otimes {u})\right) dxdt' \\
& =-\int_{0}^{t} \int_{\mathbb{R}^{n}} \sigma^{n}\nabla \mathcal{P} \dot{u}:\nabla \dot{u}\;dxdt'-\int_{0}^{t} \int_{\mathbb{R}^{n}} \sigma^{n}
\nabla \mathcal{P} \dot{u}:(\Delta u \otimes u-\nabla(u \cdot \nabla u))\;dxdt' \\
&\triangleq K_{11}+K_{12}.
\end{align*}
We infer from the orthogonality between $\mathcal{P}$ and $\mathcal{O}$ that
\begin{align*}
K_{11}&=-\int_{0}^{t} \int_{\mathbb{R}^{n}} \sigma^{n}\nabla \mathcal{P} \dot{u}:\nabla \dot{u}\;dxdt'\\
&=-\int_{0}^{t} \int_{\mathbb{R}^{n}} \sigma^{n}\nabla \mathcal{P} \dot{u}:(\nabla \mathcal{O}\dot{u}+\nabla \mathcal{P}\dot{u})\;dxdt'\\
&=\int_{0}^{t} \int_{\mathbb{R}^{n}} \sigma^{n} \mathcal{P} \Delta\dot{u}\cdot\mathcal{O}\dot{u}\;dxdt'-\int_{0}^{t} \int_{\mathbb{R}^{n}} \sigma^{n}|\nabla \mathcal{P} \dot{u}|^2\;dxdt'\\
&=-\int_{0}^{t} \int_{\mathbb{R}^{n}} \sigma^{n}|\nabla \mathcal{P} \dot{u}|^2\;dxdt'.
\end{align*}
And integration by parts leads to
\begin{align*}
K_{12}=&-\int_{0}^{t} \int_{\mathbb{R}^{n}} \sigma^{n}\nabla \mathcal{P} \dot{u}:(\Delta u \otimes u-\nabla(u \cdot \nabla u)) \;dxdt'\\
=&\int_{0}^{t} \int_{\mathbb{R}^{n}} \sigma^{n}[(u \cdot \nabla)\nabla \mathcal{P} \dot{u}:\nabla u+\nabla\mathcal{P} \dot{u}\nabla u :\nabla u +\nabla\mathcal{P} \dot{u}:(u\cdot\nabla)\nabla u+\nabla\mathcal{P} \dot{u}:(\nabla u \nabla u )]\;dxdt'\\
=&\int_{0}^{t} \int_{\mathbb{R}^{n}} \sigma^{n}[u\cdot \nabla (\nabla\mathcal{P} \dot{u}:\nabla u)+\nabla\mathcal{P} \dot{u}\nabla u :\nabla u+\nabla\mathcal{P} \dot{u}:(\nabla u \nabla u)]\;dxdt'.
\end{align*}
The first term of $K_{12}$ can be computed as
\begin{equation*}
\int_{0}^{t} \int_{\mathbb{R}^{n}} \sigma^{n}u\cdot \nabla (\nabla\mathcal{P} \dot{u}:\nabla u)\;dxdt'=-\int_{0}^{t} \int_{\mathbb{R}^{n}} \sigma^{n}\operatorname{div}u(\nabla\mathcal{P} \dot{u}:\nabla u)\;dxdt'
=0,
\end{equation*}
which, together with Young's inequality, gives
\begin{align*}
K_{12}=&\int_{0}^{t} \int_{\mathbb{R}^{n}} \sigma^{n}(\nabla\mathcal{P} \dot{u}\nabla u :\nabla u+\nabla\mathcal{P} \dot{u}:(\nabla u \nabla u))\;dxdt'\\
\leq &\frac{1}{4}\int_{0}^{t} \int_{\mathbb{R}^{n}} \sigma^{n}|\nabla \mathcal{P} \dot{u}|^2\;dxdt'+C\int_{0}^{t} \int_{\mathbb{R}^{n}} \sigma^{n}|\nabla u|^4\;dxdt'.
\end{align*}
\par Combining the bounds for $K_{11}$ and $K_{12}$ yields
\begin{align*}
K_1\leq -\frac{3}{4}\int_{0}^{t} \int_{\mathbb{R}^{n}} \sigma^{n}|\nabla \mathcal{P} \dot{u}|^2\;dxdt'+C\int_{0}^{t} \int_{\mathbb{R}^{n}} \sigma^{n}|\nabla u|^4\;dxdt'.
\end{align*}
\par For $K_2$, since $|\sigma|\leq 1, |\sigma'|\leq 1$, one has
\begin{align*}
K_2=\frac{n}{2}\int_{0}^{t}\int_{\mathbb{R}^{n}}\sigma^{n-1}\sigma'|\mathcal{P}\dot{u}|^2\;dxdt'\leq C\int_{0}^{t}\int_{\mathbb{R}^{n}}\sigma|\mathcal{P}\dot{u}|^2\;dxdt'.
\end{align*}
\par Now we deduce the bound for $K_3$. Transpose (\ref{1.1'})$_2$ to get
\begin{equation} \label{zhuanzhi}
\partial_t \mathbb F^{\top}+u\cdot \nabla \mathbb F^{\top}=\mathbb F^{\top}(\nabla u)^{\top},
\end{equation}
Adding (\ref{1.1'})$_2$ right multiplied by $\mathbb F^{\top}$ and (\ref{zhuanzhi}) left multiplied by $\mathbb F$, we obtain
\begin{align}\label{jian}
\partial_t(\mathbb F\mathbb F^{\top})+u\cdot \nabla(\mathbb F\mathbb F^{\top})=\nabla u\mathbb F\mathbb F^{\top}+\mathbb F\mathbb F^{\top}(\nabla u)^{\top}.
\end{align}
We deduce from  the orthogonality between $\mathcal{P}$ and $\mathcal{O}$ together with integration by parts that
\begin{align*}
K_3=&\int_{0}^{t}\int_{\mathbb{R}^{n}}\sigma^n(\mathcal{P}\operatorname{div}(\mathbb F\mathbb F^{\top})_t+\operatorname{div}(\mathcal{P}\operatorname{div}(\mathbb F\mathbb F^{\top})\otimes u))\cdot\mathcal P \dot u\;dxdt'\\
=&\int_{0}^{t}\int_{\mathbb{R}^{n}}\sigma^n((Id-\mathcal{O})\operatorname{div}(\mathbb F\mathbb F^{\top})_t+\operatorname{div}(\mathcal{P}\operatorname{div}(\mathbb F\mathbb F^{\top})\otimes u))\cdot\mathcal P \dot u\;dxdt'\\
=&\int_{0}^{t}\int_{\mathbb{R}^{n}}\sigma^n(\operatorname{div}(\mathbb F\mathbb F^{\top})_t+\operatorname{div}(\mathcal{P}\operatorname{div}(\mathbb F\mathbb F^{\top})\otimes u))\cdot\mathcal P \dot u\;dxdt'\\
=&-\int_{0}^{t}\int_{\mathbb{R}^{n}}\sigma^n((\mathbb F\mathbb F^{\top})_t+\operatorname{div}(\mathbb F\mathbb F^{\top})\otimes u):\nabla\mathcal P \dot u\;dxdt'\\
&+\int_{0}^{t}\int_{\mathbb{R}^{n}}\sigma^n\mathcal{O}\operatorname{div}(\mathbb F\mathbb F^{\top})\otimes u:\nabla\mathcal P \dot u\;dxdt'\\
\triangleq&K_{31}+K_{32}.
\end{align*}
Due to (\ref{jian}) and integration by parts, $K_{31}$ can be written as
\begin{align*}
K_{31}
=&-\int_{0}^{t}\int_{\mathbb{R}^{n}}\sigma^n(\mathbb F\mathbb F^{\top}(\nabla u)^{\top}+\nabla u\mathbb F\mathbb F^{\top}):\nabla\mathcal P \dot u\;dxdt'\\
&-\int_{0}^{t}\int_{\mathbb{R}^{n}}\sigma^n(-u\cdot\nabla (\mathbb F\mathbb F^{\top})+\operatorname{div}(\mathbb F\mathbb F^{\top})\otimes u):\nabla\mathcal P \dot u\;dxdt'.
\end{align*}
Owing to  \eqref{12}, \eqref{nl},  integration by parts, H\"{o}lder's inequality and Young's inequality, the first term of $K_{31}$ is bounded by
\begin{align*}
C\int_{0}^{t}\sigma^n\left(\int_{\mathbb{R}^{n}}|\nabla \mathcal P \dot{u}|^2\;dx\right)^\frac{1}{2}\left(\int_{\mathbb{R}^{n}}|\nabla u|^2\;dx\right)^\frac{1}{2}dt'
\leq C\varepsilon_{0}+\frac{1}{8}\int_{0}^{t}\int_{\mathbb{R}^{n}}\sigma^n|\nabla \mathcal P \dot{u}|^2\;dxdt',
\end{align*}
and the second term of $K_{31}$ is bounded by
\begin{align*}
&-\int_{0}^{t}\int_{\mathbb{R}^{n}}\sigma^n(-u\cdot\nabla (\mathbb F\mathbb F^{\top})+\operatorname{div}(\mathbb F\mathbb F^{\top})\otimes u):\nabla\mathcal P \dot u\;dxdt'\\
=&-\int_{0}^{t}\int_{\mathbb{R}^{n}}\sigma^n (\mathbb F\mathbb F^{\top}):(u\cdot\nabla)\nabla\mathcal P \dot u\;dxdt'+\int_{0}^{t}\int_{\mathbb{R}^{n}}\sigma^n (\mathbb F\mathbb F^{\top}):\nabla(u\cdot\nabla\mathcal P \dot u)\;dxdt'\\
=&-\int_{0}^{t}\int_{\mathbb{R}^{n}}\sigma^n (\mathbb F\mathbb F^{\top}):(u\cdot\nabla)\nabla\mathcal P \dot u\;dxdt'+\int_{0}^{t}\int_{\mathbb{R}^{n}}\sigma^n (\mathbb F\mathbb F^{\top}):(u\cdot\nabla)\nabla\mathcal P \dot u\;dxdt'\\
&+\int_{0}^{t}\int_{\mathbb{R}^{n}}\sigma^n (\mathbb F\mathbb F^{\top}):(\nabla\mathcal P \dot u\nabla  u)\;dxdt'\\
\leq&C\int_{0}^{t}\sigma^n\left(\int_{\mathbb{R}^{n}}|\nabla \mathcal P \dot{u}|^2\;dx\right)^{\frac{1}{2}}\left(\int_{\mathbb{R}^{n}}|\nabla u|^2\;dx\right)^{\frac{1}{2}}dt'\\
\leq&C\varepsilon_{0}+\frac{1}{8}\int_{0}^{t}\int_{\mathbb{R}^{n}}\sigma^n|\nabla \mathcal P \dot{u}|^2\;dxdt'.
\end{align*}
The above two inequalities imply that
\begin{align*}
K_{31}\leq C\varepsilon_{0}+\frac{1}{4}\int_{0}^{t}\int_{\mathbb{R}^{n}}\sigma^n|\nabla \mathcal P \dot{u}|^2\;dxdt'.
\end{align*}
For $K_{32}$, using $\operatorname{div}(\mathcal{P}\dot{u})=0$, \eqref{12}, \eqref{nl}, the boundedness of Riesz transforms in Hardy space $\mathcal H^1(\mathbb R^n)$, integration by parts and Young's inequality, we have,
\begin{align*}
K_{32}
=&\int_{0}^{t}\int_{\mathbb{R}^{n}}\sigma^n\nabla(-\Delta)^{-1}\operatorname{div}\operatorname{div}(\mathbb F\mathbb F^{\top})\otimes u:\nabla\mathcal P \dot u\;dxdt'\\
=&\int_{0}^{t}\int_{\mathbb{R}^{n}}\sigma^n(u\cdot\nabla)\nabla(-\Delta)^{-1}\operatorname{div}\operatorname{div}(\mathbb F\mathbb F^{\top}) \cdot\mathcal P \dot u\;dxdt'\\
=&\int_{0}^{t}\int_{\mathbb{R}^{n}}\sigma^n(-\Delta)^{-1}\operatorname{div}\operatorname{div}(\mathbb F\mathbb F^{\top}) \operatorname{div}(u\cdot\nabla\mathcal P \dot u)\;dxdt'\\
\leq&  C\int_{0}^{t}\sigma^n\|\mathbb F\mathbb F^{\top}\|_{BMO(\mathbb{R}^{n})}\|\nabla\nabla(-\Delta)^{-1}((\nabla u)^\top:\nabla\mathcal{P}\dot{u})\|_{\mathcal{H}^1(\mathbb{R}^{n})}\;dt'\\
\leq&C\int_{0}^{t}\sigma^n\|\mathbb F\mathbb F^{\top}\|_{L^\infty(\mathbb{R}^{n})}\|(\nabla u)^\top:\nabla\mathcal{P}\dot{u}\|_{\mathcal{H}^1(\mathbb{R}^{n})}\;dt'\\
\leq& C\int_{0}^{t}\sigma^n \|\nabla\mathcal P \dot u\|_{L^2(\mathbb{R}^{n})} \|\nabla u\|_{L^2(\mathbb{R}^{n})}\;dt'\\
\leq& C\varepsilon_{0}+\frac{1}{4}\int_{0}^{t}\int_{\mathbb{R}^{n}}\sigma^n|\nabla \mathcal P \dot{u}|^2\;dxdt'.
\end{align*}
Therefore,
\begin{align*}
K_{3}
\leq C\varepsilon_{0}+\frac{1}{2}\int_{0}^{t}\int_{\mathbb{R}^{n}}\sigma^n|\nabla \mathcal P \dot{u}|^2\;dxdt'.
\end{align*}
\par Combining the estimates of $K_{1}$, $K_2$ and $K_{3}$ yields
\begin{align}\label{2-s}
&\sup\limits_{0\leq t\leq T}\int_{\mathbb{R}^{n}}\sigma^n|\mathcal{P}\dot {u}|^2\;dx+
\int_0^T\int_{\mathbb{R}^{n}}\sigma^n|\nabla\mathcal{P}\dot{u}|^2\;dxdt\nonumber\\
\leq& C\left(\varepsilon_{0}+\int_0^T\int_{\mathbb{R}^{n}}\sigma^n|\nabla u|^4\;dxdt+\sum\left|\int_0^T\int_{\mathbb{R}^{n}}\sigma u_i^lu_j^mu_k^r\;dxdt\right|\right).
\end{align}
\par To sum up, we have claimed that
\begin{equation*}
A_1(T)\leq  C\left(\varepsilon_{0}+\int_0^T\int_{\mathbb{R}^{n}}\sigma^n|\nabla u|^4\;dxdt+\sum\left|\int_0^T\int_{\mathbb{R}^{n}}\sigma u_i^lu_j^mu_k^r\;dxdt\right|\right).
\end{equation*}
Thus the proof of Lemma \ref{A1d} is finished.
\end{proof}

To further estimate $A_1(T)$, we first establish in the following lemma some properties of the effective viscous fluxes.
\begin{Lemma}
Define effective viscous fluxes by
\begin{align*}
\begin{cases}
\mathcal{G}_1\triangleq\nabla u-
(-\Delta)^{-1}\nabla{\mathcal{P}}\operatorname{div}(\mathbb{F}\mathbb{F}^{\top}-\mathbb I),\\
\mathcal{G}_2\triangleq\cur u-
(-\Delta)^{-1}\cur\diver(\mathbb{F}\mathbb{F}^{\top}-\mathbb I),
\end{cases}
\end{align*}
where $\cur u=\na u-(\na u)^\top$.
Then the effective viscous fluxes satisfy
\begin{equation}\label{ck}
\|\Delta\mathcal{G}_1\|_{L^2(\mathbb{R}^{n})}^2+\|\Delta\mathcal{G}_2\|_{L^2(\mathbb{R}^{n})}^2\leq C\|\nabla \mathcal{P}\dot{u}\|_{L^2(\mathbb{R}^{n})}^2, \quad \|\nabla\mathcal{G}_1\|_{L^2(\mathbb{R}^{n})}^2+\|\nabla\mathcal{G}_2\|_{L^2(\mathbb{R}^{n})}^2\leq C   \|\mathcal{P}\dot{u}\|_{L^2(\mathbb{R}^{n})}^2.
\end{equation}
In addition, the equation
\begin{equation}\label{324}
\na u+\mathbb E=\mathcal G+\mathcal R
\end{equation}
holds with
\begin{align*}\label{325}
\begin{cases}
\mathbb E=\mathbb F-\mathbb I,\quad
\mathcal{G}=-(-\Delta)^{-1}\nabla\operatorname{div}\mathcal{G}_2,\\
\mathcal R^{ij}=-(-\Delta)^{-1}({\mathbb{E}}^{lj}{\mathbb{E}}^{ik}-
{\mathbb{E}}^{lk}{\mathbb{E}}^{ij})_{kl}-((-\Delta)^{-1}\na\diver
(-\Delta)^{-1}\cur\diver({\mathbb{E}}{\mathbb{E}}^{\top}))^{ij}.
\end{cases}
\end{align*}
\end{Lemma}
\begin{proof}
It is easily check from $\eqref{1.1'}_1$ that
\begin{equation*}
\Delta \mathcal{G}_1=\na{\mathcal{P}}\dot{u}\quad \text{and} \quad \Delta \gl_2=\cur(\mathcal P\dot u).
\end{equation*}
Then an application of the $L^p$-estimate of elliptic system to the above equations gives \eqref{ck}.

Using $\diver{\mathbb F}^{\top}=0$, $\diver u=0$ and the Helmholtz decomposition:
\begin{equation*}
\Delta =\na\diver+\diver\cur,
\end{equation*}
one has
\begin{align*}
\diver{\mathcal G}_2=&\diver\cur u-(-\Delta)^{-1}\diver\cur\diver({\mathbb{E}}{\mathbb{E}}^{\top})-(-\Delta)^{-1}\diver\cur\diver{\mathbb{E}}\nonumber\\
=&\Delta u-(-\Delta)^{-1}\Delta\diver\mathbb E+(-\Delta)^{-1}\na\diver\diver\mathbb E-(-\Delta)^{-1}\diver\cur\diver({\mathbb{E}}{\mathbb{E}}^{\top})\nonumber\\
=&\Delta u+\diver \mathbb E-(-\Delta)^{-1}\diver\cur\diver({\mathbb{E}}{\mathbb{E}}^{\top}),
\end{align*}
Then, applying the operator $\na$ to the above equation, we get
\begin{align*}
\na\diver{\mathcal G}_2=&\na\Delta u+\na\diver \mathbb E-\na(-\Delta)^{-1}\diver\cur\diver({\mathbb{E}}{\mathbb{E}}^{\top})\\
=&\Delta(\na u+\mathbb E)-\diver\cur\mathbb E-\na(-\Delta)^{-1}\diver\cur\diver({\mathbb{E}}{\mathbb{E}}^{\top}).
\end{align*}
Hence
\begin{align}\label{329}
\na u+\mathbb E=-(-\Delta)^{-1}\na\diver{\mathcal G}_2-(-\Delta)^{-1}E-(-\Delta)^{-1}\na(-\Delta)^{-1}\diver\cur\diver({\mathbb{E}}{\mathbb{E}}^{\top}),
\end{align}
where $E^{ij}=({\mathbb E}_k^{ij}-{\mathbb E}_j^{ik})_k$.
Since
\begin{equation*}
{\mathbb F}^{lk}{\mathbb F}_l^{ij}={\mathbb F}^{lj}{\mathbb F}_l^{ik},
\end{equation*}
then
\begin{equation*}
{\mathbb E}^{lk}{\mathbb E}_l^{ij}+{\mathbb E}_k^{ij}={\mathbb E}^{lj}{\mathbb E}_l^{ik}+{\mathbb E}_j^{ik},
\end{equation*}
which combined with $\diver {\mathbb E}^\top=0$ shows that
\begin{equation*}
E^{ij}=({\mathbb E}_k^{ij}-{\mathbb E}_j^{ik})_k=({\mathbb E}^{lj}{\mathbb E}^{ik}-{\mathbb E}^{lk}{\mathbb E}^{ij})_{kl}.
\end{equation*}
Substituting the above equality into \eqref{329} gives \eqref{324}.
\end{proof}
Now we give a further estimate of $A_1(T)$ using the effective viscous fluxes in the lemma below.
\begin{Lemma}\label{2.2}
Assume that the hypotheses and notations of Proposition \ref{gj} are in force. Then
\begin{equation*}
A_1(T)\leq C(\varepsilon_{0}^\alpha+A_1(T)^2+A_2(T)^{\frac{2}{3}}).
\end{equation*}
\end{Lemma}
\begin{proof}
At first, we show that
\begin{equation}\label{u4}
\int_0^T\int_{\mathbb{R}^{n}}\sigma^n|\nabla u|^4\;dxdt\leq C(\varepsilon_{0}^\alpha+A_1(T)^2).
\end{equation}
Combining $\eqref{1.1'}_2$ and \eqref{324} leads to
\begin{align}\label{qd}
\frac{D}{Dt}(\mathbb F-\mathbb I)+\mathbb F-\mathbb I=\mathcal{G}\mathbb F+\mathcal{R}\mathbb F-(\mathbb F-\mathbb I)(\mathbb F-\mathbb I).
\end{align}	 	
After multiplying (\ref{qd}) by $4(\mathbb F-\mathbb I)|\mathbb F-\mathbb I|^2$, it follows from \eqref{12} and Young's inequality that
\begin{align}\label{budeng}
\frac{D}{Dt}|\mathbb F-\mathbb I|^4+|\mathbb F-\mathbb I|^4\leq C(|\mathcal{G}|^4+|\mathcal{R}|^4).
\end{align}
We use  integration by parts to compute $\int_0^T\int_{\mathbb{R}^n}\sigma^n\frac{D}{Dt}|\mathbb F-\mathbb I|^4\;dxdt$ that
\begin{align*}
\int_0^T\int_{\mathbb{R}^{n}}\sigma^n\frac{D}{Dt}|\mathbb F-\mathbb I|^4\;dxdt
&=-\int_0^T\int_{\mathbb{R}^{n}}n\sigma'\sigma^{n-1}|\mathbb F-\mathbb I|^4\;dxdt+\int_{\mathbb{R}^{n}}\sigma^n|\mathbb F-\mathbb I|^4\;dx\bigg|_{0} ^{T}\\
&=-n\int_0^\sigma\int_{\mathbb{R}^{n}}\sigma^{n-1}|\mathbb F-\mathbb I|^4\;dxdt+\int_{\mathbb{R}^{n}}\sigma^n(T)|\mathbb F(x,T)-\mathbb I|^4\;dx.
\end{align*}
Hence, multiplying (\ref{budeng}) by $\sigma^n$ and integrating the resulting inequality over $\mathbb{R}^n\times[0, T]$, and then using \eqref{12}, \eqref{nl} together with the boundedness of Riesz transforms in $L^p(\mathbb R^n)$ ($1< p<\infty$), yields
\begin{align*}
\int_0^T\int_{\mathbb{R}^{n}}\sigma^n|\mathbb F-\mathbb I|^4\;dxdt
&\leq C\left(\varepsilon_0+\int_0^T\int_{\mathbb{R}^{n}}\sigma^n(|\mathcal{G}|^4+|\mathcal{R}|^4)\;dxdt\right)
\\&\leq C\left(\varepsilon_0+\int_0^T\int_{\mathbb{R}^{n}}\sigma^n|\mathcal{G}|^4\;dxdt\right)+
\frac12\int_0^T\int_{\mathbb{R}^{n}}\sigma^n|\mathbb F-\mathbb I|^4\;dxdt.
\end{align*}
Owing to the definition of ${\mathcal G}_1$ and \eqref{12}, we have
\begin{align}\label{4-2}
\int_0^T\int_{\mathbb{R}^{n}}\sigma^n|\nabla u|^4\;dxdt
&\leq C\left(\int_0^T\int_{\mathbb{R}^{n}}\sigma^n|(-\Delta)^{-1}\nabla \mathcal{P}\operatorname{div}(\mathbb F\mathbb F^\top-\mathbb I)|^4\;dxdt+\int_0^T\int_{\mathbb{R}^{n}}\sigma^n|\mathcal{G}_1|^4\;dxdt\right)\nonumber\\
&\leq C\left(\int_0^T\int_{\mathbb{R}^{n}}\sigma^n|\mathbb F-\mathbb I|^4\;dxdt+\int_0^T\int_{\mathbb{R}^{n}}\sigma^n|\mathcal{G}_1|^4\;dxdt\right)\nonumber\\
&\leq
C\left(\varepsilon_0+\int_0^T\int_{\mathbb{R}^{n}}\sigma^n(|\mathcal{G}_1|^4+|\mathcal{G}_2|^4)\;dxdt\right).
\end{align}
\par An application of \eqref{GN} leads to
\begin{align}\label{335}
\int_0^T\int_{\mathbb{R}^{2}}\sigma^2|\mathcal{G}_1|^4\;dxdt
\leq&\int_0^T\sigma^2\left(\int_{\mathbb{R}^{2}}|\mathcal{G}_1|^2\;dx\right)
\left(\int_{\mathbb{R}^{2}}|\nabla\mathcal{G}_1|^2\;dx\right)dt\nonumber\\
\leq&\sup_{0\leq t\leq T}\left(\sigma\int_{\mathbb{R}^{2}}|\mathcal{G}_1|^2\;dx\right)\left(
\int_0^T\int_{\mathbb{R}^{2}}\sigma|\nabla\mathcal{G}_1|^2\;dxdt\right)
\end{align}
and
\begin{align}\label{4-3}
\int_0^T\int_{\mathbb{R}^{3}}\sigma^3|\mathcal{G}_1|^4\;dxdt
\leq&\int_0^T\sigma^3\left(\int_{\mathbb{R}^{3}}|\mathcal{G}_1|^2\;dx\right)^\frac{1}{2}\left(\int_{\mathbb{R}^{3}}|\nabla\mathcal{G}_1|^2\;dx\right)
^\frac{3}{2}dt\nonumber\\
\leq&\sup_{0\leq t\leq T}\left[\left(\sigma\int_{\mathbb{R}^{3}}|\mathcal{G}_1|^2\;dx\right)^\frac{1}{2}
\left(\sigma^3\int_{\mathbb{R}^{3}}|\nabla \mathcal{G}_1|^2\;dx\right)^\frac{1}{2}\right]\times
\int_0^T\int_{\mathbb{R}^{3}}\sigma|\nabla\mathcal{G}_1|^2\;dxdt.
\end{align}
The definition of $\mathcal{G}_1$, \eqref{12}, (\ref{nl}), \eqref{ck} and the boundedness of Riesz transforms in $L^p(\mathbb R^n)(1< p<\infty)$ imply that
\begin{gather}\label{4-4}
\sigma\int_{\mathbb{R}^{n}}|\mathcal{G}_1|^2\;dx
\leq C\left(\sigma\int_{\mathbb{R}^{n}}|\mathbb F\mathbb F^\top-\mathbb I|^2\;dx+\sigma\int_{\mathbb{R}^{n}}|\nabla u|^2\;dx\right)
\leq C(\varepsilon_{0}+A_1(T)),\nonumber\\
\displaystyle\int_0^T\int_{\mathbb{R}^{n}}\sigma|\nabla\mathcal{G}_1|^2\;dxdt
\leq 	\int_0^T\int_{\mathbb{R}^{n}}\sigma|\mathcal{P}\dot{u}|^2\;dxdt
\leq A_1(T),\nonumber\\
\vspace{1ex}\displaystyle\sigma^n\int_{\mathbb{R}^{n}}|\nabla \mathcal{G}_1|^2\;dx
\leq \sigma^n\int_{\mathbb{R}^{n}}|\mathcal{P}\dot{u}|^2\;dx\leq A_1(T). \nonumber
\end{gather}
Substituting the above inequalities into \eqref{335} and (\ref{4-3}), one gets
\begin{align*}
\int_0^T\int_{\mathbb{R}^{n}}\sigma^n|\mathcal{G}_1|^4\;dxdt
\leq C(\varepsilon_{0}^2+A_1(T)^2).
\end{align*}
Similar to the estimate for $\int_0^T\int_{\mathbb{R}^{n}}\sigma^n|\mathcal{G}_1|^4\;dxdt$, we have
\begin{equation*}
\int_0^T\int_{\mathbb{R}^{n}}\sigma^n|\mathcal{G}_2|^4\;dxdt
\leq C(\varepsilon_{0}^2+A_1(T)^2).
\end{equation*}
Substituting the above two estimates back in \eqref{4-2} yields \eqref{u4}.

Next, for the third term on the right-hand side of \eqref{2.4}, when $n=2$, from \eqref{nl}, \eqref{u4} and Young's inequality, it follows that
\begin{align*}
\left|\int_0^T\int_{\mathbb{R}^{2}}\sigma u_i^lu_j^mu_k^r\;dxdt\right|\leq
 C\left(\int_0^T\int_{\mathbb{R}^{2}}|\nabla u|^2\;dxdt+\int_1^T\int_{\mathbb{R}^{2}}\sigma^2|\nabla u|^4\;dxdt\right)\leq C(\varepsilon_{0}^\alpha+A_1(T)^2),
\end{align*}
while when $n=3$ we need a more subtle computation.
When $ T>1 $, \eqref{nl}, \eqref{u4} and Young's inequality lead to
\begin{align*}
\left|\int_1^T\int_{\mathbb{R}^{3}}\sigma u_i^lu_j^mu_k^r\;dxdt\right|
&\leq C\left(\int_1^T\int_{\mathbb{R}^{3}}|\nabla u|^2\;dxdt+\int_1^T\int_{\mathbb{R}^{3}}\sigma^3|\nabla u|^4\;dxdt\right)\\
&\leq C(\varepsilon_{0}^\alpha+A_1(T)^2),
\end{align*}
thus we only need to prove
\begin{align}\label{3x}
\left|\int_0^{\sigma(T)}\int_{\mathbb{R}^{3}}\sigma u_i^lu_j^mu_k^r\;dxdt\right|
\leq C(\varepsilon_{0}^\alpha+A_1(T)^2+A_2(T)^{\frac{2}{3}}).
\end{align}
Due to the definition of $\mathcal{G}_1$, \eqref{12} and \eqref{nl}, we have
\begin{align}\label{head}
&\left|\int_0^{\sigma(T)}\int_{\mathbb{R}^{3}}\sigma u_i^lu_j^mu_k^r\;dxdt\right|\nonumber\\
=&\left|\int_0^{\sigma(T)}\int_{\mathbb{R}^{3}}\sigma u_i^l(\mathcal{G}_1^{mj}+((-\Delta)^{-1}\nabla{\mathcal{P}}\operatorname{div}(\mathbb{F}\mathbb{F}^{\top}-\mathbb I))^{mj})(\mathcal{G}_1^{rk}+((-\Delta)^{-1}\nabla{\mathcal{P}}\operatorname{div}(\mathbb{F}\mathbb{F}^{\top}-\mathbb I))^{rk}))\;dxdt\right|\nonumber\\
\leq&\left|\int_0^{\sigma(T)}\int_{\mathbb{R}^{3}}\sigma u_i^l\mathcal{G}_1^{mj}\mathcal{G}_1^{rk}\;dxdt\right|+C\int_0^{\sigma(T)}\int_{\mathbb{R}^{3}}\sigma|\nabla u|^2|(-\Delta)^{-1}\nabla{\mathcal{P}}\operatorname{div}(\mathbb{F}\mathbb{F}^{\top}-\mathbb I)|\;dxdt\nonumber\\
&+C\int_0^{\sigma(T)}\int_{\mathbb{R}^{3}}\sigma|\nabla u||(-\Delta)^{-1}\nabla{\mathcal{P}}\operatorname{div}(\mathbb{F}\mathbb{F}^{\top}-\mathbb I)|^2\;dxdt. 
\end{align}
It follows from integration by parts, H\"{o}lder's inequality, Young's inequality,  the boundedness of Riesz transforms in $L^p(\mathbb R^n)(1< p<\infty)$, \eqref{12}, \eqref{nl} and \eqref{ck} that the first term on the right-hand side of \eqref{head} can be bounded as
\begin{align*}
\left|\int_0^{\sigma(T)}\int_{\mathbb{R}^{3}}\sigma u_i^l\mathcal{G}_1^{mj}\mathcal{G}_1^{rk}\;dxdt\right|
&\leq C\left(\int_0^{\sigma(T)}\int_{\mathbb{R}^{3}}\sigma|\nabla\mathcal{G}_1|^2\;dxdt\right)^\frac{1}{2}\left(\int_0^{\sigma(T)}\int_{\mathbb{R}^{3}}\sigma |u|^2|\mathcal{G}_1|^2\;dxdt\right)^\frac{1}{2}\nonumber\\
&\leq CA_1(T)^\frac{1}{2}\left[\int_0^{\sigma(T)}\int_{\mathbb{R}^{3}}( |u|^2|\nabla u|^2+|u|^2|(-\Delta)^{-1}\nabla{\mathcal{P}}\operatorname{div}(\mathbb{F}\mathbb{F}^{\top}-\mathbb I)|^2)\;dxdt\right]^\frac{1}{2}\nonumber\\
&\leq C(\varepsilon_{0}^\alpha+A_1(T)^2+A_2(T)^{\frac{2}{3}}),
\end{align*}
and the second term on the right-hand side of \eqref{head} can be bounded as
\begin{align*}
&\int_0^{\sigma(T)}\int_{\mathbb{R}^{3}}\sigma|\nabla u|^2|(-\Delta)^{-1}\nabla{\mathcal{P}}\operatorname{div}(\mathbb{F}\mathbb{F}^{\top}-\mathbb I)|\;dxdt\nonumber\\
\leq &C\left(\int_0^{\sigma(T)}\int_{\mathbb{R}^{3}}|\nabla u|^2\;dxdt\right)^\frac12\left(\int_0^{\sigma(T)}\int_{\mathbb{R}^{3}}\sigma^3|\nabla u|^3\;dxdt\right)^\frac13\left(\int_0^{\sigma(T)}\int_{\mathbb{R}^{3}}|\mathbb{F}\mathbb{F}^{\top}-\mathbb I|^6\;dxdt\right)^\frac16\nonumber\\
\leq&C\varepsilon_{0}^\alpha\left(\int_0^{\sigma(T)}\int_{\mathbb{R}^{3}}\sigma^3(|\mathcal {G}_1|^3+|\mathbb{F}\mathbb{F}^{\top}-\mathbb I|^3)\;dxdt\right)^\frac13\nonumber\\
\leq &C\varepsilon_{0}^\alpha+C\varepsilon_{0}^\alpha\int_0^{\sigma(T)}\left(\sigma\int_{\mathbb{R}^{3}}
|\mathcal {G}_1|^2\;dx\right)^\frac34\left(\sigma\int_{\mathbb{R}^{3}}
|\na\mathcal {G}_1|^2\;dx\right)^\frac34dt\nonumber\\
\leq&C(\varepsilon_{0}^\alpha+A_1(T)^2).
\end{align*}
Finally, owing to \eqref{12}, \eqref{nl}, Young's inequality and the boundedness of Riesz transforms in $L^p(\mathbb R^n)(1< p<\infty)$, the third term on the right-hand side of \eqref{head} can be bounded as
\begin{align*}
&\int_0^{\sigma(T)}\int_{\mathbb{R}^{3}}\sigma|\nabla u||(-\Delta)^{-1}\nabla{\mathcal{P}}\operatorname{div}(\mathbb{F}\mathbb{F}^{\top}-\mathbb I)|^2\;dxdt\nonumber\\
\leq &C\left(\int_0^{\sigma(T)}\int_{\mathbb{R}^{3}}\sigma|\nabla u|^2\;dxdt+\int_0^{\sigma(T)}\int_{\mathbb{R}^{3}}\sigma|\mathbb{F}\mathbb{F}^{\top}-\mathbb I|^4\;dxdt\right)\nonumber\\
\leq &C\varepsilon_0^\alpha.
\end{align*}
Substituting these estimates back in (\ref{head}), then \eqref{3x} is proved. We thus complete the proof of Lemma \ref{2.2}.
\end{proof}

Then,  we can close the estimate on $A_1(T)$ and $A_2(T)$ using the following lemma:
\begin{Lemma}\label{2.3}
Assume that the hypotheses and notations of Proposition \ref{gj} are in force. Then
\begin{equation*}\label{e2}
A_2(T)\leq C(\varepsilon_{0}^\alpha+A_2(T)^\frac{3}{2}).
\end{equation*}
\end{Lemma}
\begin{proof}
Multiplying $\eqref{1.1'}_1$ by $|u|^2u$ and integrating the resulting equality over $\mathbb{R}^n\times[0, t]$ for $t\in[0,\sigma(T)]$, we obtain that
\begin{align}\label{sup}		&\frac{1}{4}\int_0^t\int_{\mathbb{R}^{n}}(|u|^4)_{t'}\;dxdt'+\int_0^t\int_{\mathbb{R}^{n}}|u|^2(u\cdot\nabla u-\Delta u+\nabla P)\cdot u\;dxdt'=\int_0^t\int_{\mathbb{R}^{n}}|u|^2\operatorname{div}(\mathbb F\mathbb F^{\top})\cdot u\;dxdt'.
\end{align}
According to integration by parts and $\operatorname{div}u=0$, we get that
\begin{align}\label{shou}
\int_0^t\int_{\mathbb{R}^{n}}|u|^2(u\cdot\nabla u)\cdot u\;dxdt'
=\int_0^t\int_{\mathbb{R}^{n}}\frac{1}{4}u\cdot\nabla |u|^4\;dxdt'=-\int_0^t\int_{\mathbb{R}^{n}}\frac{1}{4}|u|^4\operatorname{div}u\;dxdt'=0
\end{align}
and
\begin{align}
-\int_0^t\int_{\mathbb{R}^n}\Delta u\cdot |u|^2u\;dxdt'
&=\int_0^t\int_{\mathbb{R}^n}\nabla (|u|^2u):\nabla u\;dxdt'
=\int_0^t\int_{\mathbb{R}^n}(|u|^2|\nabla u|^2+\frac{1}{2}|\nabla|u|^2|^2)\;dxdt'.
\end{align}
Applying $\mathcal P$ to $\eqref{1.1'}_1$ yields
\begin{align*}
\mathcal{P} \dot{u} - \mathcal{P}\operatorname{div}((\mathbb F-\mathbb I)(\mathbb F-\mathbb I)^\top)= \Delta u + \mathcal{P}\operatorname{div}(\mathbb F-\mathbb I).
\end{align*}
Noticing that $\operatorname{div}\mathbb F^{\top}=0$ gives
\begin{align*}			
\mathcal{P}\operatorname{div}(\mathbb F-\mathbb I)=&\operatorname{div}(\mathbb F-\mathbb I)+\nabla (-\Delta )^{-1}\operatorname{div}\operatorname{div}(\mathbb F-\mathbb I)
=\operatorname{div}(\mathbb F-\mathbb I).
\end{align*}
Hence
\begin{align}\label{hou}
\operatorname{div}(\mathbb F\mathbb F^{\top})-\nabla P
&=u\cdot\nabla u+u_t-\Delta u\nonumber\\
&=\mathcal{O}(u\cdot\nabla u)+\mathcal{P}\dot{u}-\Delta u\nonumber\\
&=\mathcal{O}(u\cdot\nabla u)+\mathcal{P}\operatorname{div}(( \mathbb F-\mathbb I)(\mathbb F-\mathbb I)^{\top})+\operatorname{div}(\mathbb F-\mathbb I).
\end{align}
Substituting (\ref{shou})-(\ref{hou}) into (\ref{sup}) leads to
\begin{align}\label{214}
&\left.  \frac{1}{4}\int_{\mathbb{R}^{n}}|u|^4\;dx\right|_{0} ^{t}+\int_0^t\int_{\mathbb{R}^{n}}(|\nabla u|^2|u|^2+\frac{1}{2}|\nabla|u|^2|^2)\;dxdt'\nonumber\\ =&\int_0^t\int_{\mathbb{R}^{n}}[\mathcal{O}(u\cdot\nabla u)+\mathcal{P}\operatorname{div}((\mathbb F-\mathbb I)(\mathbb F-\mathbb I)^{\top})+\operatorname{div}(\mathbb F-\mathbb I)]\cdot |u|^2u\;dxdt'\nonumber\\ =&\int_0^t\int_{\mathbb{R}^{n}}[-\na(-\Delta)^{-1}\operatorname{div}(u\cdot\na u)\cdot |u|^2u+\operatorname{div}((\mathbb F-\mathbb I)(\mathbb F-\mathbb I)^{\top})\cdot|u|^2 u\nonumber\\ &+\nabla(-\Delta)^{-1}\operatorname{div}\operatorname{div}((\mathbb F-\mathbb I)(\mathbb F-\mathbb I)^{\top})\cdot |u|^2 u -(\mathbb F-\mathbb I):\nabla (|u|^2u)]\;dxdt'\nonumber\\ =&\int_0^t\int_{\mathbb{R}^{n}}[-\na(-\Delta)^{-1}\operatorname{div}(u\cdot\na u)\cdot |u|^2u-(\mathbb F-\mathbb I)(\mathbb F-\mathbb I)^{\top}:\nabla (|u|^2 u)\nonumber\\		&-(-\Delta)^{-1}\operatorname{div}\operatorname{div}((\mathbb F-\mathbb I)(\mathbb F-\mathbb I)^{\top})\nabla |u|^2\cdot u-(\mathbb F-\mathbb I):\nabla (|u|^2 u)]\;dxdt'.
\end{align}
Owing to Young's inequality, \eqref{12}, \eqref{nl} and the following fact that $|\nabla (|u|^2u)|\leq C|\nabla u ||u|^2$, the right-hand side of \eqref{214} is controlled by
\begin{align*}
\frac{1}{2}\int_0^t\int_{\mathbb{R}^{n}}|u|^2|\nabla u|^2\;dxdt'+C\int_0^t\int_{\mathbb{R}^{n}}| u|^6\;dxdt'+C\int_0^tA_2(t')dt'+C\varepsilon_{0}^\alpha.
\end{align*}
Using \eqref{GN} and H\"{o}lder's inequality, one arrives at
\begin{align*}
\int_0^t\int_{\mathbb{R}^{2}}| u|^6\;dxdt'\leq C\int_0^t\left(\int_{\mathbb{R}^{2}}| u|^4\;dx\right)\left(\int_{\mathbb{R}^{2}}| \na u|^2\;dx\right)dt'\leq C\varepsilon_{0}A_2(T)\leq C\varepsilon_{0}^\alpha+CA_2(T)^\frac32,
\end{align*}
and
\begin{align*}
\int_0^t\int_{\mathbb{R}^{3}}| u|^6\;dxdt'\leq C\int_0^t\left(\int_{\mathbb{R}^{3}}| u|^4\;dx\right)^\frac34\left(\int_{\mathbb{R}^{3}}| \na|u|^2|^2\;dx\right)^\frac34dt'\leq CA_2(T)^\frac32.
\end{align*}

Also, noticing that H\"{o}lder's inequality gives
\begin{equation*}
\int_{\mathbb{R}^{n}}|u_0|^4\;dx\leq \left(\int_{\mathbb{R}^{n}}|u_0|^2\;dx\right)^\frac{q-4}{q-2}\left(\int_{\mathbb{R}^{n}}|u_0|^q\;dx\right)^\frac
{2}{q-2}\leq C\varepsilon^\alpha.
\end{equation*}
We thus obtain Lemma \ref{2.3} directly after substituting the above estimates back in \eqref{214} and using Gronwall's inequality.
\end{proof}

We now combine the estimates in Lemma \ref{2.2} and Lemma \ref{2.3},  together with the smallness hypothesis of $\varepsilon$ and the fact that $A_1(T)$ and $A_2(T)$ are continuous in $t$ to conclude
\begin{equation}\label{a12}
A_1(T)+A_2(T)\leq C\varepsilon_{0}^\alpha.
\end{equation}

Next, the following estimate on $A_3(T)$ will play important roles in the estimate of $B(T)$.
\begin{Lemma}\label{A3}
Assume that the hypotheses and notations of Proposition \ref{gj} are in force. Then
\begin{equation*}
A_3(T)\leq M.
\end{equation*}
\end{Lemma}
\begin{proof}
Split $A_3(T)$ into $A_{31}(T)$ and $A_{32}(T)$ with
\begin{align*}
A_{31}(T)=\sup\limits_{0\leq t\leq T}\int_{\mathbb{R}^n}\sigma(t)^{1-s}|\nabla u(x, t)|^2\;dx+\int_0^T\int_{\mathbb{R}^n}\sigma(t)^{1-s}|\mathcal{P}\dot{u}(x, t)|^2\;dxdt,\\
A_{32}(T)=\sup\limits_{0\leq t\leq T}\int_{\mathbb{R}^n}\sigma(t)^{2-s}|\mathcal P\dot u(x, t)|^2\;dx+\int_0^T\int_{\mathbb{R}^n}\sigma(t)^{2-s}|\na\mathcal{P}\dot{u}(x, t)|^2\;dxdt.
\end{align*}
In fact, due to \eqref{a12}, we only need to prove
\begin{equation*}
A_{31}(\sigma (T))+ A_{32}(\sigma (T))\leq M.
\end{equation*}
Without loss of generality, assume $T\leq 1$ in the rest of the proof of this Lemma.

We first show the estimate for $A_{31}(\sigma (T))$, in which we follow the interpolation argument as given in Hoff \cite{Hoff2002}. To begin with, we define a differential operator $\mathcal L$ acting on functions $w:\mathbb R^n\times\mathbb R^+\to\mathbb R^n$ with $\diver w=0$
by
\begin{equation*}
\mathcal L w=w_t+ \mathcal P(u\cdot \na w)-\Delta w.
\end{equation*}
We define $w_1$ and $w_2$ by
\begin{equation}\label{w1}
\mathcal Lw_1=0\quad\text{with}\quad w_1(x,0)=w_{10}(x),\;\diver w_{10}=0,
\end{equation}
and
\begin{equation}\label{w2}
\mathcal Lw_2=\mathcal P\diver(\mathbb F\mathbb F^\top)\quad\text{with}\quad w_2(x,0)=0
\end{equation}
for a given $w_{10}$. Notice that if $w_{10}=u_0$, then $w_1+w_2=u$.
\\
Multiplying \eqref{w1} and \eqref{w2} by $w_1$ and $w_2$, respectively, applying standard energy estimate, we immediately obtain
\begin{gather}
\sup\limits_{0\leq t'\leq t}\int_{\mathbb{R}^n}|w_1|^2\;dx+\int_0^t\int_{\mathbb{R}^n}|\na w_1|^2\;dxdt'\leq \int_{\mathbb{R}^n}| w_{10}|^2\;dx,\label{2.24}\\
\sup\limits_{0\leq t'\leq t}\int_{\mathbb{R}^n}|w_2|^2\;dx+\int_0^t\int_{\mathbb{R}^n}|\na w_2|^2\;dxdt'\leq C\int_0^t\int_{\mathbb{R}^n}|\mathbb F \mathbb F^\top-\mathbb I|^2\;dxdt'\leq  C\varepsilon_{0}^\alpha\label{2.25}
\end{gather}
for $t\in[0,T]$, where we have used \eqref{12} and  \eqref{nl}.
\\
When $n=2$, multiplying \eqref{w1} by $ \mathcal P\dot w_{1}$, integrating over $\mathbb R^2\times [0,t]$ by parts for $t\in[0,T]$,  using \eqref{GN}, \eqref{a12} and H\"{o}lder's inequality, one derives at
\begin{align}\label{364}
&\int_0^t\int_{\mathbb{R}^2}|\mathcal P\dot w_1|^2\;dxdt'+\frac12\int_{\mathbb{R}^2}|\na w_1|^2\;dx\nonumber\\
=&\frac12\int_{\mathbb{R}^2}|\na w_{10}|^2\;dx+\int_0^t\int_{\mathbb{R}^2}\Delta w_1\cdot (u\cdot \na w_1)\;dxdt'\nonumber\\
\leq &\frac12\int_{\mathbb{R}^2}|\na w_{10}|^2\;dx+C\int_0^t\int_{\mathbb{R}^2}|\na u||\na w_1|^2\;dxdt'\nonumber\\
\leq&\frac12\int_{\mathbb{R}^2}|\na w_{10}|^2\;dx+C\varepsilon_{0}^\alpha\left(\int_0^t\left(\int_{\mathbb{R}^2}|\mathcal P\dot w_1|^2\;dx\right)\left(\int_{\mathbb{R}^2}|\na w_1|^2\;dx\right)dt'\right)^\frac12.
\end{align}
When $n=3$, multiplying \eqref{w1} by $ w_{1t}$, integrating over $\mathbb R^3\times [0,t]$ by parts for $t\in[0,T]$, using \eqref{GN}, \eqref{a12} and H\"{o}lder's inequality, one derives at
\begin{align}\label{2.26}
&\int_0^t\int_{\mathbb{R}^3}|\mathcal P\dot w_1|^2\;dxdt'+\frac12\int_{\mathbb{R}^3}|\na w_1|^2\;dx\nonumber\\
=&\frac12\int_{\mathbb{R}^3}|\na w_{10}|^2\;dx+\int_0^t\int_{\mathbb{R}^3}\mathcal P\dot w_1\cdot\mathcal P(u\cdot \na w_1)\;dxdt'\nonumber\\
\leq&\frac12\int_{\mathbb{R}^3}|\na w_{10}|^2\;dx+\int_0^t\left(\int_{\mathbb{R}^3}|\mathcal P\dot w_1|^2\;dx\right)^\frac12\left(\int_{\mathbb{R}^3}|u|^3\;dx\right)^\frac13\left(\int_{\mathbb{R}^3}|\na w_1|^6\;dx\right)^\frac16dt'\nonumber\\
\leq&\frac12\int_{\mathbb{R}^3}|\na w_{10}|^2\;dx+C\varepsilon_{0}^\alpha\int_0^t\int_{\mathbb{R}^3}|\mathcal P\dot w_1|^2\;dxdt'.
\end{align}
Similarly, multiplying \eqref{w1} by $\sigma \mathcal P\dot w_{1}$ when $n=2$ and $\sigma w_{1t}$ when $n=3$, integrating the resulting equalities over $\mathbb R^n\times [0,t]$ for $t\in[0,T]$ and using \eqref{GN}, \eqref{a12},  \eqref{2.24} and H\"{o}lder's inequality, we have
\begin{align}\label{2.27}
\int_0^t\int_{\mathbb{R}^n}\sigma|\mathcal P\dot w_1|^2\;dxdt'+\frac12\int_{\mathbb{R}^n}\sigma|\na w_1|^2\;dx
\leq&C\int_{\mathbb{R}^n}| w_{10}|^2\;dx+C\varepsilon_{0}^\alpha\int_0^t\int_{\mathbb{R}^n}\sigma|\mathcal P\dot w_1|^2\;dxdt'\nonumber\\
&+C\varepsilon_{0}^\alpha\left(\int_0^t\left(\int_{\mathbb{R}^n}\sigma|\mathcal P\dot w_1|^2\;dx\right)\left(\int_{\mathbb{R}^n}\sigma|\na w_1|^2\;dx\right)dt'\right)^\frac12.
\end{align}
The second term on the right-hand side of  \eqref{364} and \eqref{2.26}, the second and third term on the right-hand side of \eqref{2.27} can be absorbed by the left-hand side respectively due to the smallness of $\varepsilon$. Then it follows from \eqref{364}, \eqref{2.26} and \eqref{2.27} that
\begin{gather*}
\sup\limits_{0\leq t\leq T}\int_{\mathbb{R}^n}|\nabla w_1|^2dx+\int_{0}^T\int_{\mathbb{R}^n}|\mathcal P\dot w_1|^2dxdt\leq C\|w_{10}\|_{H^1(\mathbb R^n)}^2,
\\\sup\limits_{0\leq t\leq T}\left(t\int_{\mathbb{R}^n}|\nabla w_1|^2dx\right)+\int_{0}^T\int_{\mathbb{R}^n}t|\mathcal P\dot w_1|^2dxdt\leq C\|w_{10}\|_{L^2(\mathbb R^n)}^2.
\end{gather*}
Since the operator $\mathcal L$ is linear, combining the Riesz-Thorin interpolation argument (see Grafakos \cite{Grafakos2014}) with the above two inequalities yields
\begin{equation}\label{2.29}
\sup\limits_{0\leq t\leq T}\left(t^{1-s}\int_{\mathbb{R}^n}|\nabla w_1|^2dx\right)+\int_{0}^T\int_{\mathbb{R}^n}t^{1-s}|\mathcal P\dot w_1|^2dxdt\leq C \|w_{10}\|_{H^s(\mathbb R^n)}^2.
\end{equation}

Now we turn to the estimate for $w_2$.
We define the effective viscous flux $\mathcal W$ about $w_2$ that
\begin{equation*}
\mathcal{W}\triangleq\nabla w_2- (-\Delta)^{-1}\nabla{\mathcal{P}}\operatorname{div}(\mathbb{F}\mathbb{F}^{\top}-\mathbb I).
\end{equation*}
As $\mathcal {G}_1$ in \eqref{ck}, we have following estimates for $\mathcal W$ :
\begin{equation}\label{2.31}
\|\Delta\mathcal{W}\|_{L^2(\mathbb{R}^{n})}^2\leq C\|\nabla \mathcal{P}\dot{w}_2\|_{L^2(\mathbb{R}^{n})}^2, \quad \|\nabla\mathcal{W}\|_{L^2(\mathbb{R}^{n})}^2\leq C   \|\mathcal{P}\dot{w}_2\|_{L^2(\mathbb{R}^{n})}^2.
\end{equation}
Multiplying \eqref{w2} by $\mathcal P\dot w_2$, integrating the resulting equation by parts and using \eqref{a12}, \eqref{2.25}, H\"{o}lder's inequality, Young's inequality and  the boundedness of Riesz transforms in $L^p(\mathbb R^n)(1< p<\infty)$, then it follows in a similar way to \eqref{318} that
\begin{align}\label{318'}
\sup_{0\leq t\leq T}\int_{\mathbb{R}^n}|\nabla w_2|^2\;dx+\int_{0}^{T}\int_{\mathbb{R}^n}|\mathcal{P}\dot{w_2}|^2\;dxdt\leq C\left(\varepsilon_{0}+\sum\left|\int_0^T\int_{\mathbb{R}^{n}} u_i^l{w_2}_j^m{w_2}_k^r\;dxdt\right|\right).
\end{align}
Define
\begin{equation*}
D(T)=\sup_{0\leq t\leq T}\int_{\mathbb{R}^n}|\nabla w_2|^2\;dx+\int_{0}^{T}\int_{\mathbb{R}^n}|\mathcal{P}\dot{w_2}|^2\;dxdt.
\end{equation*}

With the help of \eqref{GN}, \eqref{a12}, \eqref{2.25}, \eqref{2.31}, H\"{o}lder's inequality, Cauchy's inequality and the boundedness of Riesz transforms in $L^p(\mathbb R^n)(1< p<\infty)$, the second term on the right-hand side of \eqref{318'} can be bounded as
\begin{align*}
C\left|\int_0^T\int_{\mathbb{R}^{2}} u_i^l{w_2}_j^m{w_2}_k^r\;dxdt\right|\nonumber
\leq&C\int_0^{T}\int_{\mathbb{R}^{2}}|\na u|^2\;dxdt +\epsilon\int_0^{T}\int_{\mathbb{R}^{2}}|\na w_2|^4\;dxdt\\
\leq&C\varepsilon_0^\alpha+\epsilon\int_0^{T}\left(\int_{\mathbb{R}^{2}}|\mathcal W|^2\;dx\right)\left(\int_{\mathbb{R}^{2}}|\na\mathcal W|^2\;dx\right)dt\nonumber\\
\leq&C\varepsilon_0^\alpha+\frac12D(T)+\frac12D(T)^2
\end{align*}
when $n=2$, and
\begin{align*}
&C\left|\int_0^T\int_{\mathbb{R}^{3}} u_i^l{w_2}_j^m{w_2}_k^r\;dxdt\right|\nonumber\\
=&C\left|\int_0^{T}\int_{\mathbb{R}^{3}} u_i^l(\mathcal{W}^{mj}+((-\Delta)^{-1}\nabla{\mathcal{P}}\operatorname{div}(\mathbb{F}\mathbb{F}^{\top}-\mathbb I))^{mj})(\mathcal{W}^{rk}+((-\Delta)^{-1}\nabla{\mathcal{P}}\operatorname{div}(\mathbb{F}\mathbb{F}^{\top}-\mathbb I))^{rk})\;dxdt\right|\nonumber\\
\leq &\nonumber C\int_0^{T}\int_{\mathbb{R}^{3}}|u||\mathcal W||\na\mathcal W|\;dxdt
+C\int_0^{T}\int_{\mathbb{R}^{3}}|\nabla u||(-\Delta)^{-1}\nabla{\mathcal{P}}\operatorname{div}(\mathbb{F}\mathbb{F}^{\top}-\mathbb I)|^2\;dxdt\nonumber\\
&+C\int_0^{T}\int_{\mathbb{R}^{3}}|\nabla u||\mathcal W||(-\Delta)^{-1}\nabla{\mathcal{P}}\operatorname{div}(\mathbb{F}\mathbb{F}^{\top}-\mathbb I)|\;dxdt\nonumber\\
\leq& C\varepsilon_0^\alpha+C\int_0^T\left(\int_{\mathbb{R}^3}|\na \mathcal W|^2\;dx\right)^\frac12\left(\int_{\mathbb{R}^3}|\mathcal W|^6\;dx\right)^\frac16\left(\int_{\mathbb{R}^3}| u|^3\;dx\right)^\frac13\,dt\nonumber\\
&+C\int_0^T\left(\int_{\mathbb{R}^3}|\na u|^2\;dx\right)^\frac12\left(\int_{\mathbb{R}^3}|\mathcal W|^6\;dx\right)^\frac16\left(\int_{\mathbb{R}^3}| \mathbb{F}\mathbb{F}^{\top}-\mathbb I|^3\;dx\right)^\frac13\,dt\nonumber\\
\leq &C\varepsilon_0^\alpha+\frac12 D(T)
\end{align*}
when $n=3$.
Substituting the above two inequalities into \eqref{318'}, we obtain
\begin{equation}\label{2.40}
\sup\limits_{0\leq t\leq T}\int_{\mathbb{R}^n}|\nabla w_2|^2dx+\int_{0}^T\int_{\mathbb{R}^n}|\mathcal P\dot w_2|^2dxdt\leq C\varepsilon_{0}^\alpha.
\end{equation}
Setting $w_{10}=u_0$, such that $w_1+w_2=u$, we then conclude from \eqref{2.29} and \eqref{2.40} that
\begin{equation}\label{A31}
A_{31}(\sigma (T))\leq M.
\end{equation}

To complete the proof of this lemma, we therefore need to bound $A_{32}(\sigma (T))$.
Multiplying \eqref{jifen} by $\sigma^{2-s}\mathcal P \dot u$, and integrating the resulting equality over $\mathbb{R}^n\times [0, t]$ for $t\in[0,T]$ , one obtains that
\begin{equation}\label{2.43'}
\sup\limits_{0\leq t\leq T}\int_{\mathbb{R}^{n}}\sigma^{2-s}|\mathcal{P}\dot {u}|^2\;dx+	\int_0^T\int_{\mathbb{R}^{n}}\sigma^{2-s}|\nabla\mathcal{P}\dot{u}|^2\;dxdt
\leq C\varepsilon_{0}+C\left(\int_0^T\int_{\mathbb{R}^{n}}\sigma^{2-s}|\nabla u|^4\;dxdt\right).
\end{equation}
In view of \eqref{GN}, \eqref{ck}, \eqref{a12}, \eqref{A31}, Young's inequality and the boundedness of Riesz transforms in $L^p(\mathbb R^n)(1< p<\infty)$, we deduce that the second term on the right-hand side of \eqref{2.43'} is bounded by
\begin{align}\label{2.43}
 & C\int_0^T\int_{\mathbb{R}^{2}}\sigma^{2-s}(|\mathcal {G}_1|^4+|\mathbb F{\mathbb F}^\top-\mathbb I|^4)\;dxdt\nonumber\\
\leq &C\varepsilon_{0}^\alpha+C\left(\sup\limits_{0\leq t\leq T}\sigma\int_{\mathbb{R}^{2}}|\mathcal {G}_1|^2\;dx\right)\int_0^T\int_{\mathbb{R}^{n}}\sigma^{1-s}|\na\mathcal {G}_1|^2\;dxdt\nonumber\\
\leq&C\varepsilon_{0}^\alpha+M
\end{align}
in the case  $n=2$, and
\begin{align}\label{2.44}
 & C\int_0^T\int_{\mathbb{R}^{3}}\sigma^{2-s}(|\mathcal {G}_1|^4+|\mathbb F{\mathbb F}^\top-\mathbb I|^4)\;dxdt\nonumber\\
\leq &C\varepsilon_{0}^\alpha+C\left(\sup\limits_{0\leq t\leq T}\sigma^{1-s}\int_{\mathbb{R}^{3}}|\mathcal {G}_1|^2\;dx\right)^\frac12\left(\sup\limits_{0\leq t\leq T}\sigma^{2-s}\int_{\mathbb{R}^{3}}|\na\mathcal {G}_1|^2\;dx\right)^\frac12\int_0^T\int_{\mathbb{R}^{3}}\sigma^{1-s}|\na\mathcal {G}_1|^2\;dxdt\nonumber\\
\leq&C\varepsilon_{0}^\alpha+M+\frac12\sup\limits_{0\leq t\leq T}\int_{\mathbb{R}^{3}}\sigma^{2-s}|\mathcal{P}\dot {u}|^2\;dx
\end{align}
in the case $n=3$.
Substituting \eqref{2.43} and \eqref{2.44} into \eqref{2.43'} implies
\begin{equation*}
A_{32}(\sigma (T))\leq M.
\end{equation*}
Thus we complete the proof of Lemma \ref{A3}.
\end{proof}

To proof Lemma \ref{A1}, we only need to show that
\begin{Lemma}\label{bt}
Assume that the hypotheses and notations of Proposition \ref{gj} are in force. Then
\begin{equation*}
B(T)\leq C\varepsilon_{0}^\alpha.
\end{equation*}
\end{Lemma}
\begin{proof}
Applying the operator $\dot\Delta_k$ to \eqref{qd}, we obtain from \eqref{0} and \eqref{bony} that
\begin{align}\label{377}
	(\dot\Delta_k\mathbb {E}^{ij})_t+\dot S_{k-1}u\cdot\na\dot\Delta_k\mathbb {E}^{ij}+\dot\Delta_k\mathbb {E}^{ij}=\dot\Delta_k(\mathcal{G}\mathbb F)^{ij}+\dot\Delta_k(\mathcal{R}\mathbb F)^{ij}-\dot\Delta_k(\mathbb E\mathbb E)^{ij}+f_k,
\end{align}
where
\begin{align*}
f_k=&\diver(\dot S_{k-1}u\dot\Delta_k\mathbb {E}^{ij}-\dot\Delta_k\dot T_u\mathbb {E}^{ij})-\dot\Delta_k\diver(\dot T_{\mathbb {E}^{ij}}u+\dot R(u,\mathbb {E}^{ij}))\nonumber\\
=&\diver\left(\sum_{|k'-k|\leq4}(\dot S_{k-1}u-\dot S_{k'-1}u)\dot\Delta_k\dot\Delta_{k'}\mathbb {E}^{ij}+[\dot S_{k'-1}u,\dot\Delta_k]\dot\Delta_{k'}\mathbb {E}^{ij}\right)-\dot\Delta_k\diver(\dot T_{\mathbb {E}^{ij}}u+\dot R(u,\mathbb {E}^{ij})).
\end{align*}
Then, one deduces from Lemma \ref{TB}, Lemma \ref{ber} and Lemma \ref{jh} that
\begin{align*}
\sum_{k\in\mathbb Z}2^{\frac npk}\|f_k\|_{L^p(\mathbb R^n)}& \leq C\|\nabla{u}\|_{L^\infty(\mathbb R^n)}\|\mathbb {E}\|_{\dot{B}_{p,1}^{\frac np}(\mathbb R^n)}+C\|\dot T_{\mathbb {E}^{ij}}{u}+\dot R({u},\mathbb {E}^{ij})\|_{\dot{B}_{p,1}^{1+\frac np}(\mathbb R^n)} \nonumber\\
&\leq C\|\nabla{u}\|_{\dot{B}_{p,1}^{\frac np}(\mathbb R^n)\cap L^\infty(\mathbb R^n)}\|\mathbb E\|_{\dot{B}_{p,1}^{\frac np}(\mathbb R^n)\cap L^\infty(\mathbb R^n)} \nonumber\\
&\leq C\|\nabla{u}\|_{\dot{B}_{p,1}^{\frac np}(\mathbb R^n)}\|\mathbb {E}\|_{\dot{B}_{p,1}^{\frac np}(\mathbb R^n)}.
\end{align*}
Multiplying \eqref{377} by $2^{\frac npk}\dot\Delta_k\mathbb {E}^{ij}|\dot\Delta_k\mathbb {E}^{ij}|^{p-2}$, summing with respect to $i,j ,k$, integrating the resulting equation over $\mathbb R^n$, and using the above estimate, one has
\begin{equation}\label{380}
\frac{d}{dt}\|\mathbb E\|_{{\dot{B}_{p,1}^{{\frac{n}{p}}}(\mathbb R^n)}}+\|\mathbb E\|_{{\dot{B}_{p,1}^{{\frac{n}{p}}}(\mathbb R^n)}}
\leq  C(\|\nabla{u}\|_{\dot{B}_{p,1}^{\frac np}(\mathbb R^n)}\|\mathbb {E}\|_{\dot{B}_{p,1}^{\frac np}(\mathbb R^n)}+\|\mathcal{G}\mathbb{F}\|_{{\dot{B}_{p,1}^{{\frac{n}{p}}}(\mathbb R^n)}}
+\|\mathcal{R}\mathbb{F}\|_{{\dot{B}_{p,1}^{{\frac{n}{p}}}(\mathbb R^n)}}+\|\mathbb E\mathbb E\|_{{\dot{B}_{p,1}^{{\frac{n}{p}}}(\mathbb R^n)}}).
\end{equation}
Using \eqref{12}, Lemma \ref{TB} and the boundedness of Riesz transforms in Besov spaces $\dot B_{p,1}^{\frac np}(\mathbb R^n)(1< p<\infty)$, we have following estimates above terms on the right-hand side of \eqref{380}:
\begin{gather*}\label{381}
\|\mathbb E\mathbb E\|_{{\dot{B}_{p,1}^{{\frac{n}{p}}}(\mathbb R^n)}}
\leq C\|\mathbb E\|_{L^\infty(\mathbb R^n)}\|\mathbb E\|_{{\dot{B}_{p,1}^{{\frac{n}{p}}}(\mathbb R^n)}}
\leq C\|\mathbb E\|_{{\dot{B}_{p,1}^{{\frac{n}{p}}}(\mathbb R^n)}}^2,\nonumber\\
\|\nabla{u}\|_{\dot{B}_{p,1}^{\frac np}(\mathbb R^n)}\leq \|\mathcal G_1\|_{\dot{B}_{p,1}^{\frac np}(\mathbb R^n)}+\|(-\Delta)^{-1}\na\mathcal P\diver(\mathbb F\mathbb F^\top-\mathbb I)\|_{\dot{B}_{p,1}^{\frac np}(\mathbb R^n)}\leq C\|\mathcal G_1\|_{\dot{B}_{p,1}^{\frac np}(\mathbb R^n)}+C\|\mathbb E\|_{{\dot{B}_{p,1}^{{\frac{n}{p}}}(\mathbb R^n)}},\nonumber\\
\|\mathcal{G}\mathbb{F}\|_{{\dot{B}_{p,1}^{{\frac{n}{p}}}(\mathbb R^n)}}
+\|\mathcal{R}\mathbb{F}\|_{{\dot{B}_{p,1}^{{\frac{n}{p}}}(\mathbb R^n)}}
\leq C\|\mathcal{G}\|_{{\dot{B}_{p,1}^{{\frac{n}{p}}}(\mathbb R^n)}}
+C\|\mathcal{R}\|_{{\dot{B}_{p,1}^{{\frac{n}{p}}}(\mathbb R^n)}}
\leq C\|\mathcal G_2\|_{\dot{B}_{p,1}^{\frac np}(\mathbb R^n)}+C\|\mathbb E\|_{{\dot{B}_{p,1}^{{\frac{n}{p}}}(\mathbb R^n)}}^2.
\end{gather*}
Substituting the above inequalities into \eqref{380} gives
\begin{equation}\label{384}
\frac{d}{dt}\|\mathbb E\|_{{\dot{B}_{p,1}^{{\frac{n}{p}}}(\mathbb R^n)}}+\|\mathbb E\|_{{\dot{B}_{p,1}^{{\frac{n}{p}}}(\mathbb R^n)}}
\leq C(\|\mathcal G_1\|_{\dot{B}_{p,1}^{\frac np}(\mathbb R^n)}+\|\mathcal G_2\|_{\dot{B}_{p,1}^{\frac np}(\mathbb R^n)}).
\end{equation}
Next we give the argument when $n=3$, and a similar approach can be applied in $n=2$.
We infer from \eqref{ck} and the interpolation $(\dot H_p^{s_0},\dot H_p^{s_1})_{\alpha,q}=\dot B_{p,q}^s$ for $s=(1-\alpha)s_0+\alpha s_1,0<\alpha<1,1\leq p,q\leq\infty$ (cf. Theorem 6.3.1 in \cite{Bergh1976}) that
\begin{equation}\label{383}
\|\mathcal{G}_i\|_{\dot{B}_{p,1}^{\frac3p}(\mathbb R^3)}
\leq C\|\mathcal{G}_i\|_{L^p(\mathbb R^3)}^{\frac{p-3}p}\|\nabla\mathcal{G}_i\|_{L^p(\mathbb R^3)}^{\frac3p}
\leq C\|\mathcal{G}_i\|_{L^p(\mathbb R^3)}^{\frac{p-3}p}\|\mathcal P \dot{{u}}\|_{L^p(\mathbb R^3)}^{\frac3p}.
\end{equation}
It follows from \eqref{GN} and \eqref{ck} that
\begin{gather*}
\|\mathcal P \dot{{u}}\|_{L^p(\mathbb R^3)}
\leq C\|\mathcal P \dot{{u}}\|_{L^2(\mathbb R^3)}^{\frac{6-p}{2p}}\|\na\mathcal P \dot{{u}}\|_{L^2(\mathbb R^3)}^{\frac{3p-6}{2p}},\\
\|\mathcal{G}_i\|_{L^p(\mathbb R^3)}
\leq C\|\mathcal{G}_i\|_{L^2(\mathbb R^3)}^{\frac{6-p}{2p}}\|\na\mathcal{G}_i\|_{L^2(\mathbb R^3)}^{\frac{3p-6}{2p}}
\leq C\|\mathcal{G}_i\|_{L^2(\mathbb R^3)}^{\frac{6-p}{2p}}\|\mathcal P\dot u\|_{L^2(\mathbb R^3)}^{\frac{3p-6}{2p}},
\end{gather*}
which combined with \eqref{383} shows that
\begin{equation}\label{3.84}
\|\mathcal{G}_i\|_{\dot{B}_{p,1}^{\frac3p}(\mathbb R^3)}
\leq C\|\mathcal{G}_i\|_{L^2(\mathbb R^3)}^{\frac{-p^2+9p-18}{2p^2}}
\|\mathcal P\dot u\|_{L^2(\mathbb R^3)}^{\frac{3p^2-18p+36}{2p^2}}
\|\na\mathcal P \dot{{u}}\|_{L^2(\mathbb R^3)}^{\frac{9p-18}{2p^2}}.
\end{equation}

When $T>1$,
multiplying \eqref{384} by $\|\mathbb E\|_{{\dot{B}_{p,1}^{{\frac{3}{p}}}(\mathbb R^3)}}^{\frac{4p^2}{9p-18}-1}$ and integrating over $(1,t']$ respect to $t$ for $t'\in(1,T]$, we obtain
\begin{equation}\label{3.85}
\sup_{1\leq t'\leq T}\|\mathbb E\|_{{\dot{B}_{p,1}^{{\frac{3}{p}}}(\mathbb R^3)}}^{\frac{4p^2}{9p-18}}(t')+\int_1^{T}\|\mathbb E\|_{{\dot{B}_{p,1}^{{\frac{3}{p}}}(\mathbb R^3)}}^{\frac{4p^2}{9p-18}}(t)\;dt
\leq \|\mathbb E\|_{{\dot{B}_{p,1}^{{\frac{3}{p}}}(\mathbb R^3)}}^{\frac{4p^2}{9p-18}}(1)+C\sum_{i=1}^2\int_1^T\|\mathcal G_i\|_{{\dot{B}_{p,1}^{{\frac{3}{p}}}(\mathbb R^3)}}^{\frac{4p^2}{9p-18}}(t)\;dt.
\end{equation}
Due to \eqref{a12}  and \eqref{3.84}, the second term on the right-hand side of \eqref{3.85} can be bounded by
\begin{align*}
\int_1^T\|\mathcal G_i\|_{{\dot{B}_{p,1}^{{\frac{3}{p}}}(\mathbb R^3)}}^{\frac{4p^2}{9p-18}}(t)\;dt
\leq&\int_1^T\left(\|\na u\|_{L^2(\mathbb R^3)}^2+\|\mathbb E\|_{L^2(\mathbb R^3)}^2\right)^{\frac{-p^2+9p-18}{9p-18}}\|\mathcal P\dot u\|_{L^2(\mathbb R^3)}^{\frac{2(3p^2-18p+36)}{9p-18}}\|\na\mathcal P\dot u\|^2_{L^2(\mathbb R^3)}\;dt\\
\leq &CA(T)^\alpha\leq C\varepsilon_0^\alpha.
\end{align*}
Substituting the above inequality into \eqref{3.85}, we obtain
\begin{equation}\label{3.86}
B(T)\leq C(\varepsilon_0^\alpha+B(1)),\; T>1.
\end{equation}

When $T\in[0,1]$, integrating \eqref{384} over $(0,t']$ respect to $t$ for $t'\in(0,T]$, one gets
\begin{equation}\label{3.85'}
\sup_{0\leq t'\leq T}\|\mathbb E\|_{{\dot{B}_{p,1}^{{\frac{3}{p}}}(\mathbb R^3)}}(t')+\int_0^{T}\|\mathbb E\|_{{\dot{B}_{p,1}^{{\frac{3}{p}}}(\mathbb R^3)}}(t)\;dt
\leq \|\mathbb E_0\|_{{\dot{B}_{p,1}^{{\frac{3}{p}}}(\mathbb R^3)}}+C\sum_{i=1}^2\int_0^T\|\mathcal G_i\|_{{\dot{B}_{p,1}^{{\frac{3}{p}}}(\mathbb R^3)}}(t)\;dt.
\end{equation}
With the help of  \eqref{a12}, Lemma \ref{A3}, \eqref{3.84} and H\"{o}lder's inequality, the second term on the right-hand side of \eqref{3.85'} can be estimated by
\begin{align}\label{388}
	&\int_0^T\|\mathcal G_i\|_{{\dot{B}_{p,1}^{{\frac{3}{p}}}(\mathbb R^3)}}(t)\;dt\nonumber\\
	\leq&\int_0^T\left(t||\na u||_{L^2(\mathbb R^3)}^2+||\mathbb E||_{L^2(\mathbb R^3)}^2\right)^{\frac{-p^2+9p-18}{4p^2}}\left(t^{2-s}\|\mathcal P\dot u\|_{L^2(\mathbb R^3)}^2\;\right)^{\frac{3p^2-18p+36}{4p^2}-\eta}\left(t\|\mathcal P\dot u\|_{L^2(\mathbb R^3)}^2\;\right)^{\eta}\nonumber\\
	&\times\left(t^{2-s}\int_{\mathbb R^3}|\na\mathcal P\dot u|^2\;dx\right)^{\frac{9p-18}{4p^2}}t^{-\frac{(5-3s-4\eta+4s\eta)p^2-9p+18+s(9p-18)}{4p^2}}\;dt\nonumber\\
	\leq &C\varepsilon^\alpha\left(\int_0^Tt\|\mathcal P\dot u\|_{L^2(\mathbb R^3)}^2dt\right)^{\eta}
	\left(\int_0^Tt^{-\frac{(5-3s-4\eta+4s\eta)p^2-9p+18+s(9p-18)}{4(1-\eta)p^2-9p+18}}\;dt\right)^{\frac{4(1-\eta)p^2-9p+18}{4p^2}},
\end{align}
where $\eta=\frac{15s-8}{48s}$. Since $s\in(\frac8{15},1)$ and $p\in(3,6]$, it holds that
\begin{equation*}
	\frac{(5-3s-4\eta^2+4s\eta^2)p^2-9p+18+s(9p-18)}{4p^2-9p+18}< 1.
\end{equation*}
Hence
\begin{equation*}
\int_0^T\|\mathcal G_i\|_{{\dot{B}_{p,1}^{{\frac{3}{p}}}(\mathbb R^3)}}(t)\;dt
\leq C\varepsilon_0^\alpha.
\end{equation*}
Therefore, we infer from \eqref{3.85'} that
\begin{equation}\label{3.86'}
B(T)\leq C\varepsilon_0^\alpha,\; T\in[0,1].
\end{equation}
The combination of \eqref{3.86} and \eqref{3.86'} completes the proof of this lemma for $n=3$, while the proof for $n=2$ is similar to the case $n=3$. The proof of Lemma \ref{bt} is finished.
\end{proof}
Then, we give the estimate for $A_4(T)$ in the following lemma.
\begin{Lemma}\label{theta}
Assume that the hypotheses and notations of Proposition \ref{gj} are in force. Then
\begin{equation*}
A_4(T)\leq C\varepsilon_{0}^\alpha.
\end{equation*}
\end{Lemma}
\begin{proof}
Split $A_4(T)=A_{41}(T)+A_{42}(T)+A_{43}(T)$ with
\begin{align*}
A_{41}(T)=&\sup\limits_{0\leq t\leq T}\int_{\mathbb{R}^n}|\theta(x,t)-\tilde{\theta}|^2\;dx+
\int_0^T\int_{\mathbb{R}^n}|\na \theta(x,t)|^2\;dxdt;\\
A_{42}(T)=&\sup\limits_{0\leq t\leq T}\int_{\mathbb{R}^n}\sigma(t)^2|\na \theta(x,t)|^2\;dx+
\int_0^T\int_{\mathbb{R}^n}\sigma(t)^2 |{\dot\theta(x,t)}|^2\;dxdt;\\
A_{43}(T)=&\sup\limits_{0\leq t\leq T}\int_{\mathbb{R}^n}\sigma(t)^4|\dot \theta(x,t)|^2\;dx+
\int_0^T\int_{\mathbb{R}^n}\sigma(t)^4 |\na{\dot\theta(x,t)}|^2\;dxdt.
\end{align*}
The proof of this lemma consists of Lemma \ref{Lamma3.3}, Lemma \ref{Lemma 3.5} and Lemma \ref{2.5}.
\end{proof}
\begin{Lemma}\label{Lamma3.3}
Assume that the hypotheses and notations of Proposition \ref{gj} are in force. Then
\begin{equation*}
A_{41}(T)\leq C\varepsilon_{0}^\alpha.
\end{equation*}
\end{Lemma}
\begin{proof}
We subdivide the estimate into two cases, that is, $T\leq1$ and $T>1$.

\textit{Step} 1.\quad$T\leq1$.
Adding $\eqref{1.1'}_1$ multiplied by $u$ and $\eqref{1.1'}_2$ multiplied by $\mathbb{F}$, subtracting the resulting equation from $\eqref{1.1'}_3$ and rearranging terms, we obtain the equation for $\theta$:
\begin{equation}\label{3.5'}
\dot{\theta}-\Delta\theta=u_j^iu_i^j+|\na u |^2.
\end{equation}
Multiplying the above equation by $\theta-\tilde{\theta}$, integrating over $\mathbb{R}^n\times[0,t]$ for $0\leq t\leq T$, using \eqref{3.2'}, we obtain that
\begin{equation}\label{3.6'}
\int_{\mathbb{R}^n}|\theta-\tilde{\theta}|^2\;dx+\int_0^t\int_{\mathbb{R}^n}|\na \theta|^2\;dxdt'\leq 2\int_0^t\int_{\mathbb{R}^n}|\na u|^2(\theta-\tilde{\theta})\;dxdt'+C\varepsilon_{0}^\alpha.
\end{equation}
To estimate the first term on the right-hand side of \eqref{3.6'}, we  multiply \eqref{pzy} by $u$ to obtain
\begin{equation*}
|\na u|^2=\mathcal{P}\diver((\mathbb{F}-\mathbb I)(\mathbb{F}-\mathbb I)^\top)u+\diver(\mathbb{F}-\mathbb I)u
+\Delta(\frac12|u|^2)-u\mathcal{P}\dot{u}.
\end{equation*}
Thus the first term on the right-hand side of \eqref{3.6'} becomes
\begin{equation*}
\int_0^t\int_{\mathbb{R}^n}\left[\mathcal{P}\diver((\mathbb{F}-\mathbb I)(\mathbb{F}-\mathbb I)^\top)u+\diver(\mathbb{F}-\mathbb I)u+
\Delta(\frac12|u|^2)-u\mathcal{P}\dot{u}\right](\theta-\tilde{\theta})\;dxdt'
\triangleq\sum_{i=1}^{4}Q_{i}.
\end{equation*}
By integrating by parts and combining \eqref{3.2'}, we derive that
\begin{align}\label{370}
Q_1=&\int_0^t\int_{\mathbb{R}^n}\mathcal{P}\diver((\mathbb{F}-\mathbb I)(\mathbb{F}-\mathbb I)^\top)u(\theta-\tilde{\theta})\;dxdt'\nonumber\\
=&\int_0^t\int_{\mathbb{R}^n}(\mathbb Id+\nabla(-\Delta)^{-1}\operatorname{div})\diver((\mathbb{F}-\mathbb I)(\mathbb{F}-\mathbb I)^\top)u(\theta-\tilde{\theta})\;dxdt'\nonumber\\
\leq& \int_0^t\int_{\mathbb{R}^n}|(\mathbb{F}-\mathbb I)(\mathbb{F}-\mathbb I)^\top:(\na u(\theta-\tilde{\theta})+u\otimes \na\theta)|\;dxdt'\nonumber\\
&+\int_0^t\int_{\mathbb{R}^n}|(-\Delta)^{-1}\operatorname{div}\diver((\mathbb{F}-\mathbb I)(\mathbb{F}-\mathbb I)^\top)u\cdot\na\theta|\;dxdt'\nonumber\\
\leq & C\int_0^t\int_{\mathbb{R}^n}(|\mathbb{F}-\mathbb I|^8+|u|^4+|\na u|^2+|u|^2)\;dxdt'+\frac1{16}A_{41}(t)\nonumber\\
\leq&C\varepsilon_{0}^\alpha+\frac1{16}A_{41}(t),
\end{align}
where we have used the boundedness of the Riesz transforms and Young's inequality.
In a similar way, we have
\begin{equation*}
Q_2\leq C\varepsilon_{0}^\alpha+\frac1{16}A_{41}(t).
\end{equation*}
For $Q_3$, using \eqref{3.2'}, integration by parts and Cauchy's inequality yields
\begin{align*}
Q_3=-\int_0^t\int_{\mathbb{R}^n}\na(\frac12|u|^2)\cdot\na\theta\;dxdt'
\leq C\int_0^t\int_{\mathbb{R}^n}|u||\na u||\na\theta|\;dxdt'\leq C\varepsilon_{0}^\alpha+\frac1{16}A_{41}(t).
\end{align*}
We can split $Q_4$ as
\begin{align*}
\begin{split}
Q_4=-\int_0^t\int_{\mathbb{R}^n}u u_t(\theta-\tilde{\theta})\;dxdt'-\int_0^t\int_{\mathbb{R}^n}u(\mathbb Id+\nabla(-\Delta)^{-1}\operatorname{div})(u\cdot \na u)(\theta-\tilde{\theta})\;dxdt'\triangleq Q_{41}+Q_{42}.
\end{split}
\end{align*}
For $Q_{41}$, in terms of \eqref{GN}, \eqref{3.2'} and \eqref{3.5'}, we deduce after integration by parts and using Young's inequality that
\begin{align*}
Q_{41}=&-\int_0^t\int_{\mathbb{R}^n}(\frac12 |u|^2(\theta-\tilde{\theta}))_{t'}\;dxdt'+\int_0^t\int_{\mathbb{R}^n}\frac12|u|^2\theta_{t'}\;dxdt'\nonumber\\
\leq& \epsilon\int_{\mathbb{R}^n}|\theta-\tilde{\theta}|^2\;dx+C\varepsilon_{0}^\alpha+C\int_{\mathbb{R}^n}|u|^4\;dx
+\int_0^t\int_{\mathbb{R}^n}\frac12|u|^2(\Delta\theta+u_j^iu_i^j+|\na u |^2-u\cdot\na \theta)\;dxdt'\nonumber\\
\leq&\epsilon A_{41}(t)+C\varepsilon_{0}^\alpha+C\int_0^t\int_{\mathbb{R}^n}|u||\na u||\na \theta|\;dxdt'+C\int_0^t\int_{\mathbb{R}^n}|u|^2|\na u|^2\;dxdt'+C\int_0^t\int_{\mathbb{R}^n}|u|^6\;dxdt'\nonumber\\
\leq&\epsilon A_{41}(t)+C\varepsilon_{0}^\alpha+C\int_0^t\left(\int_{\mathbb{R}^n}|u|^4dx\right)\left(\int_{\mathbb{R}^n}|\na u|^2dx\right)dt'\nonumber\\
\leq&\frac1{32}A_{41}(t)+C\varepsilon_{0}^\alpha
\end{align*}
when $n=2$. As $n=3$, one gets the above inequality by
\begin{align}\label{276}
\int_0^t\int_{\mathbb{R}^n}|u|^6\;dxdt'\leq \int_0^t\left(\int_{\mathbb{R}^n}|u|^4\;dx\right)^\frac34\left(\int_{\mathbb{R}^n}|\na|u|^2|^2\;dx\right)^\frac34dt'
\leq C\varepsilon_{0}^\alpha.
\end{align}
Based on \eqref{3.2'}, integration by parts, Young's inequality and  the boundedness of Riesz transforms in $L^p(\mathbb R^n)(1< p<\infty)$, we obtain
\begin{align}\label{3.13'}
Q_{42}\leq&\int_0^t\int_{\mathbb{R}^n}|u||\nabla(-\Delta)^{-1}\operatorname{div}(u\cdot \na u)||\theta-\tilde{\theta}|\;dxdt'+
\int_0^t\int_{\mathbb{R}^n}|u|^2|\na u||\theta-\tilde{\theta}|\;dxdt'\nonumber\\
\leq& C\int_0^t\int_{\mathbb{R}^n}|u|^2|\na u|^2\;dxdt'+C\int_0^t\int_{\mathbb{R}^n}|u|^2|\theta-\tilde{\theta}|^2\;dxdt'\nonumber\\
\leq& C\int_0^t\int_{\mathbb{R}^n}|u|^2|\theta-\tilde{\theta}|^2\;dxdt'+C\varepsilon_{0}^\alpha.
\end{align}
It follows from  H\"{o}lder's inequality, \eqref{GN}, \eqref{3.2'} and \eqref{276} that
\begin{align*}
C\int_0^t\int_{\mathbb{R}^2}|u|^2|\theta-\tilde{\theta}|^2\;dxdt'\leq&
C\left(\sup\limits_{0\leq t'\leq t}\int_{\mathbb{R}^2}|u|^4\;dx\right)^\frac12
\left(\int_0^t\int_{\mathbb{R}^2}|\theta-\tilde{\theta}|^4\;dxdt'\right)^\frac12\nonumber\\
\leq&C\varepsilon_{0}^\alpha\left(\int_0^t\left(\int_{\mathbb{R}^2}|\theta-\tilde{\theta}|^2\;dx\right)
\left(\int_{\mathbb{R}^2}|\na\theta|^2\;dx\right)dt'\right)^\frac12
\leq\frac1{32}A_{41}(t)
\end{align*}
and
\begin{align*}
C\int_0^t\int_{\mathbb{R}^3}|u|^2|\theta-\tilde{\theta}|^2\;dxdt'
\leq&C\left(\int_0^t\int_{\mathbb{R}^3}|u|^6\;dxdt'\right)^\frac13
\left(\int_0^t\int_{\mathbb{R}^3}|\theta-\tilde{\theta}|^3\;dxdt'\right)^\frac23\nonumber\\
\leq &C\varepsilon_{0}^\alpha\left(\int_0^t\left(\int_{\mathbb{R}^3}|\theta-\tilde{\theta}|^2\;dx\right)^\frac34
\left(\int_{\mathbb{R}^3}|\na\theta|^2\;dx\right)^\frac34dt'\right)^\frac23
\leq \frac1{32}A_{41}(t).
\end{align*}
Inserting the above two inequalities into \eqref{3.13'} yields
\begin{equation*}
Q_{42}\leq \frac1{32}A_{41}(t)+C\varepsilon_{0}^\alpha.
\end{equation*}
The combination of estimates for $Q_1,Q_2,Q_3$ and $Q_4$ completes the proof of Lemma \ref{Lamma3.3} when $T\leq1$.

\textit{Step} 2.\quad$T>1$.
Define
\begin{equation*}
\mathcal S=\mathsf{log}\;\frac{\theta}{\tilde{\theta}},\quad \mathcal S_1=\frac1{\tilde{\theta}}(\theta-\tilde{\theta}),
\end{equation*}
 then it is easily deduced from \eqref{3.5'} that
\begin{equation}\label{3.14'}
\dot{\mathcal S}=\frac{\dot{\theta}}{{\theta}}=\diver\left(\frac{\na u u}{\theta}\right)+\frac{\na\theta\cdot\left(\na u u\right)}{\theta^2}+\frac{|\na u|^2}{\theta}+\diver\left(\frac{\na \theta}{\theta}\right)+\frac{|\na\theta|^2}{\theta^2}.
\end{equation}
\begin{equation}\label{3.15'}
\dot{\mathcal S}_1=\frac{\dot{\theta}}{\tilde{\theta}}=\diver\left(\frac{\na u u}{\tilde{\theta}}\right)+\frac{|\na u|^2}{\tilde{\theta}}+\diver\left(\frac{\na \theta}{\tilde{\theta}}\right).
\end{equation}
Expanding $S(\theta)$ about $\tilde{\theta}$, we get
\begin{equation*}
\mathcal S=\mathcal S(\tilde{\theta})+\frac1{\tilde{\theta}}(\theta-\tilde{\theta})-\frac1{2\theta_1^2}(\theta-\tilde{\theta})^2
=\mathcal S_1-\frac1{2\theta_1^2}(\theta-\tilde{\theta})^2,
\end{equation*}
where $\theta_1$ is between $ \theta$ and $ \tilde{\theta}$.
After subtracting \eqref{3.14'} from \eqref{3.15'} and integrating the resulting equation along the trajectory, we have
\begin{equation*}
\left.\frac{|\theta(x(t'),t')-\tilde{\theta}|^2}{2\theta_1^2}\right|_{1}^{t}
+\int_1^t \frac{|\na \theta|^2}{\theta^2}(x(t'),t')dt'=\int_1^t \left[\frac{|\na u|^2}{\tilde{\theta}}-\frac{|\na u|^2}{\theta}-\frac{\na\theta\cdot\left(\na u u\right)}{\theta^2}\right]\;(x(t'),t')dt'+\int_1^t \mathcal{O}_1\;dt'
\end{equation*}
for all $1<t\leq T$, where $\mathcal{O}_1$ denotes any term which is a perfect divergence in $x$.
Integrating the above inequality over $\mathbb{R}^{n}$, using \eqref{uc}, \eqref{3.2'} and Cauchy's inequality, yields
\begin{align*}
A_{41}(t)
\leq &C\int_{\mathbb{R}^n}|\theta(x,1)-\tilde{\theta}|^2\;dx+C \int_1^t\int_{\mathbb{R}^n}(|\na u|^2+|\na \theta||\na u|)\;dxdt'\\
\leq&C\varepsilon_{0}^\alpha+\frac12A_{41}(t).
\end{align*}
Clearly, $A_{41}(t)\leq C\varepsilon_{0}^\alpha$ for all $1<t\leq T$. The proof of Lemma \ref{Lamma3.3} is therefore completed.
\end{proof}
\begin{Lemma}\label{Lemma 3.5}
Assume that the hypotheses and notations of Proposition \ref{gj} are in force. Then
\begin{equation*}
A_{42}(T)\leq C \varepsilon_{0}^\alpha.
\end{equation*}
\end{Lemma}
\begin{proof}
 Multiplying \eqref{3.5'} by $\sigma^2\dot\theta$ and  integrating over $\mathbb{R}^n\times[0,t]$ for $t\in[0,T]$, one gets
\begin{align*}
\int_0^t\int_{\mathbb{R}^n}\sigma^2|\dot\theta|^2\;dxdt'
+&\frac12\int_0^t\int_{\mathbb{R}^n}\sigma^2|\na\theta|^2_{t'}\;dxdt'
\\&=-\int_0^t\int_{\mathbb{R}^n}\sigma^2\theta_j\theta_k u_k^j\;dxdt'+\int_0^t\int_{\mathbb{R}^n}\sigma^2(u_j^ku_k^j+|\na u|^2)\dot\theta\;dxdt'.
\end{align*}
By virtue of  \eqref{3.2'}, \eqref{u4} and Lemma \ref{Lamma3.3}, we infer from integration by parts and Cauchy's inequality that
\begin{align}\label{D2}
A_{42}(t)\leq& C\int_0^{\sigma(t)}\int_{\mathbb{R}^n}\sigma|\na\theta|^2\;dxdt'+
C\int_0^t\int_{\mathbb{R}^n}\sigma^2(|\na \theta|^\frac83+|\na u|^4)\;dxdt'\nonumber\\
\leq&C\int_0^t\int_{\mathbb{R}^n}|\na \theta|^2\;dxdt'+C\varepsilon_{0}^\alpha+C\int_0^t\int_{\mathbb{R}^n}\sigma^2|\na \theta|^\frac83\;dxdt'\nonumber\\
\leq&C\varepsilon_{0}^\alpha+C\int_0^t\int_{\mathbb{R}^n}\sigma^2|\na \theta|^\frac83\;dxdt',
\end{align}
where we have used \eqref{GN}, \eqref{3.2'}, \eqref{ck} and the boundedness of Riesz transforms in $L^p(\mathbb R^n)(1< p<\infty)$ to obtain that
\begin{align}\label{269}
&\int_0^{t}\int_{\mathbb{R}^3}\sigma^2|\na u|^4\;dxdt'\nonumber\\
\leq& C\int_0^{t}\int_{\mathbb{R}^3}\sigma^2(|\mathcal G_1|^4+|\mathbb F\mathbb F^\top-\mathbb I|^4)\;dxdt'\nonumber\\
\leq& C\int_0^{t}\sigma^2\left(\int_{\mathbb{R}^3}| \mathcal G_1|^2\;dx\right)^\frac12\left(\int_{\mathbb{R}^3}|\na \mathcal G_1|^2\;dx\right)^\frac32dt'+C\varepsilon_{0}^\alpha\nonumber\\
\leq& C\left(\sup\limits_{0\leq t'\leq t}\sigma^{1-s}\int_{\mathbb{R}^3}| \mathcal G_1|^2\;dx\right)^\frac12
\left(\sup\limits_{0\leq t'\leq t}\sigma^{2-s}\int_{\mathbb{R}^3}|\na \mathcal G_1|^2\;dx\right)^\frac12
\int_0^{t}\int_{\mathbb{R}^3}\sigma| \na \mathcal G_1|^2\;dxdt'+C\varepsilon_{0}^\alpha\nonumber\\
\leq& C\varepsilon_{0}^\alpha.
\end{align}
Applying \eqref{GN}, \eqref{3.2'}, \eqref{u4}, \eqref{269}, Lemma \ref{Lamma3.3}, H\"{o}lder's inequality and Young's inequality, the second term on the right-hand side of \eqref{D2} is bounded by
\begin{align*}
\begin{split}
&C\int_0^t \sigma^2\left(\int_{\mathbb{R}^2}|\na \theta|^2\;dx\right)\left(\int_{\mathbb{R}^2}|\na^2 \theta|^2\;dx\right)^{\frac13}dt'\\
\leq&C\left(\sup\limits_{0\leq t'\leq t}\sigma^2\int_{\mathbb{R}^2}|\na \theta|^2\;dx\right)^\frac13\left(\int_0^t\int_{\mathbb{R}^2}|\na \theta|^2\;dxdt'\right)^\frac23\left(\int_0^t\int_{\mathbb{R}^2}\sigma^4(|\dot\theta|^2+|\na u|^4)\;dxdt'\right)^\frac13\\
\leq&CA_{42}(t)^\frac13\varepsilon_{0}^\alpha\left(A_{42}(t)^\frac13+\varepsilon_{0}^\alpha\right)\\
\leq&\frac12A_{42}(t)+C\varepsilon_{0}^\alpha
\end{split}
\end{align*}
in the case $n=2$, and  \begin{align*}
\begin{split}
&C
\int_0^t \sigma^2\left(\int_{\mathbb{R}^3}|\na \theta|^2\;dx\right)^{\frac56}\left(\int_{\mathbb{R}^3}|\na^2 \theta|^2\;dx\right)^{\frac12}dt'\\
\leq&C\left(\sup\limits_{0\leq t'\leq t}\sigma^2\int_{\mathbb{R}^3}|\na \theta|^2\;dx\right)^\frac13\left(\int_0^t\int_{\mathbb{R}^3}|\na \theta|^2\;dxdt'\right)^\frac12\left(\int_0^t\int_{\mathbb{R}^3}\sigma^2(|\dot\theta|^2+|\na u|^4)\;dxdt'\right)^\frac12\\
\leq&CA_{42}(t)^\frac13\varepsilon_{0}^\alpha\left(A_{42}(t)^\frac12+\varepsilon_{0}^\alpha\right)\\
\leq&\frac12A_{42}(t)+C\varepsilon_{0}^\alpha
\end{split}
\end{align*}
in the case $n=3$.
Substituting the above two inequalities into \eqref{D2}, we arrive at the desired estimate of Lemma \ref{Lemma 3.5}.
\end{proof}
\begin{Lemma}\label{2.5}
Assume that the hypotheses and notations of Proposition \ref{gj} are in force. Then
\begin{equation*}
A_{43}(T)\leq C \varepsilon_{0}^\alpha.
\end{equation*}
\end{Lemma}
\begin{proof}
Taking the material derivative in \eqref{3.5'}, we have
\begin{equation*}
\dot{\theta}_t+u\cdot\na\dot{\theta}-\Delta\theta_t-u\cdot\na\Delta\theta=(\na u:(\na u)^\top+|\na u|^2)_t+u\cdot\na(\na u:(\na u)^\top+|\na u|^2).
\end{equation*}
Multiplying the above equation by $\sigma^4 \dot{\theta}$ and  integrating over $\mathbb{R}^n\times[0,t]$ for $t
\in [0, T]$, we then obtain that
\begin{align}\label{R}
\frac12\sigma^4\int_{\mathbb{R}^n}|\dot \theta(x,t)|^2\;dx=&2\int_0^t\int_{\mathbb{R}^n}\sigma^3\sigma_{t'} |{\dot\theta(x,t')}|^2\;dxdt'+\int_0^t\int_{\mathbb{R}^n}\sigma^4 \dot{\theta}[\Delta \theta_{t'}+u\cdot\na\Delta \theta]\;dxdt'\nonumber\\
&+\int_0^t\int_{\mathbb{R}^n}\sigma^4 \dot{\theta}[(\na u:(\na u)^\top)_{t'}+u\cdot\na(\na u:(\na u)^\top)]\;dxdt'\nonumber\\
&+\int_0^t\int_{\mathbb{R}^n}\sigma^4 \dot{\theta}[(|\na u|^\mathsf 2)_{t'}+u\cdot\na|\na u|^2]\;dxdt'\nonumber\\
\triangleq&R_1+R_2+R_3+R_4.
\end{align}
For $R_1$, in terms of Lemma \ref{Lemma 3.5},
\begin{equation*}\label{J1}
R_1\leq2\int_0^{\sigma(t)}\int_{\mathbb{R}^n}\sigma^3 |{\dot\theta(x,t')}|^2\;dxdt'\leq 2\int_0^{\sigma(t)}\int_{\mathbb{R}^n}\sigma^2 |{\dot\theta(x,t')}|^2\;dxdt'\leq C \varepsilon_{0}^\alpha.
\end{equation*}
Based on integration by parts, $R_2$ can be written as
\begin{align*}
R_2=-\int_0^t\int_{\mathbb{R}^n}\sigma^4|\na\dot{\theta}|^2\;dxdt' + \int_0^t\int_{\mathbb{R}^n}\sigma^4 \na\dot{\theta}\cdot[\na(u\cdot\na\theta)-u\Delta\theta]\;dxdt.
\end{align*}
Since the second term on the right-hand side of above can be derived as
\begin{align*}
&\int_0^t\int_{\mathbb{R}^n}\sigma^4 [\na u:(\na\theta\otimes\na\dot\theta)+\na\dot\theta\cdot(u\cdot\na)\na\theta+(u\cdot\na)\na\dot\theta
\na\theta+\na u:(\na\dot\theta\otimes\na\theta)]\;dxdt'\\
=&\int_0^t\int_{\mathbb{R}^n}\sigma^4 (u\cdot\na)(\na\dot\theta\cdot\na\theta)\;dxdt'+\int_0^t\int_{\mathbb{R}^n}\sigma^4\na u:(\na\theta\otimes\na\dot\theta+\na\dot\theta\otimes\na\theta)\;dxdt'\\
=&\int_0^t\int_{\mathbb{R}^n}\sigma^4\na u:(\na\theta\otimes\na\dot\theta+\na\dot\theta\otimes\na\theta)\;dxdt',
\end{align*}
then
\begin{equation*}\label{J_2}
R_2=-\int_0^t\int_{\mathbb{R}^n}\sigma^4|\na\dot{\theta}|^2\;dxdt' + \int_0^t\int_{\mathbb{R}^n}\sigma^4\na u:(\na\theta\otimes\na\dot\theta+\na\dot\theta\otimes\na\theta)\;dxdt'.
\end{equation*}
One deduces from \eqref{GN}, \eqref{3.2'}, \eqref{ck}, \eqref{u4}, Lemma \ref{Lemma 3.5}, \eqref{269} and the boundedness of Riesz transforms in $L^p(\mathbb R^n)(1< p<\infty)$ that the second term of $R_2$ is bounded by
\begin{align*}
&C\int_0^t\int_{\mathbb{R}^2}\sigma^4|\na u||\na\theta||\na\dot\theta|\;dxdt'\nonumber\\
\leq& \frac12\int_0^t\int_{\mathbb{R}^2}\sigma^4|\na\dot\theta|^2\;dxdt'
+C\int_0^t\int_{\mathbb{R}^2}\sigma^4|\na u|^4\;dxdt'+C\int_0^t\int_{\mathbb{R}^2}\sigma^4|\na\theta|^4\;dxdt'\nonumber\\
\leq & \frac12 \int_0^t\int_{\mathbb{R}^2}\sigma^4|\na\dot\theta|^2\;dxdt'
+C\varepsilon_0^\alpha+C\int_0^t\sigma^4\left(\int_{\mathbb{R}^2}|\na\theta|^2\;dx\right)
\left(\int_{\mathbb{R}^2}|\Delta\theta|^2\;dx\right)dt'\nonumber\\
\leq & \frac12 \int_0^t\int_{\mathbb{R}^2}\sigma^4|\na\dot\theta|^2\;dxdt'
+C\varepsilon_0^\alpha+C\left(\sup\limits_{0\leq t'\leq t}\sigma^2\int_{\mathbb{R}^2}|\na\theta|^2\;dx\right)\int_0^t
\int_{\mathbb{R}^2}\sigma^2(|\dot\theta|^2+|\na u|^4)\;dxdt'\nonumber\\
\leq & \frac12 \int_0^t\int_{\mathbb{R}^2}\sigma^4|\na\dot\theta|^2\;dxdt'
+C\varepsilon_0^\alpha,
\end{align*}
and
\begin{align*}
&C\int_0^t\int_{\mathbb{R}^3}\sigma^4|\na u||\na\theta||\na\dot\theta|\;dxdt'\nonumber\\
\leq& \frac12 \int_0^t\int_{\mathbb{R}^3}\sigma^4|\na\dot\theta|^2\;dxdt'+
C\int_0^t\sigma^4\left(\int_{\mathbb{R}^3}|\na u|^3\;dx\right)^\frac23\left(\int_{\mathbb{R}^3}|\na \theta|^6\;dx\right)^\frac13dt'
\nonumber\\
\leq&\frac12 \int_0^t\int_{\mathbb{R}^3}\sigma^4|\na\dot\theta|^2\;dxdt'+
C\varepsilon^\alpha\int_0^t\int_{\mathbb{R}^3}\sigma^4|\Delta\theta|^2\;dxdt'\nonumber\\
&+C\int_0^t\left(\sigma\int_{\mathbb{R}^3}|\mathcal G_1|^2\;dx\right)^\frac12\left(\sigma^3\int_{\mathbb{R}^3}|\na\mathcal G_1|^2\;dx\right)^\frac12\left(\int_{\mathbb{R}^3}\sigma^2|\Delta\theta|^2\;dx\right)dt'\nonumber\\
\leq &\frac12 \int_0^t\int_{\mathbb{R}^3}\sigma^4|\na\dot\theta|^2\;dxdt'
+C\varepsilon_0^\alpha.
\end{align*}
Hence,
\begin{equation*}
R_{2}\leq -\frac12 \int_0^t\int_{\mathbb{R}^n}\sigma^4|\na\dot\theta|^2\;dxdt'
+C\varepsilon_0^\alpha.
\end{equation*}
For $R_3$, we have
\begin{align*}
R_3=2\int_0^t\int_{\mathbb{R}^n}\sigma^4 \dot{\theta}\na u_t:(\na u)^\top\;dxdt'+2\int_0^t\int_{\mathbb{R}^n}\sigma^4\dot{\theta}(u\cdot\na)\na u:(\na u)^\top\;dxdt=R_{31}+R_{32}.
\end{align*}
Noticing that $\na u_t=\na\mathcal P\dot u-\na\mathcal P(u\cdot \na u)$, thus
\begin{align*}
R_{31}=2\int_0^t\int_{\mathbb{R}^n}\sigma^4 \dot{\theta}\na\mathcal P\dot u:(\na u)^\top\;dxdt'-2\int_0^t\int_{\mathbb{R}^n}\sigma^4 \dot{\theta}\na\mathcal P(u\cdot \na u):(\na u)^\top\;dxdt'.
\end{align*}
Furthermore, using $\diver u=0$, the second term of $R_{31}$ becomes
\begin{align*}
-2\int_0^t\int_{\mathbb{R}^n}\sigma^4 \dot{\theta}[\na\na(-\Delta)^{-1}(\na u:(\na u)^\top)+\na u\na u+(u\cdot\na)\na u)]:(\na u)^\top\;dxdt'.
\end{align*}
Combining \eqref{3.2'}, \eqref{u4}, \eqref{269}, Young's equality and the boundedness of Riesz transforms in $L^p(\mathbb R^n)(1< p<\infty)$, one deduces that
\begin{align*}
R_3\leq& C\int_0^t\int_{\mathbb{R}^n}\sigma^4 |\dot{\theta}|\left[|\na\mathcal P\dot u||\na u|
+|\na\na(-\Delta)^{-1}(\na u:(\na u)^\top)||\na u|+|\na u|^3\right]\;dxdt'\nonumber\\
\leq&C\int_0^t\int_{\mathbb{R}^n}(\sigma^2 |\na\mathcal P\dot u|^2+\sigma^2|\na u|^4)\;dxdt'+C\int_0^t\int_{\mathbb{R}^n}\sigma^6 |\na u|^2|\dot\theta|^2\;dxdt'\\
\leq&C\varepsilon_{0}^\alpha+CA_{43}(t)^2+\frac18A_{43}(t)\nonumber,
\end{align*}
where we have used \eqref{GN}, \eqref{ck}, \eqref{u4}, \eqref{a12} and Lemma \ref{Lemma 3.5} to get
\begin{align*}
\int_0^t\int_{\mathbb{R}^2}\sigma^6 |\na u|^2|\dot\theta|^2dxdt'
\leq&C\int_0^t\int_{\mathbb{R}^2}\sigma^2|\na u|^4\;dxdt'+C\int_0^t\left(\int_{\mathbb{R}^2}\sigma^4|\dot\theta|^2\;dx\right)
\left(\int_{\mathbb{R}^2}\sigma^4|\na\dot\theta|^2\;dx\right)dt'\\
\leq &C\varepsilon_{0}^\alpha+CA_{43}(t)^2
\end{align*}
and
\begin{align*}
\int_0^t\int_{\mathbb{R}^3}\sigma^6 |\na u|^2|\dot\theta|^2dxdt'
\leq&\int_0^t\sigma^6\left(\int_{\mathbb{R}^3}|\na u|^6\;dx\right)^\frac13\left(\int_{\mathbb{R}^3}|\dot \theta|^3\;dx\right)^\frac23dt'\\
\leq&C\int_0^t\sigma^6\left(\int_{\mathbb{R}^3}(|\mathcal G_1|^6+|\mathbb F\mathbb F^\top-\mathbb I|^6)\;dx\right)^\frac13\left(\int_{\mathbb{R}^3}|\dot\theta|^2\;dx\right)^\frac12
\left(\int_{\mathbb{R}^3}|\na\dot\theta|^2\;dx\right)^\frac12dt'\\
\leq&C\int_0^t\left(\sigma^3\int_{\mathbb{R}^3}|\mathcal P\dot u|^2\;dx\right)\left(\int_{\mathbb{R}^3}\sigma^2|\dot\theta|^2\;dx\right)^\frac12
\left(\int_{\mathbb{R}^3}\sigma^4|\na\dot\theta|^2\;dx\right)^\frac12dt'\\
&+C\varepsilon^\alpha\int_0^t\left(\int_{\mathbb{R}^3}\sigma^2|\dot\theta|^2\;dx\right)^\frac12
\left(\int_{\mathbb{R}^3}\sigma^4|\na\dot\theta|^2\;dx\right)^\frac12dt'\\
\leq &C\varepsilon_{0}^\alpha+\frac18A_{43}(t).
\end{align*}
Similar to $R_3$, we have
\begin{equation*}
R_4\leq C\varepsilon_{0}^\alpha+CA_{43}(t)^2+\frac18A_{43}(t).
\end{equation*}

Summarizing estimates for $R_1,R_2,R_3$ and $R_4$ in \eqref{R}, one gets
\begin{align*}
A_{43}(t)\leq C\varepsilon_{0}^\alpha+CA_{43}(t)^2.
\end{align*}
Finally, the continuity of $A_{43}(t)$ in $t$ and the smallness of $\varepsilon_0$ imply the lemma.
\end{proof}

Next lemma is devoted to deriving a priori pointwise bounds and H\"{o}lder continuity estimates for $u$ and $\theta$. We can finish the proof of Proposition \ref{gj} by proving
\begin{Lemma}\label{pt}
Assume that the hypotheses and notations of Proposition \ref{gj} are in force. Then
\begin{gather}
\|u\|_{L^\infty(\mathbb{R}^n\times[\tau,T])}+\langle u\rangle_{\mathbb{R}^{n} \times[\tau, T]}^{\beta, \frac{\beta}{n+2\beta}} \leq C(\tau,\beta) \varepsilon_{0}^\alpha, \quad \tau >0,\label{ug2}\\
\|\theta-\tilde{\theta}\|_{L^\infty(\mathbb{R}^n\times[\tau,T])}+\langle \theta\rangle_{\mathbb{R}^{n} \times[\tau, T]}^{\beta, \frac{\beta}{n+2\beta}}\leq C(\tau,\beta) \varepsilon_{0}^\alpha, \quad \tau >0,\label{eg2}
\end{gather}
where $\beta$ and  $C(\tau,\beta)$ are as described in Proposition \ref{gj}.
\end{Lemma}
\begin{proof}
To prove \eqref{ug2}, we first recall the fact established by Lemma 3.4  in Hu-Lin  \cite{Hu-Lin2016} that
\begin{equation}\label{n2}
\langle u(\cdot,t)\rangle^{\beta}\leq C(\beta)\left[\varepsilon_{0}^\alpha+\left(\varepsilon_{0}
+\int_{\mathbb{R}^2}|\nabla u|^2\;dx\right)^\frac{1-\beta}{2}\times
\left(\int_{\mathbb{R}^2}|\mathcal{P}\dot{u}|^2\;dx\right)^\frac{\beta}{2}\right]
\end{equation}
for all $t\in(0,T]$ and $\beta\in(0,1)$ when $n=2$.

Next, we claim the bound for H\"{o}lder seminorm of $u$ when $n=3$ as follows:
\begin{align}\label{n3}
\langle u(\cdot,t)\rangle^{\beta}\leq C(\beta)\left[\varepsilon_{0}^\alpha
+\left(\varepsilon_{0}+\int_{\mathbb{R}^{3}}|\nabla u|^2\;dx\right)^\frac{1-2\beta}{4}\times\left(\int_{\mathbb{R}^{3}}|\mathcal{P}\dot{u}|^2\;dx\right)^\frac{1+2\beta}{4}\right]
\end{align}
for all $t\in(0,T]$ and $\beta\in(0,\frac12]$.\\%
Due to Sobolev's embedding inequality \eqref{em2}, the definition of $\mathcal G_1$ and the boundedness of Riesz transforms in $L^p(\mathbb R^n)(1< p<\infty)$, for $r>3, \beta=1-\frac{3}{r}$, we have
\begin{align}\label{4-35}
\langle u(\cdot,t)\rangle^{\beta}\leq C(\beta)\|\nabla u\|_{L^r(\mathbb{R}^{3})}
\leq C(\beta)(\|\mathcal G_1\|_{L^r(\mathbb{R}^{3})}+\|\mathbb F\mathbb F^\top-\mathbb I\|_{L^r(\mathbb{R}^{3})}).
\end{align}
When $p\in(3, 6]$, in view of  \eqref{GN},  \eqref{3.2'}, (\ref{ck}) and the boundedness of Riesz transforms in $L^p(\mathbb R^n)(1< p<\infty)$, we obtain
\begin{align*}
\|\mathcal{G}_1\|_{L^r(\mathbb{R}^{3})}&\leq C\left(\int_{\mathbb{R}^{3}}|\mathcal G_1|^2\;dx\right)^\frac{6-r}{4r}\left(\int_{\mathbb{R}^{3}}|\nabla\mathcal{G}_1|^2\;dx\right)^\frac{3r-6}{4r}\nonumber\\
&\leq C\left(\varepsilon_{0}+\int_{\mathbb{R}^{3}}|\nabla u|^2\;dx\right)^\frac{6-r}{4r}\left(\int_{\mathbb{R}^{3}}|\mathcal{P}\dot{u}|^2\;dx\right)^\frac{3r-6}{4r}\nonumber\\
&\leq C\left(\varepsilon_{0}+\int_{\mathbb{R}^{3}}|\nabla u|^2\;dx\right)^\frac{1-2\beta}{4}\left(\int_{\mathbb{R}^{3}}|\mathcal{P}\dot{u}|^2\;dx\right)^\frac{1+2\beta}{4}
\end{align*}
and
\begin{align*}
\|\mathbb F\mathbb F^\top-\mathbb I\|_{L^r(\mathbb{R}^{3})}
\leq C\|\mathbb F-\mathbb I\|_{L^\infty(\mathbb{R}^3\times[0, t])}^{\frac{r-2}{r}}\left(\int_{\mathbb{R}^3}|\mathbb F-\mathbb I|^2\;dx\right)^\frac{1}{r}
\leq B(t)^\frac{r-2}{2r}\varepsilon_{0}^\frac{1}{r}
\leq\varepsilon_{0}^\alpha.
\end{align*}
Substituting the above two inequalities into (\ref{4-35}) leads to \eqref{n3}.

The combination of \eqref{3.2'}, \eqref{n2} and \eqref{n3} gives
\begin{equation}\label{3.4}
\langle u(\cdot,t)\rangle^\beta\leq C(\tau,\beta)\varepsilon_{0}^\alpha,\quad t\geq\tau>0.
\end{equation}
Fixed $\tau\in(0,T]$, for $(x,t)\in \mathbb R^n\times [\tau, T]$, one gets
\begin{align*}
\left|u(x,t)-\frac1{|B_R(x)|}\int_{B_R(x)}u(y,t)dy\right|
&\leq \frac1{|B_R(x)|}\int_{B_R(x)}\frac{\left|u(x,t)-u(y,t)\right|}{|x-y|^\beta}|x-y|^\beta\;dy\nonumber\\
&\leq C(\tau,\beta)\varepsilon_{0}^\alpha \frac1{|B_R(0)|}\int_{B_R(0)}|z|^\beta\;dz\nonumber\\
&\leq C(\tau,\beta)\varepsilon_{0}^\alpha R^\beta.
\end{align*}
Taking $R = 1$, it then infers from \eqref{3.2'} that
\begin{equation*}
\|u\|_{L^\infty(\mathbb{R}^n\times[\tau,T])}\leq C(\tau,\beta)\varepsilon_{0}^\alpha.
\end{equation*}
For all $t_1,t_2\in[\tau,T]$ with $t_1\leq t_2$, we use \eqref{3.2'} and the above two estimates to obtain that
\begin{align*}
\left|u\left(x, t_{2}\right)-u\left(x, t_{1}\right)\right|
\leq & \frac{1}{\left|B_{R}(x)\right|} \int_{t_{1}}^{t_{2}} \int_{B_{R}(x)}\left|u_{t}(y, t')\right| d y d t'+C(\tau,\beta)\varepsilon_{0}^\alpha R^{\beta} \\
\leq & C R^{-\frac n2}\left|t_{2}-t_{1}\right|^{\frac{1}{2}}\left(\int_{t_{1}}^{t_{2}} \int_{\mathbb R^n}(\left|\mathcal P\dot{u}\right|^{2}+\left|\mathcal P(u \cdot \nabla u)\right|^{2} )\;d y d t'\right)^{\frac{1}{2}}+C(\tau,\beta)\varepsilon_{0}^\alpha R^{\beta} \\
\leq & C(\tau,\beta)\varepsilon_{0}^\alpha\left(R^{-\frac n2}\left|t_{2}-t_{1}\right|^{\frac{1}{2}}+R^{\beta}\right).
\end{align*}
Choosing $R=(t_2-t_1)^{1/(n+2\beta)}$, we immediately derive at \eqref{ug2}.

Next we prove \eqref{eg2}. In fact, the remaining details in proving \eqref{eg2} are similar to \eqref{ug2} once we obtain \eqref{3.4} but with $u$ being replaced by $\theta$. To begin with, let $\tau\in(0,T]$ be fixed and $t\in[\tau,T]$.
It follows from \eqref{3.2'}, \eqref{3.5'} and Lemma \ref{2.5}
that
\begin{equation}\label{4.15}
\int_{\mathbb R^n}|\Delta\theta(\cdot,t)|^2\;dx\leq C\int_{\mathbb R^n}|\na u(\cdot,t)|^4\;dx+C\int_{\mathbb R^n}|\dot\theta(\cdot,t)|^2\;dx\leq C(\tau)\varepsilon^\alpha.
\end{equation}
For any $r>n$, Sobolev's embedding inequality \eqref{em2} gives
\begin{equation*}
\langle\theta(\cdot,t)\rangle^\beta\leq C(\beta)\|\na\theta(\cdot,t)\|_{L^r(\mathbb R^n)}
\end{equation*}
with $\beta=1-n/r$. Hence, using \eqref{GN}, \eqref{4.15} and Lemma \ref{Lemma 3.5}, we derive at
\begin{align*}
\langle\theta(\cdot,t)\rangle^\beta\leq C(\beta)(\| \na\theta(\cdot,t)\|_{L^2(\mathbb R^n)}+\|\Delta\theta(\cdot,t)\|_{L^2(\mathbb R^n)})\leq C(\tau,\beta)\varepsilon^\alpha.
\end{align*}
Then, \eqref{eg2} can be proved in exactly the same way as \eqref{ug2}.
\end{proof}

\section{Proof of Theorem \ref{The1.1}}

\

In this section, we prove Theorem 1.1 by constructing weak solutions as limits of smooth
solutions. To achieve our goals, we make use of the a priori estimates derived in Section 3.
By the approximation theory, there exists $\{(u_0^\delta, \mathbb{F}_0^\delta, \theta_0^\delta)\}_{\delta>0}$ such that, for every $\delta>0$,
$u_0^\delta,\, \mathbb F_0^\delta-\mathbb I,\, \theta_0^\delta-\tilde{\theta}
\in H^\infty(\mathbb{R}^n)$,
$\operatorname{div} u_0^\delta=0$, $\theta_0^\delta\geq 0$,
and $\mathbb F_0^\delta$ satisfies all three constraints in~(1.4).
As $\delta\to 0^+$, we have
\begin{align*}
u_0^\delta \to u_0 \text{ in } H^s \cap L^q,
  \qquad
  \mathbb F_0^\delta - \mathbb F_0 \to 0 \text{ in } L^2 \cap \dot{B}_{p,1}^{\frac{n}{p}},
  \qquad
  \theta_0^\delta - \theta_0 \to 0 \text{ in } L^2.
\end{align*}
Standard local existence results now apply to showing that there is a smooth local solution $(u^\delta,\mathbb F^\delta,\theta^\delta)$ of \eqref{1.1'} with initial data $(u^\delta_{0},\mathbb F^\delta_{0},\theta^\delta_{0})$, defined up to a positive time $t_\delta$ (see  \cite{Lin-Liu-Zhang2005, Li-Zheng2023, Hu-Wang2010}). And the hypotheses of Proposition \ref{gj} hold for $t_\delta$. Then a priori estimates in Proposition \ref{gj} are applied to showing that
\begin{gather}
A_1(t_\delta)+A_2(t_\delta)+A_4(t_\delta)+B(t_\delta)\leq C\varepsilon_{0}^\alpha,\quad A_3(t_\delta)\leq M,\label{zong}\\
\|u^\delta\|_{L^\infty(\mathbb{R}^n\times[\tau,t_\delta])}+\langle u^\delta\rangle_{\mathbb{R}^{n} \times[\tau, t_\delta]}^{\beta, \frac{\beta}{n+2\beta}} \leq C(\tau,\beta) \varepsilon_{0}^\alpha, \quad \tau \in (0,t_\delta),\label{yz1}\\
\|\theta^\delta-\tilde{\theta}\|_{L^\infty(\mathbb{R}^n\times[\tau,t_\delta])}+\langle \theta^\delta \rangle_{\mathbb{R}^{n} \times[\tau, t_\delta]}^{\beta, \frac{\beta}{n+2\beta}}\leq C(\tau,\beta) \varepsilon_{0}^\alpha, \quad \tau \in (0,t_\delta)\label{yz2},
\end{gather}
where $C$ and $C(\tau,\beta)$ are as described in Theorem \ref{The1.1}, $A_i(T)$ and $B(T)$ are now as in \eqref{A} and \eqref{B}, but with $(u,\mathbb F,\theta)$ replaced by $(u^\delta,\mathbb{F}^\delta,\theta^\delta)$.
Using the standard arguments based on the local existence results together with the estimates \eqref{zong}--\eqref{yz2}, we can obtain that the approximate smooth solutions $(u^\delta,\mathbb F^\delta,\theta^\delta)$ exist and satisfy \eqref{zong}--\eqref{yz2} for all time.

The compactness of the approximate solutions $(u^\delta,\mathbb{F}^\delta,\theta^\delta)$ can then be summarised in the following lemma.
\begin{Lemma}\label{co}
There is a sequence $\delta_k\rightarrow 0$ and functions $(u,\mathbb F,\theta)$ such that as  $k\to\infty$,
\begin{gather}
u^{\delta_k},\;\theta^{\delta_k}\to u,\;\theta\textit{ uniformly on compact sets in }\mathbb{R}^n\times(0,\infty);\label{1}\\
\na u^{\delta_k},\;\na\theta^{\delta_k}\rightharpoonup \na u,\;\na \theta \;\textit{ weakly in } L^2(\mathbb R^n\times [0,\infty));\label{2}\\
\mathbb{F}^{\delta_k}\rightharpoonup \mathbb{F} \textit{ weakly* in }L^{\infty}(\mathbb R^n\times [0,\infty))\label{5};\\
\na u^{\delta_k}(\cdot,t),\;\na\theta^{\delta_k}(\cdot,t),\;\mathbb{F}^{\delta_k}(\cdot,t)-\mathbb I \rightharpoonup \na u(\cdot,t),\;\na \theta(\cdot,t),\;\mathbb{F}(\cdot,t)-\mathbb I\textit{ weakly in } L^2(\mathbb R^n),\;t>0;\label{3}\\
u^{\delta_k}\rightharpoonup u \textit{ weakly in }L^4_{\mathsf{loc}}(\mathbb R^n\times [0,\infty));\label{6}\\
P^{\delta_k}\rightharpoonup P \textit{ weakly in }L^2_{\mathsf{loc}}(\mathbb R^n\times [0,\infty));\label{7}\\
\mathbb{F}^{\delta_k}(\cdot,t)-\mathbb I \to \mathbb{F}(\cdot,t)-\mathbb I \textit{ strongly in } L^2_{\mathsf{loc}}(\mathbb R^n),\;t>0.\label{4}
\end{gather}
\end{Lemma}
\begin{proof}
First, using the bounds in \eqref{Ho} and \eqref{Ho2} for $(u^\delta,\mathbb{F}^\delta,\theta^\delta)$, with the help of the Ascoli-Arzela theorem, we immediately obtain \eqref{1}.

Next, \eqref{2}, \eqref{5} and \eqref{3} can be derived from the uniform boundedness of the relevant norms of $(u^\delta,\mathbb{F}^\delta,\theta^\delta)$ in \eqref{ab}. \eqref{6} can be deduced immediately since
\begin{equation*}
\int_{0}^{\infty} \int_{\mathbb{R}^2}|u^\delta|^4dxdt\leq \sup\limits_{t\in[0,\infty)}\|{u}^\delta(t)\|_{L^2(\mathbb R^2)}^2\int_0^\infty\|\nabla{u}^\delta(t)\|_{L^2(\mathbb R^2)}^2 dt\leq C\varepsilon_0^\alpha
\end{equation*}
and
\begin{equation*}
\int\int_{V_1}|u^\delta|^4dxdt\leq
\int_{G}\sigma^{\frac{3(s-1)}{2}}\left(\int_{S}|{u}^\delta|^2\;dx\right)^\frac12
\left(\int_{S}\sigma^{1-s}|\na{u}^\delta|^2\;dx\right)^\frac32dt\leq C(V_1)\varepsilon_0^\alpha
\end{equation*}
for the set $V_1=S\times G\subset \subset\mathbb R^3\times [0,\infty)$, where $C(V_1)$ is a positive constant depending on $V_1$.

Then, to prove \eqref{7}, we take the divergence of $\eqref{1.1'}_1$ for $(u^\delta,\mathbb{F}^\delta,\theta^\delta)$ to find
\begin{equation}\label{ty}
-\Delta P^\delta=\diver\diver(u^\delta\otimes u^\delta)-\diver\diver(\mathbb{F}^\delta{\mathbb{F}^\delta}^\top).
\end{equation}
For any set $V_2\subset \subset\mathbb R^n\times [0,\infty)$, due to \eqref{ab}, we have
\begin{equation*}
\|u^\delta\otimes u^\delta\|_{L^2(V_2)}\leq C(V_2)\varepsilon_0^\alpha, \quad \|\mathbb{F}^\delta{\mathbb{F}^\delta}^\top\|_{L^2(V_2)}\leq C(V_2),
\end{equation*}
where $C(V_2)$ is a positive constant that depends on $V_2$.
An application of the $L^p$-estimate of elliptic system to \eqref{ty} gives $\|P^\delta\|_{L^2(V_2)}\leq C(V_2)$. Thus \eqref{7} follows.

Finally, the strong convergence of $\mathbb F$ in \eqref{4} can be found in  \cite{Hu-Lin2016}(Lemma 5.1), and we omit details here for simplicity.
\end{proof}
\textbf{Proof of Theorem \ref{The1.1}} Now we complete the proof of Theorem \ref{The1.1}. It is clear that the limiting functions $(u,\mathbb{F},\theta)$ inherit the bounds in \eqref{ab}, \eqref{Ho} and \eqref{Ho2} from the sequence $(u^{\delta_k},\mathbb{F}^{\delta_k},\theta^{\delta_k})$, while \eqref{H} follows from Lemma \ref{co}.

Next we show that the limiting function $(u,\mathbb{F},\theta+\frac12 |u|^2+\frac12 |\mathbb{F}|^2)$ is a weak solution of \eqref{1.1'}--\eqref{1.3'}.
Assume $\varphi,{\psi},  {\phi}$  are arbitrary test functions with $ \operatorname{div} \varphi=0$, then we have
\begin{equation}
\int_{\mathbb{R}^n} {u}^{\delta_k}_{0} \cdot {\varphi}(x, 0) \;d x  +\int_{0}^{\infty} \int_{\mathbb{R}^n}[{u}^{\delta_k} \cdot {\varphi}_{t}+({u}^{\delta_k} \otimes {u}^{\delta_k}-{\mathbb{F}^{\delta_k} \mathbb{F}^{\delta_k}}^{\top}): \nabla {\varphi}] \;d x d t
=\int_{0}^{\infty} \int_{\mathbb{R}^n} \nabla {u}^{\delta_k}: \nabla {\varphi}\; d x d t,\label{we1'}
\end{equation}
\begin{equation}
\int_{\mathbb{R}^n}\left(\mathbb{F}^{\delta_k}_{j}\right)_{0} \cdot {\psi}(x, 0) \;d x+\int_{0}^{\infty} \int_{\mathbb{R}^n}\left[\mathbb{F}^{\delta_k}_{j} \cdot {\psi}_{t}+\left({\mathbb{\mathbb{F}}}^{\delta_k}_{j} \otimes {u}^{\delta_k}-{u}^{\delta_k} \otimes {\mathbb{F}}^{\delta_k}_{j}\right): \nabla {\psi}\;\right] d x d t=0,\label{we2'}
\end{equation}
\begin{equation}
\label{we3'}
\int_{\mathbb{R}^n}\mathcal{E}^{\delta_k}_0 {\phi}(x, 0) \;d x+\int_{0}^{\infty} \int_{\mathbb{R}^n} (\mathcal{E}^{\delta_k}\phi_t+\mathcal{E}^{\delta_k}u^{\delta_k}\cdot\na\phi)\; dxdt=\int_{0}^{\infty} \int_{\mathbb{R}^n} (\mathbb T ^{\delta_k}u ^{\delta_k} \cdot\nabla\phi+\na\theta^{\delta_k}\cdot\na\phi) \;dxdt.
\end{equation}
We should verify that as $k\to \infty$, \eqref{we1'}, \eqref{we2'} and \eqref{we3'} converge to \eqref{we1}, \eqref{we2} and \eqref{we3}, respectively. For convenience, we only show the convergence of some terms here, while the other terms are somewhat simpler.

For \eqref{we1'}, assume $\operatorname{supp}\varphi \subset\Omega_1=S_1\times G_1,\; S_1\subset\mathbb{R}^n,\;G_1\subset[0,\infty).$ Then $S_1$ and $G_1$ are compact.  We show that
\begin{align}\label{yan1}
\begin{split} \int_{G_1}\int_{S_1}{\mathbb{F}^{\delta_k}\mathbb{F}^{\delta_k}}^{\top}: \nabla {\varphi} \;d x d t\to \int_{G_1}\int_{S_1}\mathbb{F}\mathbb{F}^{\top}: \nabla {\varphi} \;d x d t.
\end{split}
\end{align}
In fact,
\begin{align}\label{Fc}
\begin{split}
\left|\int_{G_1}\int_{S_1}({\mathbb{F}^{\delta_k}\mathbb{F}^{\delta_k}}^{\top}-\mathbb{F}\mathbb{F}^{\top}) :\nabla {\varphi} \;d x d t\right|
\leq \int_{G_1}\int_{S_1}|\mathbb{F}^{\delta_k}-\mathbb{F}||\mathbb{F}^{\top}-\mathbb I||\na\varphi|\;dxdt+\int_{G_1}\int_{S_1}|\mathbb{F}^{\delta_k}-\mathbb{F}||\na\varphi|\;dxdt
\\+\int_{G_1}\int_{S_1}|\mathbb{F}^{\delta_k}-\mathbb I|
|{\mathbb{F}^{\delta_k}}^{\top}-\mathbb{F}^{\top}||\na\varphi|\;dxdt
+\int_{G_1}\int_{S_1}|{\mathbb{F}^{\delta_k}}^{\top}-\mathbb{F}^{\top}||\na\varphi|\;dxdt.
\end{split}
\end{align}
By virtue of \eqref{ab}, the first term on the right-hand side of \eqref{Fc}
\begin{align*}
\int_{G_1}\int_{S_1}|\mathbb{F}^{\delta_k}-\mathbb{F}||\mathbb{F}^{\top}-\mathbb I||\na\varphi|\;dxdt
\leq &C\varepsilon_{0}^\alpha \int_{G_1}\left(\int_{S_1}|\mathbb{F}^{\delta_k}-\mathbb{F}|^2\;dx\right)^\frac12
\left(\int_{S_1}|\na \varphi|^2\;dx\right)^\frac12dt\\
\leq &C\varepsilon_{0}^\alpha \int_{G_1}\left(\int_{S_1}|\mathbb{F}^{\delta_k}-\mathbb{F}|^2\;dx\right)^\frac12dt.
\end{align*}
Owing to \eqref{ab} for $\mathbb F^{\delta_k}$ and the strong convergence of $\mathbb{F}^{\delta_k}-\mathbb I $ in \eqref{4}, the first term on the right-hand side of \eqref{Fc} goes to zero as $k\to \infty$ by Lebesgue's dominated convergence theorem. The other three terms on the right-hand side of \eqref{Fc} converge in a similar way. Hence, \eqref{yan1} follows.

For \eqref{we2'}, assume $\operatorname{supp} \psi \subset\Omega_2=S_2\times G_2,\; S_2\subset\mathbb{R}^n,G_2\subset[0,\infty).$ Then $S_2$ and $G_2$ are compact.  We show that
\begin{gather}\label{yan2}
\int_{G_2}\int_{S_2}\mathbb{F}^{\delta_k}_{j} \otimes {u^{\delta_k}}:\na\psi\;dxdt\to
\int_{G_2}\int_{S_2}\mathbb{F}_{j} \otimes u:\na\psi\;dxdt.
\end{gather}
Notice that
\begin{equation}\label{uf2}
\begin{split}
\int_{G_2}\int_{S_2}(\mathbb{F}^{\delta_k}_{j} \otimes {u^{\delta_k}}\;-\mathbb{F}_{j} \otimes u):\na\psi\;dxdt
=&\int_{G_2}\int_{S_2}(\mathbb{F}^{\delta_k}_{j}-\mathbb F_j)\otimes {u}:\na\psi \;dxdt\\&+\int_{G_2}\int_{S_2}\mathbb{F}^{\delta_k}_{j} \otimes ({u^{\delta_k}}-u):\na\psi\;dxdt.
\end{split}
\end{equation}
The first term on the right-hand side of \eqref{uf2} goes to zero as $k\to \infty$ in terms of $\eqref{3}$, $\na\psi u\in L^2(S_2)$ and Lebesgue's dominated convergence theorem. For the second term on the right-hand side of \eqref{uf2}, 
fixing $\tau\in(0,1)$, in view of \eqref{ab}, 
we have
\begin{align*}
&\left|\int_{G_2}\int_{S_2}\mathbb{F}^{\delta_k}_{j} \otimes ({u^{\delta_k}}-u):\na\psi\;dxdt\right|
\leq C\int_{G_2\cap\{t\geq\tau\}}\int_{S_2}| u^{\delta_k}- u|\;dxdt+C\tau,
\end{align*}
which together with the uniform convergence of $u^{\delta_k}$ in \eqref{1} and Lebesgue's dominated convergence theorem leads to \eqref{yan2} after letting $ k\to\infty$ and $\tau\to 0$.

The convergence of  terms in \eqref{we3'} is similar to \eqref{yan1} or \eqref{yan2}, so we omit the details here.

Finally, it remains to prove the large-time behavior of $(u,\mathbb F, \theta)$ in \eqref{lb}, which consists of Lemma \ref{lb1} and Lemma \ref{lb2}.
\begin{Lemma}\label{lb1}
It holds that
$$
\lim\limits_{t\to\infty}\|u(\cdot,t), \theta(\cdot,t)-\tilde{\theta}\|_{L^p(\mathbb{R}^n)}=0,\;p\in(2,\infty]. $$
\end{Lemma}
\begin{proof}
This can be done as the ones in  \cite{Hoff1997}, however, for completeness we sketch the proof here.

First, assume $p\in(2,\infty)$, define
\begin{equation*}
b^{\delta_k}(t)=\int_{\mathbb R^n}|\na u^{\delta_k}|^2\;dx.
\end{equation*}
Then by the bounds in \eqref{zong} for $u^{\delta_k}$, we have
\begin{gather*}
0\leq b^{\delta_k}(t)\leq C\varepsilon_{0}^\alpha,\quad t\geq1,\\
\int_1^\infty b^{\delta_k}(t)dt\leq C\varepsilon_{0}.
\end{gather*}
Furthermore, dropping $\delta_k$'s temporarily, we also obtain by H\"{o}lder's inequality that
\begin{align*}
\int_{1}^{\infty}\left|\frac{db^{\delta_k}}{dt}\left(t\right)\right|dt
&=2\int_{1}^{\infty}\left|\int_{\mathbb R^n} \na u_t:\na u\;dx\right|dt\\
&=2\int_{1}^{\infty}\left|\int_{\mathbb R^n}\left[\na\mathcal P\dot u-\na\mathcal P(u\cdot\na u)\right]:\na u\;dx\right|dt\\
&\leq C\left(\int_{1}^{\infty}\int_{\mathbb R^n}|\nabla u|^2\;dxdt\right)^{1/2}\left(\int_{1}^{\infty}\int_{\mathbb R^n}|\nabla \mathcal P\dot u|^2\;dxdt\right)^{1/2}+C\int_{1}^{\infty}\int_{\mathbb R^n}|\nabla u|^3\;dxdt
\\&\leq C\varepsilon_{0}^{\alpha},
\end{align*}
where we have used integration by parts to deduce
\begin{align*}
\int_{1}^{\infty}\int_{\mathbb R^n}\na\mathcal P(u\cdot\na u):\na u\;dxdt
=&\int_{1}^{\infty}\int_{\mathbb R^n} \left[\na(u\cdot\na u)+\na\na(-\Delta)^{-1}\diver (u\cdot\na u)\right]:\na u\;dxdt\\
=&\int_{1}^{\infty}\int_{\mathbb R^n} [ u_k^lu_l^ju_k^j+u^lu_{kl}^ju_k^j-\diver\na\na(-\Delta)^{-1}\diver (u\cdot\na u)\cdot u]\;dxdt\\
=&\int_{1}^{\infty}\int_{\mathbb R^n} [ u_k^lu_l^ju_k^j+\frac12u^l(u_k^j)^2_l-\diver (u\cdot\na u)\diver u]\;dxdt\\
\leq&C\int_{1}^{\infty}\int_{\mathbb R^n}|\nabla u|^3\;dxdt.
\end{align*}
Thus $b^{\delta_k}$ are uniformly bounded and have uniformly bounded variations. According to Helly's selection theorem, there is a subsequence \{$b^{\delta'}$\} of \{$b^{\delta_k}$\} and a function $b$ of
bounded variation such that $b^{\delta'}(t)\to b (t)$ a.e. on $[1,\infty)$. By Lebesgue's dominated convergence theorem, for each $T>1$,
\begin{equation*}
\int_1^Tb(t) \;dt=\lim_{\delta '\to 0}\int_1^Tb^{\delta'}(t)\; dt\leq C\varepsilon_0.
\end{equation*}
So $b\in L^1[1,\infty)$, and therefore
\begin{equation}\label{b}
\lim_{t\to\infty}b(t)=0.
\end{equation}
Since $\na u^{\delta_k}(\cdot,t)\rightharpoonup \na u(\cdot,t)$ weakly in $L^2(\mathbb R ^n)$ for all $t\geq 1$,
\begin{gather*}
\int_{\mathbb R^n}|\nabla u(x,t)|^2 dx \leqslant\lim_{\overline{\delta '\to 0}}\int_{\mathbb R^n}|\nabla u^{\delta^{\prime}}(x,t)|^2dx =b(t).
\end{gather*}
For such $t$, by \eqref{GN} and the above inequality,
\begin{gather}
\int_{\mathbb R^2}| u(x,t)|^4 \;dx\leq C\int_{\mathbb R^2}\lvert u(x,t)\rvert^2 \;dx\int_{\mathbb R^2}\lvert\nabla u(x,t)\rvert^2\; dx\leq C\varepsilon_{0}b(t),\nonumber\\
\int_{\mathbb R^3}| u(x,t)|^4 \;dx\leq C\left(\int_{\mathbb R^3}\lvert u(x,t)\rvert^2 \;dx\right)^\frac12\left(\int_{\mathbb R^3}\lvert\nabla u(x,t)\rvert^2\; dx\right)^\frac32\leq C\varepsilon_{0}^\frac12(b(t))^\frac32.\label{421}
\end{gather}
The H\"{o}lder's  continuity of $u$ together with \eqref{b} and \eqref{421} shows that $u(\cdot,t)\to 0$ in $L^4(\mathbb R^n)$ as $t\to \infty$, then the convergence in $L^p(\mathbb R^n)$ follows by interpolation from the bounds \eqref{ab} and \eqref{Ho}.

Next, assume $p=\infty$, then the Sobolev's embedding inequality \eqref{em1} leads to
\begin{align}\label{422}
\|u(\cdot,t)\|_{L^\infty(\mathbb R^n)}^2
\leq &C\left(\|u(\cdot,t)\|_{L^8(\mathbb R^n)}^2+\||u(\cdot,t)|\na u(\cdot,t)\|_{L^4(\mathbb R^n)}\right)\nonumber\\
\leq& C\left(\|u(\cdot,t)\|_{L^8(\mathbb R^n)}^2+\|u(\cdot,t)\|_{L^{12}(\mathbb R^n)}\|\na u(\cdot,t)\|_{L^6(\mathbb R^n)}\right).
\end{align}
For $t\geq 1$, it follows from the bounds in \eqref{ab}, \eqref{GN}, \eqref{ck} and the boundedness of Riesz transforms in $L^p(\mathbb R^n)(1< p<\infty)$ that
\begin{gather}
\int_{\mathbb R^2}| \na u(x,t)|^6 \;dx
\leq \int_{\mathbb R^2}(|\mathcal G_1|^6+|\mathbb F\mathbb F^\top-\mathbb I|^6)\;dx
\leq C\varepsilon^\alpha+\left(\int_{\mathbb R^2}|\mathcal G_1|^2\;dx\right)\left(\int_{\mathbb R^2}|\na\mathcal G_1|^2\;dx\right)^2\leq C\varepsilon^\alpha,\label{423}\\
\int_{\mathbb R^3}| \na u(x,t)|^6 \;dx
\leq \int_{\mathbb R^3}(|\mathcal G_1|^6+|\mathbb F\mathbb F^\top-\mathbb I|^6)\;dx
\leq C\varepsilon^\alpha+\left(\int_{\mathbb R^3}|\na\mathcal G_1|^2\;dx\right)^3\leq C\varepsilon^\alpha.\label{424}
\end{gather}
Substituting \eqref{423} and \eqref{424} into \eqref{422}, which together with the large-time behavior of $u$ in $L^p(\mathbb R^n)(p\in(2,\infty))$ implies that $u(\cdot,t)\to 0$ in $L^\infty(\mathbb R^n)$. Hence, the large-time behavior for $u$ in this lemma is completed.

Finally, the proof for $\theta$ in this lemma is similar to that of  $u$, we omit the details here. Thus  we finish the proof of this lemma.
\end{proof}
\begin{Lemma}\label{lb2}
It holds that
$$ \lim_{t\to\infty}\|\mathbb F(\cdot,t)-\mathbb I\|_{L^{p}(\mathbb{R}^n)}=0,\;p\in(2,\infty).$$
\end{Lemma}
\begin{proof}
Similar to \eqref{qd}, we have
\begin{align}\label{FT}
\frac{D}{Dt}(\mathbb F^{\delta_k}-\mathbb I)+\mathbb F^{\delta_k}-\mathbb I=\mathcal{G}^{\delta_k}\mathbb F^{\delta_k}+\mathcal{R}^{\delta_k}\mathbb F^{\delta_k}-(\mathbb F^{\delta_k}-\mathbb I)(\mathbb F^{\delta_k}-\mathbb I).
\end{align}	
Multiplying the above equation by $4(\mathbb F^{\delta_k}-\mathbb I)|\mathbb F^{\delta_k}-\mathbb I|^2$, with the help of the bounds \eqref{ab} for $\mathbb F^{\delta_k}$, we obtain
\begin{align*}
\frac{D}{Dt}|\mathbb F^{\delta_k}-\mathbb I|^4+|\mathbb F^{\delta_k}-\mathbb I|^4\leq C|\mathcal{G}^{\delta_k}|^4+C|\mathcal{R}^{\delta_k}|^4.
\end{align*}
Integrating the above inequality over $\mathbb{R}^2\times[1,\infty)$, together with  \eqref{ab} for $(u^{\delta_k},\mathbb F^{\delta_k})$, \eqref{GN} and the boundedness of Riesz transforms in $L^p(\mathbb R^n)(1< p<\infty)$ leads to
\begin{align}\label{F}
\int_{1}^{\infty}\frac{d}{dt}\int_{\mathbb R^2}|\mathbb F^{\delta_k}-\mathbb I|^4\;dxdt
+\int_{1}^{\infty}\int_{\mathbb R^2}|\mathbb F^{\delta_k}-\mathbb I|^4\;dxdt
\leq & C\int_{1}^{\infty}\int_{\mathbb R^2}|\mathcal{G}_2^{\delta_k}|^4\;dxdt\nonumber\\
\leq & C\int_{1}^{\infty}\left(\int_{\mathbb R^2}|\mathcal{G}_2^{\delta_k}|^2\;dx\right)
      \left(\int_{\mathbb R^2}|\na\mathcal{G}_2^{\delta_k}|^2\;dx\right)dt\nonumber\\
\leq & C \varepsilon^\alpha.
\end{align}
Next, it follows from \eqref{zong} that for $t\geq 1$,
\begin{align}\label{F1}
\int_{\mathbb R^2}|\mathbb F^{\delta_k}(x,t)-\mathbb I|^4\;dx
\leq \|\mathbb F^{\delta_k}(\cdot,t)-\mathbb I\|^2_{L^\infty(\mathbb R^2)}\int_{\mathbb R^2}|\mathbb F^{\delta_k}(x,t)-\mathbb I|^2\;dx
\leq  C \varepsilon^\alpha.
\end{align}
Substituting \eqref{F1} into \eqref{F} yields $\mathbb F(\cdot,t)-\mathbb I\to 0$ in $L^4(\mathbb R^2)$ as $t\to \infty$.

When $n=3$, multiplying (\ref{FT}) by $\frac{8}{3}(\mathbb F^{\delta_k}-\mathbb I)|\mathbb F^{\delta_k}-\mathbb I|^\frac{2}{3}$ leads to
\begin{align*}
\frac{D}{Dt}|\mathbb F^{\delta_k}-\mathbb I|^\frac{8}{3}+|\mathbb F^{\delta_k}-\mathbb I|^\frac83\leq C|\mathcal{G}^{\delta_k}|^\frac83+C|\mathcal{R}^{\delta_k}|^\frac83.
\end{align*}
Similar to the argument in the case $n=2$, we obtain $\mathbb F(\cdot,t)-\mathbb I\to 0$ in $L^\frac83(\mathbb R^3)$ as $t\to \infty$.
Then, combining the estimates of $A_1(T)$ and $B(T)$ in \eqref{ab}, this lemma can be obtained by applying H\"{o}lder's inequality.

\end{proof}

\appendix
\renewcommand{\appendixname}{Appendix~\Alph{section}}

\section{Revisiting to viscoelastic models}\label{appendix}
\

The viscoelastic flow has been in the public eye for a long time. Most viscoelastic fluids exhibit non-Newtonian dynamic behaviors, such as rod-climbing effects, the die-swell and the tubeless siphon, since they have complex microstructures; such examples include suspensions, molten plastics, paints, blood and so on; see \cite{Renardy-Hrusa-Nohel1987,Larson1999,Fernandez-Cara-Guillen-Ortega2002,Renardy2000} for more backgrounds. Due to Oldroyd's pioneering work \cite{Oldroyd1950,Oldroyd1958}, many partial differential equation models have been established to describe the dynamic behaviors of viscoelastic flows. The key to derive a hydrodynamical model is how to formulate the extra stress tensor $\mathbb{T}_s$ (so-called the shear stress) but the pressure $P$ (so-called the normal stress), which appears in the motion equation as
\begin{align}\label{1-20250413}
u_t+u\cdot \na u+\na P=\diver \mathbb{T}_s
\end{align}
for an incompressible fluid for example. As is well-known, for incompressible Newtonian fluids, $\mathbb{T}_s$ satisfies the Newtonian viscosity law:
\begin{align}\label{2-20250413}
\mathbb{T}_s=2\mu\mathbb{D}(u),
\end{align}
where $\mu>0$ is the shear viscosity coefficient and $\mathbb{D}(u)=\tfrac{1}{2}(\nabla u+(\nabla u)^\top)$ is the stretching tensor or the rate of strain tensor. The equality \eqref{2-20250413} is such a given constitutive relation, which assumes that the shear stress tensor $\mathbb{T}_s$ depends on the strain rate tensor $\mathbb{D}(u)$ instantaneously and linearly. However, for viscoelastic non-Newtonian fluids, most of which have memory characteristics, the constitutive relation cannot be completely determined by the strain tensor or flow field at a time instant, but also depends on the history of the strain rate. So far, many efforts have been paid to formulate the constitutive relation for viscoelastic non-Newtonian fluids and thus the derived models have differential and integral forms in the literature \cite{Renardy-Hrusa-Nohel1987,Larson1999}. The integral models can also refer to \cite{Green-Rivlin1957,Green-Rivlin-Spencer1959,Green-Rivlin1960}. Here we only focus on the differential models, in particular, one kind of which obeys the rheological equations of state \cite{Oldroyd1950}:
\begin{align}\label{3-20250413}
\mathbb{T}_s+\lambda_1\frac{D_b\mathbb{T}_s}{Dt}=2\eta\Big[\mathbb{D}(u)+\lambda_2\frac{D_b\mathbb{D}(u)}{Dt}\Big],
\end{align}
where $\eta>0$ is the material viscosity, $\lambda_1>0$ is the relaxation time and $0<\lambda_2(\leq\lambda_1)$ is the retardation time. Here the objective derivative $\frac{D_b(\cdot)}{Dt}$ (sometimes called the Oldroyd derivative \cite{Fernandez-Cara-Guillen-Ortega2002})is defined by
\begin{align}\label{4-20250413}
\frac{D_b\mathbb{T}_s}{Dt}:=\frac{\pa \mathbb{T}_s}{\pa t}+u\cdot\na \mathbb{T}_s+\mathbb{T}_s\mathbb{W}(u)-\mathbb{W}(u)\mathbb{T}_s
-b[\mathbb{D}(u)\mathbb{T}_s+\mathbb{T}_s\mathbb{D}(u)],
\end{align}
where $\mathbb{W}(u)=\tfrac{1}{2}(\na u-(\na u)^\top)$ and $b\in [-1,1]$ is a constant. Collecting \eqref{1-20250413} and \eqref{3-20250413}, we can obtain the Oldroyd-B model with four parameters $(\eta,\la_1,\la_2,b)$:
\begin{align}\label{5-20250413}
\begin{cases}
\diver u=0,\\
u_t+u\cdot\na u=-\na P+\diver \mathbb{T}_s,\\
\displaystyle\mathbb{T}_s+\lambda_1\frac{D_b\mathbb{T}_s}{Dt}=2\eta\Big[\mathbb{D}(u)+\lambda_2\frac{D_b\mathbb{D}(u)}{Dt}\Big],
\end{cases}
\end{align}
which actually a simple case of the Oldroyd-B model with eight parameters \cite{Oldroyd1958}. Note that the shear stress $\mathbb{T}_s$ contains two parts as
\begin{align}\label{6-20250413}
\mathbb{T}_s=\mathbb{T}_v+\mathbb{T}_e,
\end{align}
where $\mathbb{T}_v$ and $\mathbb{T}_e$ denote the viscous and elastic stress tensor, respectively. When taking $\mathbb{T}_v=2\mu\mathbb{D}(u)$ and $\mu=\tfrac{\eta\lambda_2}{\lambda_1}$, the system \eqref{5-20250413} can be reduced to
\begin{align}\label{7-20250413}
\begin{cases}
\diver u=0,\\
u_t+u\cdot\na u=-\na P+\mu\De u+\diver \mathbb{T}_e,\\
\displaystyle\frac{1}{\lambda_1}\mathbb{T}_e+\frac{\pa \mathbb{T}_e}{\pa t}+u\cdot\na \mathbb{T}_e+Q(\mathbb{T}_e,\na u)=2\frac{\eta}{\la_1}\big(1-\frac{\la_2}{\la_1}\big)\mathbb{D}(u),
\end{cases}
\end{align}
where the quadratic form takes
\begin{align*}
Q(\mathbb{T}_e,\na u):=\mathbb{T}_e\mathbb{W}(u)-\mathbb{W}(u)\mathbb{T}_e-b[\mathbb{D}(u)\mathbb{T}_e+\mathbb{T}_e\mathbb{D}(u)].
\end{align*}
The system \eqref{7-20250413} is still complicated. In 2005, by choosing $\mathbb{T}_e=c^2\mathbb{F}{\mathbb{F}}^\top$ in \eqref{7-20250413} that is called the left Cauchy-Green strain tensor \cite{Gurtin1981}, Lin, Liu and Zhang \cite{Lin-Liu-Zhang2005} proposed a kind of simplified Oldroyd-B model,
\begin{align}\label{8-20250413}
\begin{cases}
\diver {u}=0,\\
u_t+u\cdot\na u=-\nabla P+\mu\Delta u+c^2\diver(\mathbb{F}{\mathbb{F}}^\top),\\
\mathbb{F}_t+{u}\cdot \nabla \mathbb{F}=\nabla {u}\mathbb{F},
\end{cases}
\end{align}
where $\mathbb{F}$ is the deformation gradient. On the one hand, one observes that the system \eqref{8-20250413} is a special case of system \eqref{7-20250413}. In fact, by letting $\la_1=\infty$ and $b=1$, one can recover the system \eqref{8-20250413} from \eqref{7-20250413} as stated in \cite[Lemma 1.1]{Lin-Liu-Zhang2005}. On the other hand, the system \eqref{8-20250413} can be looked on as a kind of incompressible isothermal viscoelastic equations independently, which obeys the following energy dissipation law:
\begin{align}\label{energy}
\frac{d}{dt}\int_{\Omega}\Big(\frac12|u|^2+\frac 12c^2|\mathbb{F}|^2\Big)\,dx=-\int_{\Omega}\mu|\na u|^2\,dx.
\end{align}
The transport equation for $\mathbb{F}$ in \eqref{8-20250413} can be derived by the definition of the deformation gradient $\mathbb{F}$ and the chain rule; cf. \cite{Lin-Liu-Zhang2005,Renardy-Hrusa-Nohel1987}. Based on the energy dissipation law \eqref{energy}, the authors in \cite{Lin-Liu-Zhang2005} also derived the motion equation in the system \eqref{8-20250413} by an energetic variational approach; see also \cite{Giga-Kirshtein-Liu2018}. Such the energetic variational approach cannot be used to derive a non-isothermal model in general; see \cite{Giga-Kirshtein-Liu2018} for more details. It is natural to ask what is the form of the viscoelastic heat-conducting model in contrast to the isothermal case \eqref{8-20250413}.

Next, we are devoted to exploring an alternative method to derive the version of \eqref{yi} with heat-conducting effects. We shall show that the model derived here is thermodynamically consistent, that is to say, it satisfies the first and second laws of thermodynamics.
The key point is that one makes full use of the second law of thermodynamics to determine the elastic stress tensor.

We begin with a fundamental system for incompressible (homogeneous or inhomogeneous) viscoelastic fluids, which incorporates the incompressibility, the conservation laws of mass, momentum and energy, and the transport equations for the deformation gradient, that is,
\begin{align}\label{1.1''}
\begin{cases}
\rho_t+\diver(\rho u)=0,\\
(\rho u)_t+\diver(\rho u\otimes u)=\diver\mathbb{T},\\
\mathcal{E}_t+\diver(\mathcal{E}  u)=\diver(\mathbb{T}  u)-\diver J,\\
\mathbb{F}_t+ u\cdot \nabla \mathbb{F}=\nabla  u\mathbb{F},\\
\diver u=0, \quad (x,t)\in \mathbb{R}^n\times\mathbb{R}^+,
\end{cases}
\end{align}
where the unknown variables $\rho=\rho(x,t)$, $ u=u(x,t)$, and $\mathbb{F}=\mathbb{F}(x,t)$ denotes the density, the velocity, and the deformation gradient, respectively. The internal energy per unit mass $e=e(\rho,\theta)$ is a smooth function of the density $\rho$ and the absolute temperature $\theta$. The elastic potential energy per unit mass is $W=W(\mathbb{F})$ depending only on the deformation gradient $\mathbb{F}$. Thus the total energy per unit volume is set by $\mathcal{E}=\tfrac{1}{2}\rho|u|^2+\rho e+\rho W(\mathbb{F})$. The heat flux $J$ is always assumed to satisfy Fourier's law:
\begin{align}\label{1.2''}
J=-\kappa\nabla\theta,
\end{align}
where $\kappa>0$ is the thermal conductivity. The matrix $\mathbb{T}$ is the symmetric stress tensor for viscoelastic fluids.

For the incompressible viscoelastic fluids, the stress tensor $\mathbb{T}$ includes the pressure $P$ (Lagrange multiplier), the viscous stress tensor and the elastic stress tensor. So we can decompose $\mathbb{T}$ as
\begin{align}\label{1.3''}
\mathbb{T}=-P\mathbb{I}+\mathbb{T}_v+\mathbb{T}_e,
\end{align}
where $\mathbb{T}_v$ and $\mathbb{T}_e$ denote the viscous and elastic stress tensors, respectively, and $\mathbb{I}$ is the identity matrix. For simplicity, we assume that $\mathbb{T}_v$ satisfies Newton's law of viscosity, that is,
\begin{align}\label{1.4''}
\mathbb{T}_v=2\mu\mathbb{D}( u),\quad \mathbb{D}( u)=\tfrac{1}{2}(\nabla  u+(\nabla  u)^\top),
\end{align}
where $\mu>0$ is the shear viscosity coefficient. By \eqref{1.2''}--\eqref{1.4''} and $\diver  u=0$, the system \eqref{1.1''} is reduced to
\begin{align}\label{1.5'}
\begin{cases}
\rho_t+ u\cdot\na \rho=0,\\
\rho u_t+\rho u\cdot\na  u=-\nabla P+\mu\Delta  u+\diver \mathbb{T}_e,\\
\rho(e_t+ u\cdot\nabla e)=2\mu\mathbb{D}(u):\mathbb{D}(u)+(\mathbb{T}_e-\rho W_{\mathbb{F}}\mathbb{F}^\top):\nabla  u+\kappa\Delta\theta,\\
\mathbb{F}_t+ u\cdot \nabla \mathbb{F}=\nabla  u\mathbb{F},\\
\diver  u=0, \quad (x,t)\in \mathbb{R}^n\times\mathbb{R}^+,
\end{cases}
\end{align}
which automatically satisfies the first law of thermodynamics. Next, we determine the internal energy $e$ and the elastic stress tensor $\mathbb{T}_e$ by using the second law of thermodynamics to obtain a closed system.

We introduce the entropy $S$ per unit mass. The entropy $S$ is determined by the fundamental thermodynamic relation (so-called the Gibbs' relation) \cite{Giga-Novotny2018,Feireisl-Novotny2009}:
\begin{align}\label{1.6'}
dS=\frac{1}{\theta}\Big(de+p_{\rm th}d\Big(\frac{1}{\rho}\Big)\Big),
\end{align}
where $p_{\rm th}$ is the thermodynamic pressure. Note that for incompressible fluids
\begin{align*}
p\neq p_{\rm th},
\end{align*}
while they are equal for compressible fluids, cf. \cite{Rajagopal2015}.
Then we can establish the entropy equation in the following proposition.

\begin{Proposition}\label{prop1-230923}
The system \eqref{1.5'} leads to the entropy equation
\begin{align}\label{1.7'}
S_t+ u\cdot\na S-\diver\Big(\frac{\kappa\nabla\theta}{\theta}\Big)=\frac{2\mu\mathbb{D}(u):\mathbb{D}(u)}{\theta}+\frac{\kappa|\nabla\theta|^2}{\theta^2}
+\frac{1}{\theta}(\mathbb{T}_e-\rho W_{\sss\mathbb{F}} \mathbb{F}^T):\nabla u,\quad(x,t)\in \mathbb{R}^n\times\mathbb{R}^+.
\end{align}
\end{Proposition}
\begin{proof}
Fix any element $V\subset \mathbb{R}^n$. Multiplying $\eqref{1.5'}_1,\eqref{1.5'}_3$ by $(S-\tfrac{p_{\rm th}}{\rho\theta}),\frac{1}{\theta}$, respectively, and integrating them over $V$, by \eqref{1.6'} and $\diver u=0$, we obtain
\begin{align*}
\int_{V}\Big[S_t+ u\cdot\na S-\diver\Big(\frac{\kappa\nabla\theta}{\theta}\Big)\Big]\,dx=\int_{V}\Big[\frac{2\mu\mathbb{D}( u):\mathbb{D}( u)}{\theta}+\frac{\kappa|\nabla\theta|^2}{\theta^2}\Big]\,dx+\int_{V}\frac{1}{\theta}(\mathbb{T}_e- \rho W_{\sss\mathbb{F}} \mathbb{F}^T):\nabla  u\,dx,
\end{align*}
which implies \eqref{1.7'} due to the arbitrariness of $V$.
\end{proof}

Next, we give the relation between the internal energy $e$ and the elastic stress $\mathbb{T}_e$ by using the second law of thermodynamics.

\begin{Proposition}\label{prop2-230923}
The second law of thermodynamics implies that the elastic potential energy $W$ and elastic stress $\mathbb{T}_e$ in the system \eqref{1.5'} can be simply determined such that
\begin{align}\label{1.8'}
\int_{\mathbb{R}^n} \frac{1}{\theta}(\mathbb{T}_e-\rho W_{\mathbb{F}} \mathbb{F}^T):\nabla  u\,dx\geq0.
\end{align}
\end{Proposition}
\begin{proof}
Integrating the entropy equation \eqref{1.7'} in Proposition \ref{prop1-230923} over $\mathbb{R}^n$ by parts, we immediately obtain
\begin{align*}
\frac{d}{dt}\int_{\mathbb{R}^n}\rho S\,dx&=\int_{\mathbb{R}^n}\Big[\frac{2\mu\mathbb{D}( u):\mathbb{D}( u)}{\theta}+\frac{\kappa|\nabla\theta|^2}{\theta^2}\Big]\,dx
+\int_{\mathbb{R}^n} \frac{1}{\theta}(\mathbb{T}_e-\rho W_{\sss\mathbb{F}} \mathbb{F}^T):\nabla u\,dx\\
&:=\triangle,
\end{align*}
where $\triangle$ is the entropy production.
The total entropy of the whole system is
\begin{align*}
\int_{\mathbb{R}^n}\rho S\,dx.
\end{align*}
Then the second law of thermodynamics requires $\triangle\geq0$, which in turn tells us that it is satisfied if the undetermined elastic potential energy $W(\mathbb{F})$ and elastic stress $\mathbb{T}_e$ satisfy \eqref{1.8'}.
\end{proof}

For simplicity, let
\begin{align}\label{1.9'}
\begin{cases}
W(\mathbb{F})=\frac12 c^2|\mathbb{F}|^2 &--\ \mbox{Hookean\ linear\ elasticity},\\
\mathbb{T}_e=c^2\rho \mathbb{F}\mathbb{F}^T& --\ \mbox{left\ Cauchy-Green\ tensor},
\end{cases}
\end{align}
where $c>0$ is the speed of propagation of shear waves. It is easy to check that the integral inequality \eqref{1.8'} in Proposition \ref{prop2-230923} holds from \eqref{1.9'}.

Then, by \eqref{1.9'}, we deduce that the incompressible homogeneous (say $\rho\equiv 1$) viscoelastic heat-conducting model \eqref{1.5'} becomes
\begin{align}\label{2-230923}
\begin{cases}
 u_t+ u\cdot\na  u+\nabla P=\mu\Delta  u+c^2\diver(\mathbb{F}\mathbb{F}^T),\\
e_t+ u\cdot\nabla e=2\mu\mathbb{D}( u):\mathbb{D}( u)+\kappa\Delta\theta,\\
\mathbb{F}_t+ u\cdot \nabla\mathbb{F}=\nabla  u\mathbb{F},\\
\diver u=0, \quad (x,t)\in \mathbb{R}^n\times\mathbb{R}^+,
\end{cases}
\end{align}
which satisfies the first and second laws of thermodynamics. We can easily derive the conservative form \eqref{1.1''} from \eqref{2-230923}.

\section*{Acknowledgements}

The work of the first author was partially supported by National Key R\&D Program of China (No. 2021YFA1002900) and Guangzhou City Basic and Applied Basic Research Fund (No. 2024A04J6336). The work of the second author was partially supported by the RGC of Hong Kong-SAR (Nos. 14300917, 14300819, 14302819, 14301421), the Basic and Applied Basic Research Foundation of Guangdong Province (No. 2020131515310002), and the Key Project of National Nature Science Foundation of China (No. 12131010).

Yong Wang would like to thank The Institute of Mathematical Sciences of The Chinese University of Hong Kong for the hospitality during his visit.

\bigskip

{\bf Data Availability:} Data sharing is not applicable to this article.

\bigskip

{\bf Conflict of Interest:} The authors declare that they have no conflict of interest.


\begin{thebibliography}{10}

\bibitem{Bahouri-Chemin-Danchin2011}
H.~Bahouri, J.~Y. Chemin, and R.~Danchin.
\newblock {\em Fourier {A}nalysis and {N}onlinear {P}artial {D}ifferential
  {E}quations}, volume 343 of {\em Grundlehren der Mathematischen
  Wissenschaften [Fundamental Principles of Mathematical Sciences]}.
\newblock Springer, Heidelberg, 2011.

\bibitem{Bergh1976}
J.~Bergh and J.~L\"{o}fstr\"{o}m.
\newblock {\em Interpolation {S}paces. {A}n {I}ntroduction}.
\newblock Grundlehren der Mathematischen Wissenschaften, No. 223.
  Springer-Verlag, Berlin-New York, 1976.

\bibitem{Chen-Zhang2006}
Y.~M. Chen and P.~Zhang.
\newblock The global existence of small solutions to the incompressible
  viscoelastic fluid system in 2 and 3 space dimensions.
\newblock {\em Comm. Partial Differential Equations}, 31(10-12):1793--1810,
  2006.

\bibitem{Coifman-Lions-Meyer-Semmes1993}
R.~Coifman, P.~L. Lions, Y.~Meyer, and S.~Semmes.
\newblock Compensated compactness and {H}ardy spaces.
\newblock {\em J. Math. Pures Appl. (9)}, 72(3):247--286, 1993.

\bibitem{Evans2010}
L.~C. Evans.
\newblock {\em Partial {D}ifferential {E}quations}, volume~19 of {\em Graduate
  Studies in Mathematics}.
\newblock American Mathematical Society, Providence, RI, second edition, 2010.

\bibitem{Feireisl-Novotny2009}
E.~Feireisl and A.~Novotn\'{y}.
\newblock {\em Singular {L}imits in {T}hermodynamics of {V}iscous {F}luids}.
\newblock Advances in Mathematical Fluid Mechanics. Birkh\"{a}user Verlag,
  Basel, 2009.

\bibitem{Feireisl-Novotny-Petzeltov2001}
E.~Feireisl, A.~Novotn{\'y}, and H.~Petzeltov{\'a}.
\newblock On the existence of globally defined weak solutions to the
  {N}avier-{S}tokes equations.
\newblock {\em J. Math. Fluid Mech.}, 3(4):358--392, 2001.

\bibitem{Fernandez-Cara-Guillen-Ortega2002}
E.~Fern\'{a}ndez-Cara, F.~Guill\'{e}n, and R.~R. Ortega.
\newblock Mathematical modeling and analysis of viscoelastic fluids of the
  {O}ldroyd kind.
\newblock In {\em Handbook of {N}umerical {A}nalysis, {V}ol. {VIII}}, Handb.
  Numer. Anal., VIII, pages 543--661. North-Holland, Amsterdam, 2002.

\bibitem{Giga-Kirshtein-Liu2018}
M.~Giga, A.~Kirshtein, and C.~Liu.
\newblock Variational modeling and complex fluids.
\newblock In {\em Handbook of {M}athematical {A}nalysis in {M}echanics of
  {V}iscous {F}luids}, pages 73--113. Springer, Cham, 2018.

\bibitem{Giga-Novotny2018}
Y.~Giga and A.~Novotn\'{y}, editors.
\newblock {\em Handbook of {M}athematical {A}nalysis in {M}echanics of
  {V}iscous {F}luids}.
\newblock Springer, Cham, 2018.

\bibitem{Grafakos2014}
L.~Grafakos.
\newblock {\em Classical {F}ourier {A}nalysis}, volume 249 of {\em Graduate
  Texts in Mathematics}.
\newblock Springer, New York, third edition, 2014.

\bibitem{Green-Rivlin1957}
A.~E. Green and R.~S. Rivlin.
\newblock The mechanics of non-linear materials with memory.
\newblock {\em Arch. Rational Mech. Anal.}, 1:1--21, 1957.

\bibitem{Green-Rivlin1960}
A.~E. Green and R.~S. Rivlin.
\newblock The mechanics of non-linear materials with memory. {III}.
\newblock {\em Arch. Rational Mech. Anal.}, 4:387--404 (1960), 1960.

\bibitem{Green-Rivlin-Spencer1959}
A.~E. Green, R.~S. Rivlin, and A.~J.~M. Spencer.
\newblock The mechanics of non-linear materials with memory. {II}.
\newblock {\em Arch. Rational Mech. Anal.}, 3:82--90, 1959.

\bibitem{Gurtin1981}
M.~E Gurtin.
\newblock {\em An {I}ntroduction to {C}ontinuum {M}echanics}, volume 158.
\newblock Academic press, 1981.

\bibitem{Hoff1995_1}
D.~Hoff.
\newblock Global solutions of the {N}avier-{S}tokes equations for
  multidimensional compressible flow with discontinuous initial data.
\newblock {\em J. Differential Equations}, 120(1):215--254, 1995.

\bibitem{Hoff1995_2}
D.~Hoff.
\newblock Strong convergence to global solutions for multidimensional flows of
  compressible, viscous fluids with polytropic equations of state and
  discontinuous initial data.
\newblock {\em Arch. Rational Mech. Anal.}, 132(1):1--14, 1995.

\bibitem{Hoff1997}
D.~Hoff.
\newblock Discontinuous solutions of the {N}avier-{S}tokes equations for
  multidimensional flows of heat-conducting fluids.
\newblock {\em Arch. Rational Mech. Anal.}, 139(4):303--354, 1997.

\bibitem{Hoff2002}
D.~Hoff.
\newblock Dynamics of singularity surfaces for compressible, viscous flows in
  two space dimensions.
\newblock {\em Comm. Pure Appl. Math.}, 55(11):1365--1407, 2002.

\bibitem{Hoff2005}
D.~Hoff.
\newblock Compressible flow in a half-space with {N}avier boundary conditions.
\newblock {\em J. Math. Fluid Mech.}, 7(3):315--338, 2005.

\bibitem{Hoff2006}
D.~Hoff.
\newblock Uniqueness of weak solutions of the {N}avier-{S}tokes equations of
  multidimensional, compressible flow.
\newblock {\em SIAM J. Math. Anal.}, 37(6):1742--1760, 2006.

\bibitem{Hu2018}
X.~P. Hu.
\newblock Global existence of weak solutions to two dimensional compressible
  viscoelastic flows.
\newblock {\em J. Differential Equations}, 265(7):3130--3167, 2018.

\bibitem{Hu-Lin2016}
X.~P. Hu and F.~H. Lin.
\newblock Global solutions of two-dimensional incompressible viscoelastic flows
  with discontinuous initial data.
\newblock {\em Comm. Pure Appl. Math.}, 69(2):372--404, 2016.

\bibitem{Hu-Ou-Wang-Lu2023}
X.~P. Hu, Y.~B. Ou, D.~H. Wang, and L.~Yang.
\newblock Incompressible limit for compressible viscoelastic flows with large
  velocity.
\newblock {\em Adv. Nonlinear Anal.}, 12(1):Paper No. 20220324, 20, 2023.

\bibitem{Hu-Wang2010}
X.~P. Hu and D.~H. Wang.
\newblock Local strong solution to the compressible viscoelastic flow with
  large data.
\newblock {\em J. Differential Equations}, 249(5):1179--1198, 2010.

\bibitem{Hu-Wu2013}
X.~P. Hu and G.~C. Wu.
\newblock Global existence and optimal decay rates for three-dimensional
  compressible viscoelastic flows.
\newblock {\em SIAM J. Math. Anal.}, 45(5):2815--2833, 2013.

\bibitem{Huang-Li2018}
X.~D. Huang and J.~Li.
\newblock Global classical and weak solutions to the three-dimensional full
  compressible {N}avier-{S}tokes system with vacuum and large oscillations.
\newblock {\em Arch. Ration. Mech. Anal.}, 227(3):995--1059, 2018.

\bibitem{Huang-Li-Xin2012}
X.~D. Huang, J.~Li, and Z.~P. Xin.
\newblock Global well-posedness of classical solutions with large oscillations
  and vacuum to the three-dimensional isentropic compressible {N}avier-{S}tokes
  equations.
\newblock {\em Comm. Pure Appl. Math.}, 65(4):549--585, 2012.

\bibitem{Jiang-Zhang2001}
S.~Jiang and P.~Zhang.
\newblock On spherically symmetric solutions of the compressible isentropic
  {N}avier-{S}tokes equations.
\newblock {\em Comm. Math. Phys.}, 215(3):559--581, 2001.

\bibitem{Jiang-Zhang2003}
S.~Jiang and P.~Zhang.
\newblock Axisymmetric solutions of the 3{D} {N}avier-{S}tokes equations for
  compressible isentropic fluids.
\newblock {\em J. Math. Pures Appl. (9)}, 82(8):949--973, 2003.

\bibitem{Larson1999}
R.~G. Larson.
\newblock {\em The {S}tructure and {R}heology of {C}omplex {F}luids}.
\newblock Oxford University Press, New York, 1999.

\bibitem{Lei-Liu-Zhou2008}
Z.~Lei, C.~Liu, and Y.~Zhou.
\newblock Global solutions for incompressible viscoelastic fluids.
\newblock {\em Arch. Ration. Mech. Anal.}, 188(3):371--398, 2008.

\bibitem{Lei-Zhou2005}
Z.~Lei and Y.~Zhou.
\newblock Global existence of classical solutions for the two-dimensional
  {O}ldroyd model via the incompressible limit.
\newblock {\em SIAM J. Math. Anal.}, 37(3):797--814, 2005.

\bibitem{Li-Zheng2023}
J.~K. Li and Y.~S. Zheng.
\newblock Local existence and uniqueness of heat conductive compressible
  {N}avier-{S}tokes equations in the presence of vacuum without initial
  compatibility conditions.
\newblock {\em J. Math. Fluid Mech.}, 25(1):Paper No. 14, 42, 2023.

\bibitem{Li-Liu-Ye2022}
Z.~L. Li, H.~Liu, and Y.~L. Ye.
\newblock Global classical solutions to the viscous two-phase flow model with
  {N}avier-type slip boundary condition in 2{D} bounded domains.
\newblock {\em J. Math. Fluid Mech.}, 24(3):Paper No. 85, 34, 2022.

\bibitem{Lin-Liu-Zhang2005}
F.~H. Lin, C.~Liu, and P.~Zhang.
\newblock On hydrodynamics of viscoelastic fluids.
\newblock {\em Comm. Pure Appl. Math.}, 58(11):1437--1471, 2005.

\bibitem{Lions1996}
P.~L. Lions.
\newblock {\em Mathematical {T}opics in {F}luid {M}echanics. {V}ol. 1},
  volume~3 of {\em Oxford Lecture Series in Mathematics and its Applications}.
\newblock The Clarendon Press, Oxford University Press, New York, 1996.
\newblock Incompressible models, Oxford Science Publications.

\bibitem{Lions1998}
P.~L. Lions.
\newblock {\em Mathematical {T}opics in {F}luid {M}echanics. {V}ol. 2},
  volume~10 of {\em Oxford Lecture Series in Mathematics and its Applications}.
\newblock The Clarendon Press, Oxford University Press, New York, 1998.
\newblock Compressible models, Oxford Science Publications.

\bibitem{Matsumura-Nishida1980}
A.~Matsumura and T.~Nishida.
\newblock The initial value problem for the equations of motion of viscous and
  heat-conductive gases.
\newblock {\em J. Math. Kyoto Univ.}, 20(1):67--104, 1980.

\bibitem{Nirenberg1959}
L.~Nirenberg.
\newblock On elliptic partial differential equations.
\newblock {\em Ann. Scuola Norm. Sup. Pisa}, 13:115--162, 1959.

\bibitem{Oldroyd1950}
J.~G. Oldroyd.
\newblock On the formulation of rheological equations of state.
\newblock {\em Proc. Roy. Soc. London Ser. A}, 200:523--541, 1950.

\bibitem{Oldroyd1958}
J.~G. Oldroyd.
\newblock Non-{N}ewtonian effects in steady motion of some idealized
  elastico-viscous liquids.
\newblock {\em Proc. Roy. Soc. London. Ser. A}, 245:278--297, 1958.

\bibitem{Qian-Zhang2010}
J.~Z. Qian and Z.~F. Zhang.
\newblock Global well-posedness for compressible viscoelastic fluids near
  equilibrium.
\newblock {\em Arch. Ration. Mech. Anal.}, 198(3):835--868, 2010.

\bibitem{Rajagopal2015}
K.~R. Rajagopa.
\newblock Remarks on the notion of ``pressure".
\newblock {\em Int. J. Non-Linear Mech.}, 71:165--172, 2015.

\bibitem{Renardy2000}
M.~Renardy.
\newblock {\em Mathematical {A}nalysis of {V}iscoelastic {F}lows}, volume~73 of
  {\em CBMS-NSF Regional Conference Series in Applied Mathematics}.
\newblock Society for Industrial and Applied Mathematics (SIAM), Philadelphia,
  PA, 2000.

\bibitem{Renardy-Hrusa-Nohel1987}
M.~Renardy, W.~J. Hrusa, and J.~A. Nohel.
\newblock {\em Mathematical {P}roblems in {V}iscoelasticity}, volume~35 of {\em
  Pitman Monographs and Surveys in Pure and Applied Mathematics}.
\newblock Longman Scientific \& Technical, Harlow; John Wiley \& Sons, Inc.,
  New York, 1987.

\bibitem{Suen2020}
A.~Suen.
\newblock Existence and uniqueness of low-energy weak solutions to the
  compressible 3{D} magnetohydrodynamics equations.
\newblock {\em J. Differential Equations}, 268(6):2622--2671, 2020.

\bibitem{Suen-Hoff2012}
A.~Suen and D.~Hoff.
\newblock Global low-energy weak solutions of the equations of
  three-dimensional compressible magnetohydrodynamics.
\newblock {\em Arch. Ration. Mech. Anal.}, 205(1):27--58, 2012.

\bibitem{Wei-Lin-Yao2016}
R.~Y. Wei, Y.~Li, and Z.~A. Yao.
\newblock Decay of the compressible viscoelastic flows.
\newblock {\em Commun. Pure Appl. Anal.}, 15(5):1603--1624, 2016.

\bibitem{Wu-Tan2018}
G.~C. Wu and Z.~Tan.
\newblock Global low-energy weak solution and large-time behavior for the
  compressible flow of liquid crystals.
\newblock {\em J. Differential Equations}, 264(11):6603--6632, 2018.

\end{thebibliography}

\def\cprime{$'$}

\end{document}